\RequirePackage{fix-cm}
\documentclass[smallcondensed]{svjour3}     
\smartqed  
\usepackage{mathptmx}      
\usepackage{bbding}
\usepackage{amsmath,amssymb,dsfont,amstext,amsopn,graphicx,epstopdf,subfigure,url,multirow,rotating,wasysym}
\usepackage{algorithmic}
\usepackage{hyperref}
\hypersetup{
	colorlinks=true,
	linkcolor=blue,
	filecolor=blue,      
	urlcolor=blue,
	citecolor=green,
}

\newcommand{\R}{\mathbb{R}}
\newcommand{\N}{\mathbb{N}}
\newcommand{\oRn}{\vec{0}_{\R^n}}
\newcommand{\Bor}{\overline{B}(\oRn;r)}

\newcommand{\argmin}{\mathop{\rm argmin}}

\newcommand{\X}{\vec{X}}
\newcommand{\Y}{\vec{Y}}
\newcommand{\U}{\vec{U}}
\newcommand{\Ubar}{\overline{\U}}
\newcommand{\V}{\vec{V}}
\newcommand{\vl}{\vec{l}}
\newcommand{\Ueps}{\U_{\varepsilon}}
\newcommand{\Buepssigma}{B(\U_\varepsilon; \sigma)}

\newcommand{\one}{\vec{1}_{\R^n}}
\newcommand{\tolf}{\texttt{Tolf}}
\newcommand{\tolu}{\texttt{TolU}}
\newcommand{\cref}[1]{\eqref{#1}}

\graphicspath{{figs/}}

\journalname{Journal of Scientific Computing}

\begin{document}

\title{Discrete Dynamical System Approaches for Boolean Polynomial Optimization
	\thanks{In memoriam of the co-author Roland Glowinski, who passed away on January 26th 2022 during the peer-review of the manuscript. This is one of his last works. The first author was supported by the National Natural Science Foundation of China (Grant 11601327), and the second author was supported by the Hong Kong Kennedy Wong Foundation.}
}


\author{Yi-Shuai Niu         \and
        Roland Glowinski
}


\institute{Yi-Shuai Niu \at 
			  Department of Applied Mathematics, The Hong Kong Polytechnic University, Hong
			  Kong; and School of Mathematical Sciences \&  SJTU-Paristech Institute, Shanghai Jiao Tong University, Shanghai \\
              \email{yi-shuai.niu@polyu.edu.hk, niuyishuai@sjtu.edu.cn}           
           \and
           Roland Glowinski \at
              Department of Mathematics, University of Houston, 4800 Calhoun Road, Houston, TX 77204;
              and Department of Mathematics, Hong Kong Baptist University, Hong Kong \\
              \email{roland@math.uh.edu}
}

\date{Received: date / Accepted: date}

\maketitle

\begin{abstract}
In this article, we discuss the numerical solution of Boolean polynomial programs by algorithms borrowing from numerical methods for differential equations, namely the Houbolt scheme, the Lie scheme, and a Runge-Kutta scheme. We first introduce a quartic penalty functional (of Ginzburg-Landau type) to approximate the Boolean program by a continuous one and prove some convergence results as the penalty parameter $\varepsilon$ converges to $0$. We prove also that, under reasonable assumptions, the distance between local minimizers of the penalized problem and the set $\{\pm1\}^n$ is of order $O(\sqrt{n}\varepsilon)$. Next, we introduce algorithms for the numerical solution of the penalized problem, these algorithms relying on the Houbolt, Lie and Runge-Kutta schemes, classical methods for the numerical solution of ordinary or partial differential equations. We performed numerical experiments to investigate the impact of various parameters on the convergence of the algorithms. We have tested our ODE approaches and compared with the classical nonlinear optimization solver IPOPT and a quadratic binary formulation approach (QB-G) as well as an exhaustive method using parallel computing techniques. The numerical results on various datasets (including small and large-scale randomly generated synthetic datasets of general Boolean polynomial optimization problems, and a large-scale heterogeneous MQLib benchmark dataset of Max-Cut and Quadratic Unconstrained Binary Optimization (QUBO) problems) show good performances for our ODE approaches. As a result, our ODE algorithms often converge faster than the other compared methods to better integer solutions of the Boolean program.
\keywords{Boolean polynomial programs \and Houbolt scheme \and Lie scheme \and Runge-Kutta scheme}
\subclass{90C09 \and 90C10 \and 90C27 \and 90C30 \and 35Q90 \and 65K05}
\end{abstract}

\section{Introduction}
\label{sec:intro}
Integer polynomial programming has a history of more than 60 years, back to the early 1950s with the birth of combinatorial integer programming through investigations of the well-known traveling salesman problem (TSP) initialized by several pioneers: Hassler Whitney, George Dantzig, Karl Menger, Julia Robinson, Ray Fulkerson, Selmer Johnson \cite{Robinson1949,Dantzig1954,Dantzig1959} etc. Many nonlinearities in integer programming appear in form of polynomial functions (e.g., via polynomial approximations), and thus this problem has rich applications as in forms of integer linear programs: traveling salesman problem \cite{Robinson1949,Dantzig1954,Dantzig1959}, bin packing problem \cite{Man1996,De1999}, graph coloring problem  \cite{Pardalos1998}, knapsack problem \cite{Freville2004}; quadratic integer programs: capital budgeting \cite{Markowitz1952}, scheduling and allocations \cite{Moder1983}, maximum independent set problem \cite{Tarjan1977,Feo1994} and maximum cut problem \cite{Delorme1993,Goemans1995}; and general polynomial integer programs: set covering problem \cite{Chvatal1979,Caprara2000}, maximum satisfiability problem \cite{Hansen1990,Goemans1995}, vector partitioning and clustering \cite{Hansen1997}, AI and neural networks \cite{Amit1992,Personnaz1986,Glover1986,Niu2019}, and portfolio optimization with cardinality constraints and minimum transaction lots \cite{Chang2000,Li2006,Le2009,Gao2013} etc. In our paper, we are interested in minimizing a high degree (e.g., degree higher than $3$) multivariate real valued polynomial functional $P(\Y)$ with a Boolean decision variable $\Y\in \{0,1\}^n$, often referred to as \emph{Pseudo-Boolean Optimization} \cite{Lovasz1991,Boros2002}. 
This problem is obviously NP-hard in general (since its particular case - TSP is a well-known NP-complete problem, see e.g. \cite{Karp1972}).

There are three main ways to approach this problem: First, the combinatorial approaches by the development of specific algorithms for directly tackling Boolean nonlinear programs, such as branch-and-bound algorithms \cite{Dantzig1959,Efroymson1966,Kolesar1967}, cutting plane methods \cite{Gomory1958,Balas1971,Balas1993,Cornuejols2008,Niu2019dccut}, enumerative approaches \cite{Geoffrion1969}, Bender's decomposition and column generations \cite{Benders1962}, and economical linear representations for simplifying the Boolean nonlinear polynomial program (e.g., replacing cross-product terms of polynomial by additional continuous variables \cite{Balas1964,Balas1965,Watters1967,Glover1974,Glover1975}). These methods are currently most commonly used frameworks in almost all existing integer optimization solvers such as commercial solvers: BARON \cite{Sahinidis1996}, XPRESS \cite{Xpress}, LINGO/LINDO \cite{Lindo1996}, GUROBI \cite{Gurobi}, CPLEX \cite{Cplex}, MOSEK \cite{Mosek} and MATLAB intlinprog \cite{MATLAB}; and open-source solvers: SCIP \cite{Achterberg2009}, BONMIN \cite{Bonmin}, COUENNE \cite{Belotti2009}, CBC \cite{Forrest2012}, GLPK \cite{Glpk}, as well as many other COIN-OR projects \cite{Coinor}. Most of these solvers are designed for solving integer or mixed-integer linear and quadratic programs, and few of them for integer nonlinear program (e.g., BARON, BONMIN and SCIP). All of them are trying to find global optimal solutions for integer programs with moderate size. Solving very large-scale cases will be often intractable or very computationally expensive due to their inherent NP-hardness.

Second, the continuous approaches by introducing \emph{continuous reformulations} of Boolean variable as nonconvex constraints by replacing discrete variables as continuous variables with additional nonconvex constraints involving continuous functions (often polynomials). The later continuous reformulations are in general nonconvex which will be handled by classical nonlinear optimization approaches, such as: Newton-type methods \cite{Bertsekas1997}, gradient-type methods \cite{Bertsekas1997}, and difference-of-convex approaches \cite{Le2001,Niu2008,Niu2010,Pham2016,Niu2018} etc. These techniques are focused on inexpensive local optimization algorithms for finding potentially ``good" local optimal solutions, which are popularly used in many previously mentioned integer optimization solvers for local searches and bounds' estimations. Particularly, convex relaxations techniques such as linear relaxation, Lagrangian relaxation, and SDP relaxation (e.g., Lasserre's moment relaxation) \cite{Sherali1990,Lasserre2002,Parrilo2003} are often very useful either to provide good lower bounds for optimal solutions, or to construct good initializations for further local searches.

Third, the heuristic and meta-heuristic approaches such as tabu search \cite{Lokketangen1998}, genetic algorithms \cite{Gen2007}, simulated annealing \cite{Connolly1992}, ant colony \cite{Dorigo1999}, neural networks \cite{Looi1992,Smith1999}, mimic some activities in nature which does not necessarily need to be perfect, but could speed up the process of reaching a satisfactory computed solution (may not even be a local optimum) in a reasonable time frame. These methods are popular in situations where there are no known working algorithms, and often used in optimization solvers as heuristics for providing potentially good initial candidates. 

The relative effectiveness of all these approaches depends on the problem at hand, and neither approach claims a uniform advantage over the other in all cases. The reader is refereed to an excellent survey book \emph{50 Years of Integer Programming 1958-2008} \cite{Junger2009} on more histories and related topics on essential techniques for integer programming.

Clearly related to the second approach, the methodology we discuss in this paper can be summarized as follows: (i) We use \emph{penalty} to approximate the Boolean optimization problem by a continuous one. (ii) We associate with the optimality condition of the penalized problem a first or second order initial value problem (flow in Dynamical System terminology) that we time-discretize by appropriate numerical schemes (in this paper, we focused on the Houbolt and Lie schemes, and to a particular Runge-Kutta scheme; these schemes are commonly used for the solution of problems modeled by ordinary or partial differential equations). In fact, our methods are not limited for problems with polynomials, but apply also to non-polynomial differentiable cases. 

\emph{Our contributions} are mainly focused on: 
(1) Establishing a quartic penalty approximation for Boolean polynomial program, and prove that all converging sub-sequences of minimizers for penalty problems converge to a global minimizer of the Boolean problem. 
(2) Providing an error bound for the distance from solutions of the penalty problem to the set $\{\pm 1\}^n$, which is of order $O(\sqrt{n}\varepsilon)$ at best with penalty parameter $\varepsilon$. 
(3) Associating with the optimality condition of the above penalized problem a first or second order initial value problem, that we discretize by appropriate numerical schemes in order to capture its steady state solutions. The time-discretization schemes we consider are the Houbolt scheme, the Lie scheme, and a particular Runge-Kutta scheme, all classical methods for the numerical solution of ordinary or partial differential equations.
The choice of the parameters involved in these schemes is discussed, and practical strategies on choosing suitable parameter values are proposed. 
(4) A MATLAB optimization toolbox, namely \verb|DEMIPP| (Differential Equation Methods for Integer Polynomial Programming), using an effective multivariate polynomial package \verb|POLYLAB| \cite{Polylab} has been developed. Some numerical tests for solving both randomly generated high order Boolean polynomial optimization problems and MQLib large-scale heterogeneous benchmark dataset of Max-Cut and Quadratic Unconstrained Binary Optimization (QUBO) problems \cite{dunning2018} are reported. The impact of parameters to the quality of the numerical results is investigated. Our ODE approaches are compared with the classical nonlinear optimization solver IPOPT \cite{Ipopt}, a quadratic binary formulation approach (QB-G) proposed in \cite{Boros2002} combining with GUROBI solver \cite{Gurobi}, as well as an exhaustive method using parallel computing techniques. As a result, our ODE methods appear to be promising approaches which provide stable and fast convergence, especially in large-scale cases, and often obtain better integer solutions of the Boolean program.

\section{Problem Formulations}\label{sec:formulations}
Our goal in this article is to discuss the numerical solution of the following \emph{Boolean optimization} problem:
\begin{equation}
\label{prob:Boolean_opt}
\X \in \argmin_{\Y\in \{0,1\}^n} P(\Y),
\end{equation}
where in \eqref{prob:Boolean_opt}, functional $P:\R^n\to \R$ is a polynomial function of degree $d$ of the Boolean variable $\Y=(y_i)_{i=1}^n$. Note that through the article, we will use a capital letter as $\X$ for a vector and a lower case letter as $x$ for a scalar. We use the notation $\{x_i\}_i$ for a set, a family or a sequence, and the notation $(x_i)_i$ for a vector whose $i$-th coordinate is $x_i$. 

Using the transformation $\Y\to \V$ defined by
\begin{equation}
\label{eq:lintrans}
v_i = 2y_i-1, \forall i=1,\ldots,n,
\end{equation}
we can reformulate all $\{0,1\}$ variables as $\{\pm 1\}$ variables in problem \eqref{prob:Boolean_opt}, then solve
\begin{equation}
\label{prob:int_opt}
\U \in \argmin_{\V\in \{\pm 1\}^n} \Pi (\V),
\end{equation}
where the functional $\Pi$ is defined by
\[\Pi: \V \mapsto P\left(\frac{\one+\V}{2}\right) := P\left( \frac{1+v_1}{2},\ldots, \frac{1+v_n}{2} \right), ~\forall ~\V=(v_i)_{i=1}^n \in \R^n, \]
with $\one=(1)_{i=1}^n\in \R^n$, and recover $\X$ from $\U$ via
\begin{equation}
\label{eq:XfromU}
\X = (x_i)_{i=1}^n = \left( \frac{1+u_i}{2} \right)_{i=1}^n.
\end{equation}

From now on, we will consider problem \eqref{prob:int_opt} only. Actually, problem \eqref{prob:int_opt} is equivalent to 
\begin{equation}
\label{prob:int_cv_opt}
\U \in \argmin_{\V\in \{\pm 1\}^n}\left[ \frac{c}{2} \|\V\|_2^2 + \Pi(\V) \right],
\end{equation}
with $\|\V\|_2 = \sqrt{\sum_{i=1}^n v_i^2}, \forall \V\in \R^n$, and $c$ being a positive constant (indeed, since $\|\V\|_2^2 = n$, problems \eqref{prob:int_cv_opt} and \eqref{prob:int_opt} are equivalent).

\begin{remark}\label{rmk:1}
	Concerning $c$, we suggest taking it not too large but large enough, so that the functional $\V\to \frac{c}{2}\|\V\|_2^2 + \Pi(\V)$ is convex over the ball of radius $r (> \sqrt n)$ centered at $\oRn$. In the next section, we will introduce a quartic penalty approximation of problem \eqref{prob:int_cv_opt}. This penalization may cancel convexity, but the regularization term $\frac{c}{2}\|\V\|_2^2$ may help solving the nonconvex problem and induces small norm solutions (when the problem is relaxed into
	a continuous one). Details on the choice of $c$ and its impact on the numerical results will be given in Subsections \ref{subsec:choiceofc} and \ref{subsec:testsonrandomdataset}.
\end{remark}

\section{A Quartic Penalty Approximation of Problem (\ref{prob:int_cv_opt})}\label{sec:penalty_approximation}

Let $\varepsilon$ be a positive parameter, let $\U\circ \V$ be the Hadamard product of the vector $\U$ and $\V$ \footnote{$\forall~ \U=(u_i)_{i=1}^n\in \R^n, \forall~ \V=(v_i)_{i=1}^n\in \R^n, \U\circ \V = (u_i v_i)_{i=1}^n\in \R^n.$}. We approximate problem \cref{prob:int_cv_opt} by penalizing the constraint $\V\in\{\pm 1\}^n$ using the quartic penalty functional $$\V\mapsto \frac{1}{4\varepsilon} \|\V\circ \V - \one\|_2^2 = \frac{1}{4\varepsilon} \sum_{i=1}^n(v_i^2-1)^2.$$ The resulting penalized problem reads as
\begin{equation}
\label{prob:approx_int_cv_opt}
\Ueps \in \argmin_{\V\in \R^n} ~ J_{\varepsilon}(\V).
\end{equation}
where $$J_{\varepsilon}(\V) = \frac{1}{4\varepsilon} \|\V\circ \V -\one \|_2^2  +\frac{c}{2} \|\V\|_2^2 + \Pi(\V).$$ 
Let us denote by $d$ the degree of the polynomial functional $\Pi$. We can easily show that as long as $d\leq 4$, then we have
\begin{equation}
\label{eq:limJesp}
\lim_{\|\V\|_2\to +\infty} J_{\varepsilon} (\V) = +\infty,
\end{equation}
if $\varepsilon$ is sufficiently small. Relation \cref{eq:limJesp} implies that problem \cref{prob:approx_int_cv_opt} has a solution, possibly non-unique. Actually, if $\Pi$ is convex, relation \cref{eq:limJesp} holds $\forall \varepsilon>0$ and $c\geq 0$, implying that the associated problem \cref{prob:approx_int_cv_opt} has a solution. However, when $d>4$, relation \cref{eq:limJesp} does not hold in general. In this case, as $J_{\varepsilon}(\V)$ is a polynomial, then (except for constant case which never occurs when $d> 4$) $J_{\varepsilon}(\vec{V})$ goes to $+\infty$ or $-\infty$ as $\|\vec{V}\|\to +\infty$. Hence, if the relation \cref{eq:limJesp} does not hold, then problem \cref{prob:approx_int_cv_opt} has no solution. To overcome this difficulty, we observe that any solution to problem \cref{prob:int_cv_opt} belongs to the closed ball $\Bor$ of radius $r \geq \sqrt n$ centered at $\oRn$. This observation suggests approximating \cref{prob:int_cv_opt} by the following constrained variant of problem \cref{prob:approx_int_cv_opt}:
\begin{equation}
\label{prob:penaltyoverball}
\U_{\varepsilon} \in \argmin_{V\in \Bor} J_{\varepsilon}(\V).
\end{equation}
with $r> \sqrt{n}$ (one can take $r=\sqrt{n}$, but taking $r>\sqrt{n}$ makes things simpler mathematically). 
\begin{remark}\label{rmk:2}
	Problems \cref{prob:approx_int_cv_opt} and \cref{prob:penaltyoverball} are not equivalent in general. However, if problem \cref{prob:approx_int_cv_opt} has a solution, then this solution is also solution to problem \cref{prob:penaltyoverball} for $r$ large enough. Note that the optimal solution $\U_{\varepsilon}$ in both problem \cref{prob:approx_int_cv_opt} and problem \cref{prob:penaltyoverball} may be non-unique.
\end{remark}

Based on the classical convergence result for penalization methods in finite dimensional optimization, we can prove that problem \cref{prob:penaltyoverball} has a solution, and that we can extract from the family $\{\U_{\varepsilon} \}_{\varepsilon >0}$ (with $\varepsilon$ decreasing to $0^+$) a sub-sequence converging to a solution of problem \cref{prob:int_cv_opt}. Although, the above results are classical, we decided to include their proof in this article since we will use related techniques in Section \ref{sec:errordiscussion} to estimate the distance between exact and approximate solutions as a function of $\varepsilon$ (and other parameters in $J_{\varepsilon}$). We have then the following theorem:
\begin{theorem}
	\label{thm:penalization}
	Let $\varepsilon>0$ and $r>\sqrt{n}$. Then the penalized problem \cref{prob:penaltyoverball} has a solution $\Ueps$ such that one can extract from the family $\{\U_{\varepsilon}\}_{\varepsilon>0}$ (with $\varepsilon$ decreasing to $0$) a sub-sequence converging to a solution of problem \cref{prob:int_cv_opt} (a global minimizer of functional $\Pi$ over $\{\pm 1\}^n$). 
\end{theorem}
\begin{proof}
	$\rhd$ \textit{Existence of a solution to problem \cref{prob:penaltyoverball}:} The set $\Bor$ being compact and non-empty, and the functional $\Pi$ being continuous over $\Bor$, it follows from the Weierstra\ss\ extreme value theorem that problem \cref{prob:penaltyoverball} has a solution.
	
	\noindent $\rhd$ \textit{Existence of sub-sequences extracted from $\{\U_{\varepsilon}\}_{\varepsilon>0}$ converging to solutions of problems \cref{prob:int_opt} and \cref{prob:int_cv_opt}:} The set $\Bor$ being compact and non-empty, it follows from the theorem of Bolzano-Weierstra\ss\ that one can extract from $\{\Ueps\}_{\varepsilon>0}$ a converging sub-sequence (that we still denote by $\{\Ueps\}_{\varepsilon>0}$ for simplicity) such that  
	\begin{equation}
	\label{eq:limUeps}
	\lim_{\varepsilon\to 0^+} \Ueps = \U \in \Bor.
	\end{equation}
	The functions $\|\cdot\|_2$ and $\Pi$ being continuous over $\R^n$, it follows from \cref{eq:limUeps} that
	\begin{equation}
	\label{eq:limits}
	\left\lbrace 
	\begin{aligned}
	& \lim_{\varepsilon \to 0^+} \|\Ueps\|_2 = \|\U\|_2,\\
	& \lim_{\varepsilon \to 0^+} \Pi(\Ueps) = \Pi(\U).
	\end{aligned}
	\right.
	\end{equation}
	Consider now a solution $\U^*$ to problem \cref{prob:int_cv_opt}. From the equivalence between \cref{prob:int_cv_opt} and \cref{prob:int_opt}, $\U^*$ is also a minimizer of functional $\Pi$ over $\{\pm 1\}^n \subset \Bor$.  We have then 
	\begin{equation}
	\label{eq:ineqJeps}
	\begin{aligned}
	\frac{1}{4\varepsilon} \|\Ueps \circ \Ueps - \one\|_2^2 +\frac{c}{2} \|\Ueps\|_2^2 + \Pi(\Ueps)\leq \\
	\frac{1}{4\varepsilon} \|\U^* \circ \U^* - \one\|_2^2 +\frac{c}{2} \|\U^*\|_2^2 + \Pi(\U^*).
	\end{aligned}
	\end{equation}
	Since $\U^*\in \{\pm 1\}^n$, we have $\|\U^*\circ \U^*-\one\|_2=0$ and $\|\U^*\|_2^2 = n$, implying that \cref{eq:ineqJeps} reduces to 
	\begin{equation}
	\label{eq:bndforueps^2-1}
	\|\Ueps \circ \Ueps - \one\|_2^2  \leq 4\varepsilon \left[ \frac{cn}{2} + \Pi(\U^*)-\Pi(\Ueps) \right].
	\end{equation}
	Combining \cref{eq:limUeps}, \cref{eq:limits} and \cref{eq:bndforueps^2-1}, we obtain
	\[ \lim_{\varepsilon\to 0^+} \|\Ueps \circ \Ueps - \one\|_2^2 = \|\U \circ \U - \one\|_2^2  = 0,\]
	which implies 
	\begin{equation}
	\label{eq:Uinpm1}
	\U\in \{\pm 1\}^n.
	\end{equation}
	It follows from \cref{eq:ineqJeps} that 
	\begin{equation}
	\label{eq:ineq2}
	\frac{c}{2} \|\Ueps\|_2^2 + \Pi(\Ueps) \leq  \frac{cn}{2} + \Pi(\U^*).
	\end{equation}
	Combining \cref{eq:limUeps}, \cref{eq:limits}, \cref{eq:Uinpm1} and \cref{eq:ineq2}, one obtains
	\begin{equation}
	\label{eq:convresult}
	\left\lbrace 
	\begin{aligned}
	& \U\in \{\pm 1\}^n,\\
	& \Pi(\U) \leq \Pi(\U^*) \leq \Pi(\V), ~\forall \V\in \{\pm 1\}^n.
	\end{aligned}
	\right.
	\end{equation}
	It follows from \cref{eq:convresult} that $\U$ is a global minimizer of $\Pi$ over $\{\pm 1\}^n$, that is a solution of problems \cref{prob:int_opt} and \cref{prob:int_cv_opt}.
	\qed
\end{proof}

\begin{remark}\label{rmk:3}
	Quadratic penalty methods are known to be inexact (we mean by inexact that, in general, the solution(s) of the associated penalized problem do not verify the constraints one has penalized). The quartic penalties used in this article as in problem \cref{prob:approx_int_cv_opt} and \cref{prob:penaltyoverball} are not exceptions, as we shall see below. Indeed, a very simple and convincing example is the following one:
	\begin{example}(Inexact penalty) 
		Suppose that the penalized problem \cref{prob:penaltyoverball} (with $r>\sqrt{n}$) has a solution $\Ueps\in \{\pm 1\}^n$. Then, the pair $(c,\Ueps)$ verifies the KKT optimality conditions for problem \cref{prob:penaltyoverball} as:
		\begin{equation}
		\label{eq:kkt}
		\left\lbrace 
		\begin{aligned}
		& c>0, \Ueps \in \{\pm 1\}^n, \lambda \geq 0,\\
		& \frac{1}{\varepsilon} (\Ueps\circ \Ueps - \one)\circ \Ueps + c \Ueps + \nabla \Pi(\Ueps) + 2\lambda \Ueps = \oRn,\\
		& \lambda (r^2 - \|\Ueps\|_2^2) = 0.
		\end{aligned}
		\right.
		\end{equation}
		Since $r>\sqrt{n}$, the complementarity condition $\lambda (r^2 - \|\Ueps\|_2^2) = 0$ implies that $\lambda = 0$, then system \cref{eq:kkt} reduces to 
		\begin{equation}
		\label{eq:optcond2}
		\left\lbrace 
		\begin{aligned}
		&c >0, \Ueps \in \{\pm 1\}^n,\\
		&c\Ueps + \nabla \Pi(\Ueps) = \oRn,
		\end{aligned}
		\right.
		\end{equation}
		which is (a kind of) generalized (nonlinear if $d\geq 3$) eigenvalue problem. System \cref{eq:optcond2} has no solution in general, implying that $\Ueps \notin\{\pm 1\}^n$, contradicting the initial assumption. For example, suppose that we choose 
		$$c > \frac{\max_{\V\in \Bor }\|\nabla \Pi (\V)\|_2}{\sqrt{n}},$$
		then the related system \cref{eq:optcond2} has no solution, implying inexact penalty.\qed
	\end{example}
\end{remark}

In the following sections we are going to investigate the solution of the unconstrained penalized problem, the main reason, being (besides simplicity) that the numerical experiments we
performed, with random initializations, never encountered any trouble associated with the
possible non-verification of the property \cref{eq:limJesp}, that is
$$\lim_{\|\V\|_2\to +\infty} J_{\varepsilon} (\V) = +\infty.$$

The optimality condition associated with the penalized problem \cref{prob:approx_int_cv_opt} reads as:
\begin{equation}
\label{eq:optsys1}
\frac{1}{\varepsilon} (\Ueps\circ \Ueps- \one) \circ \Ueps + c\Ueps + \nabla \Pi(\Ueps) = \oRn.	
\end{equation}
There is no doubt that many methods are applicable to the solution of problem  \cref{eq:optsys1}. In this article, we took the following (classical) approaches: we associated with \cref{eq:optsys1} a first order or second order accurate time-stepping method such as the Houbolt scheme, the Lie scheme, and a MATLAB explicit Runge-Kutta ODE solver.

Among the three schemes we are going to employ, the Lie scheme is, as shown in Subsection \ref{subsec:Lie}, the one that can handle more easily the constraint $\|\V\|_2\leq r$ encountered in the penalized problem \cref{prob:penaltyoverball}.

\begin{remark}
	Problem \cref{eq:optsys1} being a system of $n$ nonlinear equations, it makes sense trying solving it by existing solvers for such systems, before moving to more dedicated solvers. Actually, it is what we did by applying MATLAB \verb|fsolve| to the solution of problem \cref{eq:optsys1} (\verb|fsolve| relies on a trust-region dogleg algorithm, a variant of the Powell dogleg method). However, numerical experiments showed that the quality of the computed solutions is rather poor, since (using random initializations) only $15\% \sim 20\%$ of the computed solutions were close to elements of the set $\{\pm 1\}^n$. These poor results drove us to look for the alternative approaches discussed in the following sections.
\end{remark}

\section{ODE Approaches for The Solution of System \eqref{eq:optsys1}} \label{sec:ODEapproaches}
In this section, we are going to investigate ODE approaches (namely, the Lie scheme and the Houbolt scheme) for the solution of system \eqref{prob:1storderode}. Basically, they are based on the numerical solution of the first order or the second order ODE associated with system \eqref{eq:optsys1}.

\subsection{A Lie Scheme} \label{subsec:Lie}
The first order in time ordinary differential equation associated with \cref{eq:optsys1} reads:
\begin{equation}\label{eq:prob_odeforLie}
\left\lbrace 
\begin{aligned}
&\U(t)=(u_i(t))_{i=1}^n,\\
&\left\lbrace 
\begin{aligned}
&\dot{u}_i(t) + \frac{1}{\varepsilon} (u_i^2(t)-1)u_i(t) + cu_i(t) + \frac{\partial \Pi}{\partial v_i}(\U(t)) = 0, \forall t>0,\\
&i=1,\ldots, n,\\
\end{aligned}
\right.\\
&\U(0)=\U_0.
\end{aligned}
\right.
\end{equation}
The initial value problem \cref{eq:prob_odeforLie} (a gradient flow) can be written also as
\begin{equation}\label{eq:splitted}
\left\lbrace 
\begin{aligned}
&\U(t)=(u_i(t))_{i=1}^n,\\
&\dot{\U}(t) + A_1(\U(t)) + A_2(\U(t)) = 0, \forall t>0,\\
&\U(0)=\U_0.
\end{aligned}
\right.
\end{equation}
with operators $A_1$ and $A_2$ defined by
\begin{equation}
\label{eq:splitting-operators}
\left\lbrace
\begin{aligned}
&A_1(\V)= \nabla \Pi (\V), ~\forall \V\in \R^n,\\
&A_2(\V) = \frac{1}{\varepsilon}(\V\circ \V-\one)\circ \V + c\V, ~\forall \V\in \R^n.
\end{aligned}
\right.
\end{equation}
The structure of \cref{eq:splitted} suggests using operator-splitting for its time-integration from $t=0$ to $t=+\infty$. The simplest operator splitting scheme we can think about for the time-integration of problem \cref{eq:splitted} is clearly the following variant of the \emph{Lie scheme} (known as the \emph{Marchuk-Yanenko scheme}, very popular in Computational Mechanics and Physics for its simplicity and robustness (see, e.g., \cite{Glowinski2017} and the references therein for details and applications)). The \textbf{Marchuk-Yanenko scheme} we are going to use reads as follows (with as usual $\tau$ ($>0$) a time-discretization step and $\U^k$ an approximation of $\U(k\tau)$): 

\noindent\textbf{1) Initialization:}
\begin{equation}
\label{eq:init_Lie}
\boxed{\U^0=\U_0.}
\end{equation}
\noindent\textbf{2) Computations of $\U^k\to \U^{k+\frac{1}{2}}\to \U^{k+1}$ for $k\geq 0$:} 
Solve
\begin{equation}
\label{eq:prob_uk0.5_Lie}
\frac{\U^{k+\frac{1}{2}} - \U^k}{\tau} + A_1(\U^{k+\frac{1}{2}}) = \oRn.
\end{equation}
and 
\begin{equation}
\label{eq:prob_uk+1_Lie}
\frac{\U^{k+1} - \U^{k+\frac{1}{2}}}{\tau} + A_2(\U^{k+1}) = \oRn.
\end{equation}
One can rewrite \cref{eq:prob_uk0.5_Lie} using \cref{eq:splitting-operators} as a nonlinear system:
\begin{equation}
\label{eq:system_uk0.5_lie}
\boxed{\U^{k+\frac{1}{2}} + \tau \nabla \Pi \left(\U^{k+\frac{1}{2}}\right) = \U^k,}
\end{equation}
which can be solved by employing, for example, the MATLAB function \verb|fsolve|. Again, from \cref{eq:splitting-operators}, one can write problem \cref{eq:prob_uk+1_Lie} as
\begin{equation}\label{eq:system_uk+1_lie}
\left\lbrace 
\begin{aligned}
&\U^{k+1}=(u_i^{k+1})_{i=1}^n,\\
&\left\lbrace \begin{aligned}
& \frac{\tau}{\varepsilon} \left(u_i^{k+1}\right)^3 + \left(1+c\tau - \frac{\tau}{\varepsilon}\right) u_i^{k+1} = u_i^{k+\frac{1}{2}},\\
&i=1,\ldots,n,
\end{aligned}
\right.
\end{aligned}
\right.
\end{equation}
a system of $n$ uncoupled cubic equations, which has a unique solution $\U^{k+1}$ if 
\begin{equation}
\label{eq:convcond_Lie}
c + \frac{1}{\tau} \geq \frac{1}{\varepsilon}.
\end{equation}
For $c>0$, condition \cref{eq:convcond_Lie} is verified if 
\begin{equation}
\label{eq:boundsoftauandeps}
\varepsilon \leq \frac{1}{c} \text{ and } 0< \tau \leq \frac{\varepsilon}{1-\varepsilon c}.
\end{equation}
For $c=0$, condition \cref{eq:convcond_Lie} is simplified as 
\begin{equation}
\label{eq:boundsoftauforc=0}
0< \tau\leq \varepsilon.
\end{equation}	
Assuming that \cref{eq:convcond_Lie} holds, one can solve the $n$ equations in \cref{eq:system_uk+1_lie} by using Newton’s method initialized with $\U^{k+\frac{1}{2}}$.
	Moreover, it is also possible to solve the system \cref{eq:system_uk+1_lie} by using the well-known Cardano's formula, then the unique real root is given by
	\begin{equation}
	\label{eq:uk+1lie_cardano}
	\boxed{u_i^{k+1} = \sqrt[3]{-\frac{q_i}{2} + \sqrt{\frac{q_i^2}{4} + \frac{p^3}{27}}} + \sqrt[3]{-\frac{q_i}{2} - \sqrt{\frac{q_i^2}{4} + \frac{p^3}{27}}}, ~i=1,\ldots,n,}
	\end{equation}
	where $$p = \frac{\varepsilon}{\tau} + c \varepsilon -1;~ q_i = -\frac{\varepsilon}{\tau} u_i^{k+\frac{1}{2}}.$$
\begin{remark}
		Note that the cubic roots are also computed by the Newton's methods, despite the explicit formulation for the solution of $\U^{k+1}$ by Cardano's formula, these two approaches for solving \cref{eq:system_uk+1_lie} should be as fast as each other. 
\end{remark}
\begin{remark}
	Rewriting $J_{\varepsilon}(\V)$ as $J_1(\V) + \Pi(\V)$ where $J_1(\V)=\frac{1}{4\varepsilon} \|\V\circ \V -\one \|_2^2  +\frac{c}{2} \|\V\|_2^2$. Then the discretizations in \eqref{eq:prob_uk0.5_Lie} and \eqref{eq:prob_uk+1_Lie} are nothing but two steps of the proximal point method as
	$$\U^{k+\frac{1}{2}}=\text{prox}_{\tau \Pi}(\U^k),\qquad  \U^{k+1} = \text{prox}_{\tau J_1}(\U^{k+\frac{1}{2}})$$
	where prox is the standard proximal operator defined by
	$$\text{prox}_f(x) = \argmin_{y}\{f(y) + \frac{1}{2}\|y-x\|_2^2\}.$$
	The second proximal operator has a unique closed-form solution as in \eqref{eq:uk+1lie_cardano}, while the first proximal operator involves solving a polynomial equation \eqref{eq:system_uk0.5_lie}.
\end{remark}

\subsection{A Houbolt Scheme}\label{subsec:houbolt}
The algorithm we are going to investigate in this subsection is related to the \emph{B.T. Polyak Heavy Ball Method} \cite{polyak1964}. Actually, as shown by Su, Boyd and Candes \cite{Su2015}, the celebrated \emph{Nesterov minimization algorithm} belongs to the heavy ball family. The Nesterov algorithm reads as
\begin{equation}
\label{eq:nesterovalgo}
\left\lbrace 
\begin{aligned}
&\Y_0=\X_0\in \R^n,\\
&\left\lbrace 
\begin{aligned}
&\Y_k = \X_k + \beta_k (\X_k - \X_{k-1}),\\
&\X_{k+1} = \Y_k - s \nabla f(\Y_k), \forall k\in \N^*,
\end{aligned}
\right.
\end{aligned}
\right.
\end{equation}
where the step size parameter $s>0$ (often, one takes $s\leq 1/L$) and $\beta_k=\frac{k-1}{k+2}>0$ and the function $f:\R^n\to \R$ is $L$-smooth, i.e., $$\|\nabla f(\X) - \nabla f(\Y)\|_2 \leq L\|\X-\Y\|_2,~ \forall \X, \Y\in \R^n.$$ In addition, if $f$ is convex, then this scheme enjoys inverse quadratic convergence rate:
	$$f(\X_k) - f^* \leq O\left( \frac{\|\X_0 - \X^*\|}{s k^2} \right),$$
where $\X^*$ is any minimizer of $f$ and $f^* = f(\X^*)$.

According to \cite{Su2015}, the \emph{first-order} Nesterov algorithm \cref{eq:nesterovalgo} is nothing but a discrete form of the following \emph{second-order} ODE:
\begin{equation}
\label{eq:nesterovode}
\left\lbrace 
\begin{aligned}
&\ddot{\X}(t) + a(t)\dot{\X}(t)+\nabla f(\X(t)) = \oRn, \forall t>0,\\
&\X(0)=\X_0, \dot{\X}(0)=\oRn,
\end{aligned}
\right.
\end{equation}
with $a(t)=\frac{3}{t}$. The differential system \cref{eq:nesterovode} is definitely of the Polyak Heavy Ball type, with a vanishing damping term $a(t)$ as $t\to +\infty$. 

The second-order in time differential system we associate with \cref{eq:optsys1} reads as:
\begin{equation}
\label{eq:prob_2ndordeode}
\left\lbrace 
\begin{aligned}
&\U(t)=(u_i(t))_{i=1}^n,\\
&\left\lbrace 
\begin{aligned}
&m\ddot{u}_i(t) + \gamma \dot{u}_i(t) + \frac{1}{\varepsilon} (u_i^2(t)-1)u_i(t) + cu_i(t) + \frac{\partial \Pi}{\partial v_i}(\U(t)) = 0, \forall t>0,\\
&i=1,\ldots, n,\\
\end{aligned}
\right.\\
&\U(0)=\U_0, \dot{\U}(0) = \V_0,
\end{aligned}
\right.
\end{equation}
with $m$ and $\gamma$ two positive constants.

Many methods are applicable to the numerical solution of system \cref{eq:prob_2ndordeode}, the most obvious one is to introduce $v_i = \dot{u}_i, \forall i=1,\ldots,n$, and solve the resulting equivalent first order system, namely
\begin{equation}
\label{prob:1storderode}
\left\lbrace 
\begin{aligned}
&(\U(t),\V(t))=(u_i(t),v_i(t))_{i=1}^n,\\
&\left\lbrace 
\begin{aligned}
&m\dot{v}_i(t) + \gamma v_i(t) + \frac{1}{\varepsilon} (u_i^2(t)-1)u_i(t) + cu_i(t) + \frac{\partial \Pi}{\partial v_i}(\U(t)) = 0, \forall t>0,\\
&\dot{u}_i(t) = v_i(t), \forall t>0,\\
&i=1,\ldots, n,\\
\end{aligned}
\right.\\
&\U(0)=\U_0, \V(0) = \V_0,
\end{aligned}
\right.
\end{equation}
using one of these user friendly ODE solvers, such as the Runge-Kutta scheme, from MATLAB or other scientific computing packages. The formulations of Runge-Kutta schemes are classical, thus will not be described in this paper. We will test a MATLAB implementation (namely \texttt{ode45} \cite{shampine1997matlab}, an explicit Runge-Kutta (4,5) formula, the Dormand-Prince pair) in Section \ref{sec:simulations} of numerical simulations. For more details about this Runge-Kutta scheme, the reader may refer to \cite{dormand1980family}.

The scheme discussed below is a semi-implicit finite difference time-discretization scheme of the Houbolt type, a classical and popular scheme in \emph{structural dynamics} (we used it, coupled to ADMM, to simulate the vibrations of nonlinear elastic beams (see \cite{Bourgat1980} for details)). Its main drawback is that it requires a starting procedure, but this is not a difficult issue to overcome. The \textbf{Houbolt scheme} reads as follows:

\noindent\textbf{1) Initialization:} 
\begin{equation}
\label{eq:init}
\boxed{\U^0=\U_0, (\U^1 - \U^{-1}) = 2\tau \V_0,}
\end{equation}
where the second relation is derived from the Taylor expansion at $t=0$:
\begin{equation}
\label{eq:taylor1st}
\frac{\U(t+\tau) - \U(t-\tau)}{2\tau} = \dot{\U}(t) + O(\tau^2),
\end{equation}
\noindent\textbf{2) Computations of $\U^1$ and $\U^{-1}$:}
At $t=0$, taking advantage of \cref{eq:taylor1st} and
\begin{equation}
\label{eq:taylor2nd}
\frac{\U(t+\tau) + \U(t-\tau) -2 \U(t)}{\tau^2} = \ddot{\U}(t) + O(\tau^2),
\end{equation}
one time-discretizes \cref{eq:prob_2ndordeode} by 
\begin{equation}
\label{eq:prob_u1}
\left\lbrace 
\begin{aligned}
&\U^1 =(u_i^1)_{i=1}^n\in \R^n,\\
&\left\lbrace 
\begin{aligned}
&m\frac{u_i^1 + u_i^{-1} - 2u_i^0}{\tau^2} + \gamma \frac{u_i^1 - u_i^{-1}}{2\tau} + \frac{1}{\varepsilon} ((u_i^0)^2-1)u_i^0 + cu_i^0 + \frac{\partial \Pi}{\partial v_i}(\U_0) = 0,\\
&i=1,\ldots, n.\\
\end{aligned}
\right.
\end{aligned}
\right.
\end{equation}
Using \cref{eq:init}, we can eliminate $\U^{-1}$ in \cref{eq:prob_u1}, obtaining thus 
\begin{equation}
\label{eq:prob_u1_simple}
\left\lbrace 
\begin{aligned}
&\U^1 =(u_i^1)_{i=1}^n\in \R^n,\\
&\left\lbrace 
\begin{aligned}
&2m\frac{u_i^1 - u_i^0 - \tau v_i^0}{\tau^2} + \gamma v_i^0 + \frac{1}{\varepsilon} ((u_i^0)^2-1)u_i^0 + cu_i^0 + \frac{\partial \Pi}{\partial v_i}(\U_0) = 0,\\
&i=1,\ldots, n,\\
\end{aligned}
\right.
\end{aligned}
\right.
\end{equation}
which implies in turn 
\begin{equation}
\label{eq:u1}
\boxed{\U^1 = \left(\tau - \frac{\tau^2 \gamma}{2m}\right)\V_0 + \left(1+\frac{\tau^2 }{2m}\left(\frac{1}{\varepsilon}-c\right)\right)\U_0 - \frac{\tau^2}{2m\varepsilon} \U_0\circ \U_0\circ \U_0 - \frac{\tau^2}{2m}\nabla \Pi(\U_0).}
\end{equation}
Once $\U^1$ is known, it follows from \cref{eq:init} that
\begin{equation}
\label{eq:um1}
\boxed{\U^{-1} = \U^1 - 2\tau \V_0.}
\end{equation}
\noindent\textbf{3) Computations of $\U^{k+1}$ for $k\geq 1$:}
Assuming that $\U^k$,$\U^{k-1}$ and $\U^{k-2}$ are known, we obtain $\U^{k+1}$ using the following second order accurate semi-implicit scheme of the Houbolt type:
\begin{equation}
\label{eq:prob_uk+1}
\left\lbrace 
\begin{aligned}
&\U^{k+1} =(u_i^{k+1})_{i=1}^n,\\
&\left\lbrace
\begin{aligned}
&m\frac{2u_i^{k+1} - 5u_i^{k} + 4 u_i^{k-1} - u_i^{k-2}}{\tau^2} + \gamma \frac{3u_i^{k+1}-4u_i^k + u_i^{k-1}}{2\tau} + \\
&\frac{1}{\varepsilon} ((u_i^{k+1})^2-1)u_i^{k+1} + c(2u_i^{k}-u_i^{k-1}) + \frac{\partial \Pi}{\partial v_i}(2\U^k-\U^{k-1}) = 0,\\
&i=1,\ldots, n.\\
\end{aligned}
\right. 
\end{aligned}
\right. 
\end{equation}
System \cref{eq:prob_uk+1} has been obtained by discretizing \cref{eq:prob_2ndordeode} at $t=(k+1)\tau$, taking the following relations into account:
\begin{equation}\label{eq:taylor2ndbis}
\left\lbrace 
\begin{aligned}
&\frac{2u_i(t)-5u_i(t-\tau) + 4u_i(t-2\tau)-u_i(t-3\tau)}{\tau^2} = \ddot{u}_i(t) + O(\tau^2),\\
&\frac{3u_i(t)-4u_i(t-\tau) + u_i(t-2\tau)}{2\tau} = \dot{u}_i(t) + O(\tau^2),\\
&2u_i(t-\tau) - u_i(t-2\tau) = u_i(t) + O(\tau^2),\\
& i=1,\ldots,n.
\end{aligned}
\right. 
\end{equation}
obtained by Taylor's expansions at $t$.
It follows from \cref{eq:prob_uk+1} that, $\forall i=1,\ldots,n$, the term $u_i^{k+1}$ is solution of a cubic equation which has a unique solution if condition
\begin{equation}
\label{eq:convcond}
\frac{2m}{\tau^2} + \frac{3\gamma}{2\tau}\geq \frac{1}{\varepsilon}
\end{equation}
holds. If we assume that
$$0 < \tau \leq \frac{3\gamma \varepsilon + \sqrt{9\gamma^2\varepsilon^2 + 32m\varepsilon}}{4},$$
then condition \cref{eq:convcond} will be automatically verified. Or we can simply take 
\begin{equation}
\label{eq:boundsoftau}
0 < \tau \leq \sqrt{2m\varepsilon}.
\end{equation}
From now on, we will assume that condition \cref{eq:boundsoftau} holds. In this case, the unique solution of system \eqref{eq:prob_uk+1} can be computed again using Newton's methods, or using Cardano's formula as
	\begin{equation}
	\label{eq:uk+1_houbolt_cardano}
	\boxed{u_i^{k+1} = \sqrt[3]{-\frac{q_i}{2} + \sqrt{\frac{q_i^2}{4} + \frac{p^3}{27}}} + \sqrt[3]{-\frac{q_i}{2} - \sqrt{\frac{q_i^2}{4} + \frac{p^3}{27}}}, ~i=1,\ldots,n,}
	\end{equation}
	where $$\left\{\begin{aligned}p =& \left( \frac{2m}{\tau} + \frac{3\gamma}{2} \right)\frac{\varepsilon}{\tau} - 1;\\	
	q_i =& \frac{m\varepsilon}{\tau^2}\left(- 5u_i^{k} + 4 u_i^{k-1} - u_i^{k-2}\right) +  \frac{\gamma\varepsilon}{2\tau}\left(-4u_i^k + u_i^{k-1}\right) \\
	&+ c\varepsilon(2u_i^{k}-u_i^{k-1}) + \varepsilon\frac{\partial \Pi}{\partial v_i}(2\U^k-\U^{k-1}).
	\end{aligned}\right.$$
	\begin{remark}
		Note that the first order ODE system \eqref{eq:prob_odeforLie} and the second order ODE system \eqref{eq:prob_2ndordeode} are in fact discretization of the well-known gradient flow \begin{equation}
			\label{eq:firstorderode}
			\dot{U}=-\nabla J_{\varepsilon}(U)
		\end{equation} and its accelerated version \begin{equation}
		\label{eq:secondorderode}
		\ddot{U} + \gamma(t)\dot{U}=-\nabla J_{\varepsilon}(U).
	\end{equation} Actually, pretty much all
		first-order optimization methods (with their respective accelerated variants) can be seen as
		discretizations of these systems as well (see \cite{francca2021gradient} suggested by a referee). For example, ADMM can be seen as a discretization of \eqref{eq:firstorderode} and an accelerated ADMM can be seen as a discretization of  \eqref{eq:secondorderode} (see \cite{franca2018admm}).
	\end{remark}
\subsection{Stopping Criteria}\label{subsec:stoppingcond}
	There are several possibilities to terminate our proposed ODE schemes (i.e., the Houbolt, Lie and Runge-Kutta schemes) as :
	\begin{description}
		\item[1)] Terminate algorithms at $\U^{k}$ if $|\Pi(\U^{k})-\Pi(\U^{k-1})|\leq \texttt{Tolf}$ with given \texttt{Tolf}$>0$. 
		\item[2)] Terminate algorithms at $\U^{k}$ if $\|\U^{k}-\U^{k-1}\|\leq \texttt{TolU}$ with given \texttt{TolU}$>0$.
\end{description}
These two criteria are natural ways to terminate many iterative methods, and our ODE methods could be terminated when one of the two criteria is verified. The first one implies that the sequence $\{\Pi(\U^{k})\}_k$ converges in a desired precision \texttt{Tolf}, and the second one indicates that the sequence $\{ \U^k \}_k$ converges in a desired precision \texttt{TolU}. The convergences of theses sequences are guaranteed by Theorem \ref{thm:penalization}. Note that based on Theorem \ref{thm:errorbound}, we can also terminate algorithms when $\delta$ is smaller enough.

\section{Estimating The Distance Between The Minimizers of Problem \cref{prob:penaltyoverball} And $\{\pm 1\}^n$}
\label{sec:errordiscussion}

One important question is how close (with respect to $\varepsilon$) are the local minimizers of the penalized problem \cref{prob:penaltyoverball}, to the closest element(s) of $\{\pm 1\}^n$. We will prove below, that penalization leads to approximate solutions whose distance to $\{\pm 1\}^n$ is $O(\sqrt{n}\varepsilon)$ if some reasonable conditions are verified.
\begin{theorem}\label{thm:errorbound}
	The positive parameter $\varepsilon$ being fixed. Let us denote by $J_{\varepsilon}$ the functional 
	\[\V\mapsto \frac{1}{4\varepsilon} \|\V \circ \V - \one \|_2^2 + \frac{c}{2} \|\V\|_2^2 + \Pi (\V).\] Let $\Ueps$ be a minimizer of functional $J_{\varepsilon}$ over the set $\Bor$ ($r>\sqrt{n}$), and denote by $\delta$ the Euclidean distance from $\Ueps$ to the set $\{\pm 1\}^n$. 
	Then we have
	\begin{equation}
		\label{eq:errorbound}
		\delta \leq \frac{4\varepsilon}{1+2c\varepsilon} \left( c\sqrt{n} + \max_{\V\in \Bor} \|\nabla \Pi(\V)\|_2 \right).
	\end{equation}
\end{theorem}
\begin{proof} First, we will show that there exists an open ball $\Buepssigma$ such that  
		\begin{equation}
			\label{eq:localminimizer}
			J_{\varepsilon}(\Ueps) \leq J_{\varepsilon}(\V), ~\forall \V\in \Buepssigma \cap \Bor,
		\end{equation}
		and the set $\Buepssigma \cap \{\pm 1\}^n$ is non-empty. In fact, by the definition of $U_{\varepsilon}$ (a minimizer of $J_{\varepsilon}$ over $\Bor$, see equation \eqref{prob:penaltyoverball}), the inequality 
		$$J_{\varepsilon}(\Ueps) \leq J_{\varepsilon}(\V), ~\forall \V\in \Buepssigma \cap \Bor$$
		holds for all $\sigma>0$. Then by the boundedness of $U_{\varepsilon}$, there exists $\sigma>0$ such that 
		$\Buepssigma \cap \{\pm 1\}^n\neq \emptyset$.
	Now, let us consider $$\U^*\in \argmin\{ \|\V - \Ueps\|_2: \V\in \Buepssigma \cap \{\pm 1\}^n\}.$$
	Since $\Buepssigma$ is a ball centered at $\Ueps$, we have necessarily (see Proposition \ref{prop:01})
	\[ \U^* \circ \Ueps \geq \oRn.\]
	It follows from $\U^*\in \{\pm 1\}^n$ and $\U^* \circ \Ueps \geq \oRn$ that
	\begin{align*}
		\|\Ueps \circ \Ueps - \one \|_2^2 &= \|\Ueps \circ \Ueps - \U^*\circ \U^* \|_2^2\\
		&=\|(\Ueps - \U^*)\circ (\Ueps + \U^*) \|_2^2\\
		&=\sum_{i=1}^{n} (\Ueps - \U^*)_i^2 (\Ueps + \U^*)_i^2\\
		&=\sum_{i=1}^{n} (\Ueps - \U^*)_i^2 \left((\Ueps)_i^2 + 2(\Ueps)_i(\U^*)_i + (\U^*)_i^2\right)\\
		&\geq \sum_{i=1}^{n} (\Ueps - \U^*)_i^2 =\|\Ueps - \U^*\|_2^2,
	\end{align*}
	that is 
	\begin{equation}
		\label{eq:normineq}
		\|\Ueps - \U^*\|_2 \leq \|\Ueps \circ \Ueps - \one \|_2.
	\end{equation}
	Taking $\V=\U^*$ in \cref{eq:localminimizer} and combining with \cref{eq:normineq}, we obtain 
	\begin{equation}
		\label{eq:ineq1}
		\frac{1}{4\varepsilon}\|\Ueps - \U^*\|_2^2 \leq \frac{c}{2} \left(\|\U^*\|_2^2 - \|\Ueps\|_2^2 \right) + \Pi(\U^*) -\Pi(\Ueps). 
	\end{equation}
	The mean value theorem and Cauchy-Schwartz inequality give
	\begin{equation}
		\label{eq:meanvalue}
		\begin{aligned}
			\left|\Pi(\U^*) - \Pi(\Ueps)\right|  &\leq  \max_{\V\in [\Ueps,\U^*]} | \langle \nabla \Pi(\V), \U^*-\Ueps\rangle | \\
			&\leq \max_{\V\in [\Ueps,\U^*]} \|\nabla \Pi(\V)\|_2 \|\U^*-\Ueps\|_2,
		\end{aligned}
	\end{equation}
	where $[\Ueps,\U^*]$ is the line segment between $\Ueps$ and $\U^*$. It follows that 
	\begin{eqnarray*}
		\frac{1}{4\varepsilon}\|\Ueps - \U^*\|_2^2 &\overset{\cref{eq:ineq1}}{\leq}& \frac{c}{2} \left(\|\U^*\|_2^2 - \|\Ueps\|_2^2 \right) + \left| \Pi(\U^*) -\Pi(\Ueps) \right|\\
		&\overset{\cref{eq:meanvalue}}{\leq}& \frac{c}{2} \left( 2\langle \U^*,  \U^*- \Ueps\rangle - \|\U^*-\Ueps\|_2^2\right) + \max_{\V\in [\Ueps,\U^*]} \|\nabla \Pi(\V)\|_2 \|\U^*-\Ueps\|_2\\	
		&\leq& c\|\U^*\|_2 \|\U^*- \Ueps\|_2 - \frac{c}{2}  \|\U^*-\Ueps\|_2^2  +  \max_{\V\in [\Ueps,\U^*]} \|\nabla \Pi(\V)\|_2 \|\U^*-\Ueps\|_2\\
		&=& c\sqrt{n} \|\U^*- \Ueps\|_2 - \frac{c}{2}  \|\U^*-\Ueps\|_2^2  +  \max_{\V\in [\Ueps,\U^*]} \|\nabla \Pi(\V)\|_2 \|\U^*-\Ueps\|_2,
	\end{eqnarray*}
	where the third inequality is due to the Cauchy-Schwartz inequality. 
	Dividing by $\|\U^*-\Ueps\|_2$, we obtain
	\[ \delta \leq \|\U^*-\Ueps\|_2 \leq \left(\frac{1}{4\varepsilon} + \frac{c}{2}\right)^{-1} \left( c\sqrt{n} +  \max_{\V\in [\Ueps,\U^*]} \|\nabla \Pi(\V)\|_2 \right),\]
	which (from $[\Ueps,\U^*]\subset \Bor$) proves that
	\[\delta \leq \frac{4\varepsilon}{1+2c\varepsilon} \left( c\sqrt{n} + \max_{\V\in \Bor} \|\nabla \Pi(\V)\|_2 \right).\] \qed
\end{proof}
\begin{remark}\label{rmk:6}
	Due to $\frac{4\varepsilon}{1+2c\varepsilon}\leq 4\varepsilon$ ($\forall \varepsilon>0$ and $c\geq 0$), one can simplify relation \cref{eq:errorbound}, obtaining $$\delta  \leq L(n,c,r,\Pi) \varepsilon$$
	with 
	$$L(n,c,r,\Pi) = 4\left( c\sqrt{n} +  \max_{\V\in \Bor} \|\nabla \Pi(\V)\|_2 \right).$$
	Since the polynomial functional $\Pi$ is Lipschitz continuous over $\Bor$, then there exists a Lipschitz constant $l_r$ such that $$\|\nabla \Pi(\V)\|_2 \leq l_r \sqrt{n}, ~\forall \V\in \Bor.$$
	implying  
	\begin{equation}
		\label{eq:simperrorbound}
		\delta \leq 4(c+l_r)\sqrt{n} \varepsilon.
	\end{equation}
\end{remark}

The error bounds in \eqref{eq:errorbound} and \eqref{eq:simperrorbound} are not able to be used directly in practice since the Lipschitz constant $l_r$ is not known and the optimization problem $\max_{\V\in \Bor} \|\nabla \Pi(\V)\|_2$ is non-convex and difficult to be solved. Fortunately, based on the simple fact that 
\begin{equation}
	\label{eq:norm2leqnorminf}
	\forall \X\in \R^n, ~ \|\X\|_2 \leq \|\X\|_{\infty} \sqrt{n},
\end{equation} since $$\|\X\|_2 = \sqrt{\sum_{i=1}^{n} x_i^2} \leq \sqrt{\sum_{i=1}^{n} \max_{i\in \{1,\ldots,n\}}(x_i^2)} = \max_{i\in \{1,\ldots,n\}}(|x_i|) \sqrt{n} = \|\X\|_{\infty}\sqrt{n},$$
then we can prove an explicit formulation of the error bound in the next theorem.
\begin{theorem}\label{thm:errorbound_explicitformula}
	The assumptions are those of Theorem \ref{thm:errorbound}. Let $\Pi(\V)$ be expressed as sum of monomials as $\sum_{\alpha} c_{\alpha} \V^{\alpha}$ where $\alpha = (\alpha_1,\ldots,\alpha_n)$ denotes the exponents of the monomial $\V^{\alpha} = \prod_{i=1}^{n}v_i^{\alpha_i}$, and $c_{\alpha}$ is the coefficient of the monomial $\V^{\alpha}$. Let $|\alpha|=\sum_{i=1}^{n}\alpha_i$. Then  
	\begin{equation}
		\label{eq:errorbound_explicitformula}
		\delta \leq \frac{4\varepsilon}{1+2c\varepsilon} \left( c + L_r \right) \sqrt{n}.
	\end{equation}
	where \begin{equation}
		\label{eq:L_r}
		L_r = \max_{i\in \{1,\ldots,n\}} \sum_{\alpha} |c_{\alpha}| |\alpha_i| r^{|\alpha|-1}.
	\end{equation}
\end{theorem}
\begin{proof}
	For any $\V\in \Bor$ ($r>\sqrt{n}$), it follows from \eqref{eq:norm2leqnorminf} that 
	\begin{equation}
		\label{eq:normineqofPi}
		\|\nabla \Pi(\V)\|_2 \leq \|\nabla \Pi(\V)\|_{\infty} \sqrt{n}.
	\end{equation}
	Then by the convention that for all $v_i\in \R$, $v_i^0=1$ and $v_i^{-1}=0$, we get the expression 
	\begin{equation}
		\label{eq:partialdiffofPi}
		\frac{\partial \Pi(\V)}{\partial v_i} = \sum_{\alpha} c_{\alpha} \alpha_i v_i^{\alpha_i-1} \prod_{j\neq i}v_j^{\alpha_j}.
	\end{equation}
	It follows from \eqref{eq:normineqofPi}, \eqref{eq:partialdiffofPi} and $|v_i|\leq r, \forall v\in \Bor$ that  
	\begin{align*}
		\|\nabla \Pi(\V)\|_{\infty} &= \max_{i\in \{1,\ldots,n\}} \left| \frac{\partial \Pi(\V)}{\partial v_i} \right|\\
		&= \max_{i\in \{1,\ldots,n\}} \left| \sum_{\alpha} c_{\alpha} \alpha_i v_i^{\alpha_i-1} \prod_{j\neq i}v_j^{\alpha_j} \right|\\
		&\leq \max_{i\in \{1,\ldots,n\}} \sum_{\alpha} |c_{\alpha}| |\alpha_i| |v_i|^{\alpha_i-1} \prod_{j\neq i}|v_j|^{\alpha_j} \\
		&\leq \max_{i\in \{1,\ldots,n\}} \sum_{\alpha} |c_{\alpha}| |\alpha_i| r^{\alpha_i-1} \prod_{j\neq i}r^{\alpha_j}\\
		&= \max_{i\in \{1,\ldots,n\}} \sum_{\alpha} |c_{\alpha}| |\alpha_i| r^{|\alpha|-1} = L_r.
	\end{align*}
	Then,
	\begin{equation}
		\label{eq:boundofnorm2}
		\|\nabla \Pi(\V)\|_2 \leq L_r \sqrt{n}, ~\forall \V\in \Bor.
	\end{equation}
	The desired error bound \eqref{eq:errorbound_explicitformula} is derived immediately from \eqref{eq:boundofnorm2} and \eqref{eq:errorbound}. \qed
\end{proof}
\begin{remark}\label{rmk:7}
	Similar to Remark \ref{rmk:6}, we can obtain a simplified error bound as  
	\begin{equation}
		\label{eq:simperrorbound_explicitformula}
		\delta \leq 4(c+L_r)\sqrt{n} \varepsilon,
	\end{equation}
	where $L_r$ is given by \eqref{eq:L_r}.
\end{remark}

Now, let us show the reason why we use an open ball $\Buepssigma$ to get the inequality $\U^* \circ \Ueps \geq \oRn$ in the proof of Theorem \ref{thm:errorbound}.
\begin{proposition}\label{prop:01}
	The assumptions are those of Theorem \ref{thm:errorbound}. Let $\Ueps$ be a local minimizer of functional $J_{\varepsilon}$ over the set $\Bor$ ($r>\sqrt{n}$), and \begin{equation}
		\label{eq:minimizerU*}
		\U^*\in \argmin\{ \|\V - \Ueps\|_2: \V\in \Buepssigma \cap \{\pm 1\}^n\}.
	\end{equation}
	Then we have 
	$$\U^* \circ \Ueps \geq \oRn.$$
\end{proposition} 
\begin{proof}	
	By contradiction, suppose that $\U^* \circ \Ueps \ngeq \oRn$, then there exists a nonempty index set $I = \{i\in \{1,\ldots,n\}: (\U^*)_i (\Ueps)_i < 0\}$. We will prove that there exists a point $\Ubar^* \in \Buepssigma \cap \{\pm 1\}^n$ such that $\|\Ubar^* - \Ueps\|_2 < \|\U^* - \Ueps\|_2$, which contradicts the assumption that $\U^*$ is a minimizer of \cref{eq:minimizerU*}. Such a point $\Ubar^*$ could be given by:
	$$(\Ubar^*)_i = \begin{cases}
		(\U^*)_i, &\text{if } i\in \{1,\ldots,n\}\setminus I,\\
		-(\U^*)_i, &\text{if } i\in I.
	\end{cases}$$
	We can verify that:\\
	$i)$ $\Ubar^* \in \{\pm 1\}^n$ is obvious.\\
	$ii)$ $\|\Ubar^* - \Ueps\|_2 < \|\U^* - \Ueps\|_2$ because:  
	$$\|\Ubar^* - \Ueps\|_2^2 = \|\Ubar^*\|_2^2 + \|\Ueps\|_2^2 - 2\sum_{i=1}^{n} (\Ubar^*)_i (\Ueps)_i,$$ 
	and by definition of $\Ubar^*$, we have $\|\Ubar^*\|_2^2 = \|\U^*\|_2^2$ and
	$$-2\sum_{i=1}^{n} (\Ubar^*)_i (\Ueps)_i < -2 \sum_{i=1}^{n} (\U^*)_i (\Ueps)_i,$$
	where the strict inequality is guaranteed by the non-emptiness of the index set $I$. Then  
	\begin{equation}
		\label{eq:eee}
		\|\Ubar^* - \Ueps\|_2^2  < \|\U^* - \Ueps\|_2^2.
	\end{equation}
	$iii)$ $\Ubar^*\in \Buepssigma$ because: $\U^*\in \Buepssigma$ implies $\|\U^* - \Ueps\|_2^2\leq \sigma^2$, then we get from inequality \cref{eq:eee} that $\|\Ubar^* - \Ueps\|_2^2 < \sigma^2$, which implies $\Ubar^*\in \Buepssigma.$
	
	Therefore, $\Ubar^*\in \Buepssigma \cap \{\pm 1\}^n$ and $\|\Ubar^* - \Ueps\|_2 < \|\U^* - \Ueps\|_2$ which contradicts the assumption that $\U^*$ is a minimizer of \cref{eq:minimizerU*}. Thus $\U^* \circ \Ueps \geq \oRn$. \qed
\end{proof}  
\begin{remark}
	The open ball $\Buepssigma$ can not be replaced by an arbitrary neighborhood around $\Ueps$, since the inequality $\|\Ubar^* - \Ueps\|_2 < \sigma$ cannot imply that $\Ubar^*$ belongs to this neighborhood.
\end{remark}

\section{On The Choice of The Parameters}\label{sec:choiceofparameters}
\subsection{Choice of $\varepsilon$}\label{subsec:choiceofepsilon}
Based on the error bound \eqref{eq:simperrorbound_explicitformula} in Theorem \ref{thm:errorbound_explicitformula}:
$$\delta \leq 4(c+L_r)\sqrt{n}\varepsilon,$$ for converging to a solution such that the error bound $\delta$ is smaller than any given desired precision $\zeta$ ($>0$), the parameter $\varepsilon$ can be chosen as 
\begin{equation}
	\label{eq:choiceofepsilon}
	\varepsilon \leq \frac{\zeta}{4(c+L_r)\sqrt{n}},
\end{equation}
where $$L_r = \max_{i\in \{1,\ldots,n\}} \sum_{\alpha} |c_{\alpha}| |\alpha_i| r^{|\alpha|-1}$$
is defined in \eqref{eq:L_r}. Note that decreasing (reasonably) $\varepsilon$ will lead to solutions closer to $\{\pm 1\}^n$. 

\subsection{Choice of $\gamma$ and $m$}\label{subsec:choiceofgamma&m}
As a matter of fact, the two parameters $\gamma$ and $m$ can be reduced to one, the ratio $\gamma/m$, which can be explained in physical damped system \cref{eq:prob_2ndordeode}. In fact, the Houbolt or the Runge-Kutta scheme applying to system \cref{eq:prob_2ndordeode} simulates the trajectory of a heavy ball rolling on the frictional hyper-surface (the shape of the surface is defined by the objective functional of the penalty problem \cref{prob:approx_int_cv_opt}) with $\gamma$ being the friction factor of the surface, and $m$ being the mass of the heavy ball. Therefore, the oscillations of the heavy ball will be suppressed when increasing $\gamma$ and decreasing $m$, because of the clear fact that a less heavy ball rolling on a more frictional surface will stop more quickly to a stationary point. Thus, we propose two choices for these parameters: (i) Fix the ratio $\gamma/m$ to some values (e.g., $m=1$, $\gamma=50$); (ii) Vary the ratio $\gamma/m$ by the way that 
\begin{equation}
	\label{eq:choiceofgamma&m}
	m=1 \text{ and } \gamma=\left\lbrace\begin{array}{ll}
		3, & \text{in \eqref{eq:u1}};\\
		3/(k+1), & \text{in \eqref{eq:prob_uk+1}}.
	\end{array}\right.
\end{equation}
The second choice is derived from the formulation $a(t)=3/t$ given in \cite{Su2015}, since $\gamma/m$ corresponds exactly to the damping term $a(t)$ in the Polyak Heavy Ball system \eqref{eq:nesterovode}. Note that a larger ratio $\gamma/m$ often leads to faster convergence.

\subsection{Choice of $c$}\label{subsec:choiceofc}
As suggested in Remark \ref{rmk:1}, we use parameter $c$ in order to enhance the convexity of the functional $\V\mapsto \frac{c}{2}\|\V\|_2^2 + \Pi(\V)$ over the ball $\Bor$ (with $r\geq \sqrt{n}$). For guaranteeing the convexity, we may take  
\begin{equation}
	\label{eq:infnormbound}
	c \geq \max_{\V\in \Bor} \rho( \nabla^2 \Pi(\V)),
\end{equation}
where $\rho( \nabla^2 \Pi(\V))$ stands for the spectral radius of $\nabla^2 \Pi(\V)$ (the Hessian matrix of functional $\Pi$ at $\V$). If $d>2$, the above maximization problem is neither convex nor concave in general, making its solution computationally expensive for large values of $n$. Since 
$$\rho( \nabla^2 \Pi(\V)) \leq \|\nabla^2 \Pi(\V)\|_{\infty},$$ 
a cheaper (but not as sharp) alternative to \cref{eq:infnormbound} is provided by
\begin{equation}
	\label{eq:estimofc}
	c\geq \max_{\V\in \Bor} \|\nabla^2 \Pi(\V)\|_{\infty}.
\end{equation}
Driven by \cref{eq:estimofc}, we obtain an estimation of $c$ in Proposition \ref{prop:estimatec}.
\begin{proposition}\label{prop:estimatec}
	Let $\pi_{i,j}(\V)$ be the element of the $i$th row and the $j$th column ($\forall i,j\in \{1,\ldots,n\}$) of the matrix $\nabla^2 \Pi(\V)$, which can be expanded as sum of monomials as $$\pi_{i,j}(\V) = \alpha_{(i,j),0} +  \sum_{k=1}^{k_{i,j}}\alpha_{(i,j),k} \psi_{(i,j),k}(\V)$$ with $\alpha_{(i,j),k}\in \R, \forall k=0,\ldots,k_{i,j}$ and $\psi_{(i,j),k}$ a monomial of degree $\beta_{(i,j),k}$ (with $1\leq \beta_{(i,j),k}\leq d-2$), then we can take
	\begin{equation}
		\label{eq:lowerboundofc}
		c\geq \max_{i\in \{1,\ldots,n\}}\sum_{j=1}^{n} |\alpha_{(i,j),0}| +  \sum_{k=1}^{k_{i,j}}|\alpha_{(i,j),k}| r^{\beta_{(i,j),k}}
	\end{equation}
	to guarantee the convexity of the functional $\V\mapsto \frac{c}{2}\|\V\|_2^2 + \Pi(\V)$ over the ball $\Bor$ with $r\geq \sqrt{n}$.
\end{proposition}
\begin{proof}
	By computing an upper bound of $|\pi_{i,j}(\V)|$ over $\Bor$ as
	$$|\pi_{i,j}(\V)|\leq |\alpha_{(i,j),0}| +  \sum_{k=1}^{k_{i,j}}|\alpha_{(i,j),k}| r^{\beta_{(i,j),k}}, \forall \V\in \Bor,$$
	we get 
	$$\begin{aligned}
		\max_{\V\in \Bor} \|\nabla^2 \Pi(\V)\|_{\infty} &= \max_{\V\in \Bor} \max_{i\in \{1,\ldots,n\}}\sum_{j=1}^{n}|\pi_{i,j}(\V)|\\
		&\leq \max_{i\in \{1,\ldots,n\}}\sum_{j=1}^{n} |\alpha_{(i,j),0}| +  \sum_{k=1}^{k_{i,j}}|\alpha_{(i,j),k}| r^{\beta_{(i,j),k}}.
	\end{aligned}$$
	It follows from \eqref{eq:estimofc} that, for guaranteeing the convexity of the functional $\V\to \frac{c}{2}\|\V\|_2^2 + \Pi(\V)$, we may choose $c$ as
	\[c\geq \max_{i\in \{1,\ldots,n\}}\sum_{j=1}^{n} |\alpha_{(i,j),0}| +  \sum_{k=1}^{k_{i,j}}|\alpha_{(i,j),k}| r^{\beta_{(i,j),k}}.\]\qed
\end{proof}
Recall (see Remark \ref{rmk:3}) that for $c$ large enough, penalty is necessarily inexact. Numerical results presented in Section \ref{sec:simulations} demonstrate that the parameter $c$ (not too large) has very small impact to our numerical schemes nevertheless the convexity of $\V\mapsto \frac{c}{2}\|\V\|_2^2 + \Pi(\V)$ over the ball $\Bor$. Therefore, we suggest to fix $c=0$ in most of practical applications.

\subsection{Choice of $\tau$}\label{subsec:choiceoftau}

There is no difficulty with the choice of $\tau$ for the Houbolt scheme. Indeed, relation \cref{eq:boundsoftau}, namely $\tau \leq \sqrt{2m\varepsilon}$ suggests taking $\tau$ as \begin{equation}\label{eq:tau_Houbolt}
	0< \tau = \sqrt{2m\varepsilon}.
\end{equation}
The choice of $\tau$ for the Lie-Marchuk-Yanenko scheme is a little more complicated since a small $\tau$ is required to reduce the splitting error; on the other hand a small $\tau$ may imply a large number of iterations to reach convergence. There are two ways, in practice, to choose parameter $\tau$. The first method consists taking $\tau$ verifying \cref{eq:convcond_Lie} ($\tau = \varepsilon/(1-\varepsilon c)$, for example). The second method relies on variable time steps; for example, we can use a sequence $\{\tau_k\}_{k\geq 0}$ of variable time steps verifying
\begin{equation}
	\label{eq:steptau}
	\left\lbrace 
	\begin{aligned}
		&\forall k\geq 0, \tau_k>0 \text{ and } \tau_k > \tau_{k+1},\\
		&\lim\limits_{k\to +\infty} \tau_k = 0,\\
		&\sum_{k=1}^{+\infty} \tau_k = +\infty.
	\end{aligned}
	\right.
\end{equation}
However, this method often leads to small time steps and slow down the convergence when $k$ increases. The method we used in this article is a compromise between the two above approaches: Starting with $\tau_0$ as
\begin{equation}
	\label{eq:tau0_Lie}
	0< \tau_0 = \min\left\{\frac{\varepsilon}{1-\varepsilon c},0.1\right\} \leq \frac{\varepsilon}{1-\varepsilon c},
\end{equation}
for $k\geq 0$, we compute $\tau_{k+1}$ from $\tau_{k}$ via the relation
\begin{equation}
	\label{eq:updatetau}
	\tau_{k+1} = \left\lbrace \begin{array}{ll}
		\theta \tau_k, & \text{if }\tau_k \geq \tau^*,\\
		\tau_k, & \text{otherwise,}
	\end{array} \right.
\end{equation}
where $\theta\in (0,1)$ is a reduction ratio, e.g., $\theta = 0.8$, and $\tau^*$ is a threshold (i.e., the sequence $\{\tau_k\}_{k\geq 0}$ is decreasing and belongs to the interval $[\tau^*,\tau_0]$).

\begin{remark}\label{rmk:5}
	The Houbolt scheme does not suffer from a splitting error, authorizing therefore larger time steps than the Lie-Marchuk-Yanenko scheme for the same value of $\varepsilon$.
\end{remark}

\section{Numerical Experiments}\label{sec:simulations}
Our code, namely \texttt{DEMIPP} \cite{DEMIPP} (Differential Equation Methods for Integer Polynomial Programs) is developed on MATLAB and shared on Github\footnote{DEMIPP is available at \url{https://github.com/niuyishuai/DEMIPP}.}. The numerical simulations are performed on a cluster in the department of mathematics at Shanghai Jiao Tong University, where $80$ CPUs (Intel Xeon Gold 6148 CPU @ 2.40GHz) are used for parallel computation. The proposed three ODE approaches (namely the Houbolt scheme, the Lie scheme and the Runge-Kutta (4,5) scheme) for the formulation \eqref{prob:approx_int_cv_opt} are compared with two nonlinear optimization approaches (a quadratic binary formulation approach \cite{Boros2002} for the formulation \eqref{prob:Boolean_opt}; and the IPOPT solver \cite{Ipopt} for the formulation \eqref{prob:penaltyoverball}). 
\paragraph{\textbf{Setups}}\label{para:setups} The setups of these methods are summarized as follows:
\begin{enumerate}
	\item The Houbolt scheme (cf., Houbolt): the method is described in Subsection \ref{subsec:houbolt}. We set parameters $\varepsilon\in \{10^{-4}, 10^{-5}, 10^{-6}\}$, $\gamma=50$, $m=1$, $\tau=\sqrt{2m\varepsilon}$ and $c=100$. 
	\item The Lie scheme (cf., Lie): the method is given in Subsection \ref{subsec:Lie}, and the parameters are given by $\varepsilon\in \{10^{-4}, 10^{-5}, 10^{-6}\}$, $\tau=\min\{\varepsilon/(1-c\varepsilon), 0.1\}$, and $c=100$. 
	\item The Runge-Kutta (4,5) (cf., RK(4,5)) scheme: the method is used to handle the first-order ODE formulation \cref{prob:1storderode}. We use an existing implementation in MATLAB, function \verb|ode45|, based on an explicit Runge-Kutta (4,5) formula and the Dormand-Prince pair, with default \verb|ode45|'s MATLAB settings. The parameters required in formulation \cref{prob:1storderode} is given by $\varepsilon\in \{10^{-4}, 10^{-5}, 10^{-6}\}$, $\gamma=50$, $m=1$, and $c=100$. The initial time $t$ for RK(4,5) is set to be $0$, and the final time is often set to be smaller than $1$. It should be pointed out that the higher order Runge-Kutta schemes require more gradient computations per iteration, thus could be more expensive methods.
	\item Quadratic binary formulation + GUROBI (cf., QB-G): we use a quadratic binary program (QBP) formulation for problem \eqref{prob:Boolean_opt} proposed in \cite{Boros2002}. The QBPs are solved using GUROBI solver \cite{Gurobi} with default GUROBI settings.
	\item IPOPT: A software package for large-scale nonlinear optimization using interior point method. We use IPOPT with default settings for problem \eqref{prob:penaltyoverball}.
\end{enumerate} 

Note that the system \cref{eq:system_uk0.5_lie} (nonlinear in general) required in the Lie scheme is solved by MATLAB nonlinear equation solver \texttt{fsolve} with default parameters. The stopping tolerances \tolf\ and \tolu\ are only needed in ODE approaches and fixed to $10^{-4}$ and $10^{-2}$ as default settings throughout our tests. Moreover, all tested methods are deterministic without randomness, so multiple runs of a method with the same initial point and the same setting will provide almost the same numerical result with very slight difference. Hence, we report the numerical results of one run for each tested problem.

Table \ref{tab:params} summarizes all required parameters and the corresponding problem formulations for different solution methods:
\begin{table}[htb!]
	\begin{center}
		\caption{Required parameters and problem formulations for different methods.}
		\label{tab:params}
		\begin{tabular}{r|ccccc|ccc|cc} \hline
			Method & $\varepsilon$ & $m$ & $\gamma$ & $c$ & $\tau$ & $t$ & \tolf & \tolu & Formulation & Variable \\
			\hline
			Houbolt & \Checkmark & \Checkmark & \Checkmark & \Checkmark & \Checkmark & \XSolidBrush & \Checkmark & \Checkmark & \eqref{prob:approx_int_cv_opt} & $\{\pm 1\}^n$\\
			Lie & \Checkmark & \XSolidBrush & \XSolidBrush & \Checkmark & \Checkmark & \XSolidBrush & \Checkmark & \Checkmark & \eqref{prob:approx_int_cv_opt} & $\{\pm 1\}^n$\\
			RK(4,5) & \Checkmark & \Checkmark & \Checkmark & \Checkmark & \XSolidBrush & \Checkmark & \Checkmark & \Checkmark & \eqref{prob:approx_int_cv_opt} & $\{\pm 1\}^n$\\
			IPOPT & \Checkmark & \XSolidBrush & \XSolidBrush & \Checkmark & \XSolidBrush & \XSolidBrush & \XSolidBrush & \XSolidBrush & \eqref{prob:penaltyoverball} & $\{\pm 1\}^n$\\
			QB-G & \XSolidBrush & \XSolidBrush & \XSolidBrush & \XSolidBrush & \XSolidBrush & \XSolidBrush & \XSolidBrush & \XSolidBrush & \eqref{prob:Boolean_opt} & $\{0,1\}^n$\\
			\hline
		\end{tabular}
	\end{center}	
\end{table}
Particularly, the QB-G method only handles the problem with variable $\Y\in\{0,1\}^n$, while the other methods are designed for problems with variable $\V\in\{\pm 1\}^n$. Therefore, when testing on the same set of problems for all methods, there is a need to convert variables from $\{0,1\}^n$ to $\{\pm 1\}^n$ using the transformation \eqref{eq:lintrans}, and the corresponding polynomial functional will be converted accordingly.
\paragraph{\textbf{Datasets}} We perform our experiments on two datasets: (i) A randomly generated synthetic dataset where the coefficients of polynomials are integers randomly chosen in the interval $[-10,10]$, and the decision variable $\V\in \{\pm 1\}^n$ with $2\leq n\leq 10$ for small-scale cases and $11\leq n\leq 20$ for relatively large-scale cases. The degree of polynomials $d$ is chosen in the interval $[2,6]$ which covers the most frequently used polynomials in real-world applications. (ii) A benchmark dataset \texttt{MQLib} collected by Dunning, Gupta and Silberholz in \cite{dunning2018}. This dataset consists
	large-scale, heterogeneous ``Max-Cut" and ``Quadratic Unconstrained Binary Optimization" (QUBO) problems, combining real-world problem instances and random problem instances from multiple random generators. Both of them are quadratic Boolean programs.
\paragraph{\textbf{Polynomial modeling}}
	For efficiently modeling multivariate polynomials on MATLAB, we use POLYLAB \cite{Polylab} (a multivariate polynomial modeling toolbox on MATLAB developed by Y.S. Niu, code available on Github) which is chosen for its high efficiency on polynomial construction and elementary operations (such as multiplications, additions, and derivatives of polynomials). It is observed to be more efficient than other existing MATLAB toolboxes such as the MATLAB official symbolic toolbox (more than $100$ times slower than POLYLAB) and Yalmip \cite{Yalmip} (more than $5$ times slower than POLYLAB).

\subsection{Tests on Randomly Generated Synthetic Datasets}\label{subsec:testsonrandomdataset}
In this subsection, we will test on two randomly generated datasets with $2\leq n\leq 10$, $2\leq d\leq 4$ (small-scale cases), and with $10\leq n\leq 20$, $5\leq d\leq 6$ (relatively large-scale cases). We are interested in the computation time, the quality of the solution measured by $\delta$, as well as the impacts of parameters $\varepsilon$, $m$, $\gamma$, $c$ and $\tau$ for each method.

\subsubsection{Numerical Results on Small-scale Randomly Generated Dataset}\label{subsec:resultonsmall-scale}
The tested dataset is randomly generated with $2\leq n\leq 10$ and $2\leq d\leq 4$. Let us try an easy setting of parameters as $\varepsilon=10^{-4}$, $m=1$, $\gamma=50$, $c=100$, $\tau=\sqrt{2m\varepsilon} \approx 0.0141$ for the Houbolt scheme, $\tau = \min\{\varepsilon/(1-c\varepsilon), 0.1\} \approx 1.01\times 10^{-4}$ for the Lie scheme, and taking $t\in [0,0.3]$ for the RK(4,5) scheme. The starting points are same for different method and randomly chosen on unit sphere centered at $\oRn$. The stopping tolerances $\tolf=10^{-4}$ and $\tolu=10^{-2}$. The numerical results are summarized in Table \ref{tab:results1} where the columns $n$ and $d$ indicate the number of variables and the degree of polynomial $\Pi$ in problem \eqref{prob:int_opt}; the column \texttt{time} is the wall-clock computation time (unit in second); the column \texttt{iter} is the number of iterations; the column $\delta$ is the Euclidean distance between the computed solution $\Ueps$ and its closest integer point $\texttt{round}(\Ueps)$, i.e., $$\delta = \|\Ueps - \texttt{round}(\Ueps)\|;$$
and the column \texttt{obj} stands for the value of $\Pi(\texttt{round}(\Ueps))$.

\begin{table}[htb!]
	\begin{center}
		\caption{Numerical results on small-scale dataset with $\varepsilon=10^{-4}$, $m=1$, $\gamma=50$, $c=100$, $\tau = 0.0141$ for the Houbolt scheme, $\tau=1.01\times 10^{-4}$ for the Lie scheme, $t\in [0,0.3]$ for the RK(4,5) scheme, \texttt{Tolf}=$10^{-4}$, \texttt{TolU}=$10^{-2}$, and with default parameters for IPOPT and GUROBI solvers.}
		\label{tab:results1}
		\resizebox{\linewidth}{!}{%
			\tabcolsep=2pt
			\begin{tabular}{c|c|cccc|cccc|cccc|cccc|cc} \hline
				\multirow{2}{*}{n} & \multirow{2}{*}{d} & \multicolumn{4}{c}{Houbolt} & \multicolumn{4}{|c}{Lie} & \multicolumn{4}{|c}{RK(4,5)} & \multicolumn{4}{|c}{IPOPT} & \multicolumn{2}{|c}{QB-G}\\
				\cline{3-20}
				& & obj & iter & time & $\delta$ & obj & iter & time & $\delta$ & obj & iter & time & $\delta$ & obj & iter & time & $\delta$ & obj & time \\
				\hline\hline
				$2$ & $2$ & $-19$ & $12$ & $0.00$ & $1.00e-02$ &$-19$ & $6$ & $0.01$ & $8.66e-03$ &$-19$ & $125$ & $0.01$ & $1.00e-02$ &$-19$ & $11$ & $0.02$ & $6.10e-03$ &$-19$ & $0.01$ \\
				$2$ & $3$ & $35$ & $13$ & $0.00$ & $1.12e-02$ &$35$ & $7$ & $0.01$ & $1.04e-02$ &$35$ & $137$ & $0.02$ & $1.00e-02$ &$35$ & $15$ & $0.02$ & $8.82e-03$ &$21$ & $0.00$ \\
				$2$ & $4$ & $28$ & $12$ & $0.00$ & $2.79e-03$ &$28$ & $6$ & $0.01$ & $1.43e-02$ &$28$ & $153$ & $0.02$ & $1.00e-02$ &$28$ & $11$ & $0.02$ & $1.07e-02$ &$-6$ & $0.00$ \\
				$4$ & $2$ & $16$ & $12$ & $0.00$ & $9.05e-03$ &$16$ & $7$ & $0.01$ & $1.23e-02$ &$16$ & $157$ & $0.02$ & $1.00e-02$ &$16$ & $13$ & $0.02$ & $1.06e-02$ &$-20$ & $0.01$ \\
				$4$ & $3$ & $-23$ & $14$ & $0.00$ & $1.11e-02$ &$-23$ & $7$ & $0.01$ & $1.10e-02$ &$-23$ & $169$ & $0.02$ & $1.00e-02$ &$-23$ & $19$ & $0.03$ & $8.54e-03$ &$9$ & $0.02$ \\
				$4$ & $4$ & $-29$ & $14$ & $0.00$ & $1.33e-02$ &$-29$ & $6$ & $0.01$ & $1.13e-02$ &$-29$ & $165$ & $0.07$ & $1.00e-02$ &$-29$ & $12$ & $0.02$ & $7.14e-03$ &$-49$ & $0.02$ \\
				$6$ & $2$ & $-23$ & $14$ & $0.00$ & $2.05e-02$ &$-23$ & $7$ & $0.02$ & $1.32e-02$ &$-23$ & $161$ & $0.07$ & $1.00e-02$ &$-23$ & $12$ & $0.02$ & $1.14e-02$ &$-35$ & $0.01$ \\
				$6$ & $3$ & $-78$ & $14$ & $0.00$ & $1.65e-02$ &$-78$ & $7$ & $0.02$ & $9.90e-03$ &$-78$ & $161$ & $0.02$ & $1.00e-02$ &$-78$ & $14$ & $0.03$ & $8.22e-03$ &$32$ & $0.02$ \\
				$6$ & $4$ & $-32$ & $15$ & $0.00$ & $1.66e-02$ &$-32$ & $7$ & $0.02$ & $1.41e-02$ &$-32$ & $257$ & $0.03$ & $1.22e-02$ &$-32$ & $15$ & $0.02$ & $1.24e-02$ &$-76$ & $0.03$ \\
				$8$ & $2$ & $-15$ & $15$ & $0.00$ & $1.93e-02$ &$-15$ & $7$ & $0.02$ & $1.63e-02$ &$-15$ & $261$ & $0.03$ & $1.35e-02$ &$-15$ & $16$ & $0.02$ & $1.36e-02$ &$-89$ & $0.02$ \\
				$8$ & $3$ & $16$ & $15$ & $0.00$ & $2.24e-02$ &$16$ & $7$ & $0.02$ & $1.86e-02$ &$16$ & $257$ & $0.03$ & $1.62e-02$ &$16$ & $14$ & $0.02$ & $1.61e-02$ &$-80$ & $0.04$ \\
				$8$ & $4$ & $21$ & $16$ & $0.00$ & $1.79e-02$ &$21$ & $7$ & $0.03$ & $2.04e-02$ &$21$ & $257$ & $0.06$ & $1.88e-02$ &$21$ & $20$ & $0.05$ & $1.89e-02$ &$-45$ & $0.07$ \\
				$10$ & $2$ & $35$ & $16$ & $0.00$ & $2.00e-02$ &$35$ & $7$ & $0.03$ & $2.05e-02$ &$35$ & $265$ & $0.04$ & $1.66e-02$ &$35$ & $17$ & $0.03$ & $1.69e-02$ &$-95$ & $0.02$ \\
				$10$ & $3$ & $23$ & $15$ & $0.00$ & $2.19e-02$ &$23$ & $7$ & $0.02$ & $2.06e-02$ &$23$ & $261$ & $0.04$ & $1.87e-02$ &$23$ & $13$ & $0.02$ & $1.88e-02$ &$-89$ & $0.05$ \\
				$10$ & $4$ & $47$ & $15$ & $0.00$ & $3.02e-02$ &$47$ & $7$ & $0.04$ & $2.70e-02$ &$47$ & $265$ & $0.06$ & $2.58e-02$ &$47$ & $16$ & $0.03$ & $2.54e-02$ &$-225$ & $0.08$ \\
				\hline
				\multicolumn{2}{c|}{average} & & $14$ & $0.00$ & $1.62e-02$ & & $7$ & $0.02$ & $1.52e-02$ & & $203$ & $0.04$ & $1.35e-02$ & & $15$ & $0.02$ & $1.29e-02$ & & $0.03$ \\
				\hline
		\end{tabular}}
	\end{center}	
\end{table}

\paragraph{\textbf{Observations}}
We can observe in Table \ref{tab:results1} that the three ODE methods and IPOPT almost always converge to the same computed solutions very closely to the set $\{\pm 1\}^n$ with average $\delta$ of order $O(10^{-2})$. Here, we have checked that the computed solutions for these methods on the test dataset are always the same when the objective values are the same obtained from the same initial point. The QB-G method quite often provides best computed solutions for some small-scale problems. 
However, it is interesting to observe that our ODE approaches may provide better solutions than QB-G, e.g., the cases $(n,d)=(4,3)$ and $(6,3)$. These cases seem to occur more frequently for the larger $n$ and $d$ (see test results on large-scale dataset in Table \ref{tab:results2} for more details). Despite these differences, the computation time for all tested methods are very comparable and often less than $0.05$ second in average. Note that the computation time for the Houbolt scheme is always 0.00 since it is less than $0.01$ seconds, and this precision is enough to show the difference in the computation time of different methods. The RK(4,5) requires more time steps than the Houbolt and Lie schemes, which is not a surprise since the Runge-Kutta scheme is an explicit scheme which often requires more time steps than semi-implicit or implicit schemes. 

It is worth noting that our ODE methods are coded in MATLAB while GUROBI and IPOPT are developed in C++ and Fortran. In general, a MATLAB implementation will be slower than a C++ implementation, so that our ODE methods should be more faster than we presented for C++ implementation.

Next, we will discuss the impact of all parameters $\varepsilon$, $m$, $\gamma$, $c$ and $\tau$ for each algorithm on small-scale dataset.

\paragraph{\textbf{Impact of $\varepsilon$}}
As we can see in Table \ref{tab:params} that the parameter $\varepsilon$ will affect all methods except QB-G, since the QB-G method solves a quadratic Boolean formulation without introducing the parameter $\varepsilon$, so that its numerical solution will be always integers, i.e., $\delta=0$. Differently, as we observed in Table \ref{tab:results1} that all other methods using continuous penalty formulations (i.e., the Houbolt scheme, the Lie scheme, the RK(4,5) scheme, and IPOPT) yield numerical results close but not equal to integers whose gap is measured by $\delta$.

In general, a smaller $\varepsilon$ will lead to results with a smaller $\delta$ which is guaranteed theoretically by the inequality in Theorem \ref{thm:errorbound} as: $$\delta \leq 4(c+l_r)\sqrt{n} \varepsilon.$$ 
Numerical tests show that choosing $\varepsilon$ small enough, such that $\sqrt{n}\varepsilon=O(10^{-4})$ at least, will often provide good numerical results.

Table \ref{tab:results2} illustrates the numerical results on the same dataset and with the same parameters as in Table \ref{tab:results1} except that $\varepsilon$ is reduced from $10^{-4}$ to $10^{-6}$. Comparing Tables \ref{tab:results1} and \ref{tab:results2}, the computation time for all methods are almost in the same level, but the precision of the computed solutions measured by $\delta$ is improved from $O(10^{-2})$ to $O(10^{-3})$ for Houbolt, Lie and RK(4,5) schemes, and to $O(10^{-4})$ for IPOPT. Surprisingly, the number of iterations for the Houbolt and Lie schemes increase very slightly with the decrease of $\varepsilon$, but increased a lot for the RK(4,5) scheme. It seems that the RK(4,5) scheme is much more sensitive to the parameter $\varepsilon$.
\begin{table}[htb!]
	\begin{center}
		\caption{Numerical results on small-scale dataset with $\varepsilon=10^{-6}$, and same other parameters as in Table \ref{tab:results1}.}
		\label{tab:results2}
		\resizebox{\linewidth}{!}{%
			\tabcolsep=2pt
			\begin{tabular}{c|c|cccc|cccc|cccc|cccc|cc} \hline
				\multirow{2}{*}{n} & \multirow{2}{*}{d} & \multicolumn{4}{c}{Houbolt} & \multicolumn{4}{|c}{Lie} & \multicolumn{4}{|c}{RK(4,5)} & \multicolumn{4}{|c}{IPOPT} & \multicolumn{2}{|c}{QB-G}\\
				\cline{3-20}
				& & obj & iter & time & $\delta$ & obj & iter & time & $\delta$ & obj & iter & time & $\delta$ & obj & iter & time & $\delta$ & obj & time \\
				\hline\hline
				$2$ & $2$ & $-17$ & $15$ & $0.00$ & $6.23e-03$ &$-17$ & $5$ & $0.01$ & $2.85e-03$ &$-17$ & $929$ & $0.08$ & $1.00e-02$ &$-17$ & $16$ & $0.02$ & $6.16e-05$ &$-19$ & $0.00$ \\
				$2$ & $3$ & $21$ & $16$ & $0.00$ & $5.32e-03$ &$21$ & $7$ & $0.01$ & $1.82e-03$ &$21$ & $949$ & $0.07$ & $1.00e-02$ &$21$ & $21$ & $0.02$ & $8.18e-05$ &$21$ & $0.00$ \\
				$2$ & $4$ & $34$ & $16$ & $0.00$ & $5.13e-03$ &$34$ & $6$ & $0.01$ & $3.00e-03$ &$34$ & $1029$ & $0.09$ & $1.00e-02$ &$34$ & $15$ & $0.02$ & $1.02e-04$ &$-6$ & $0.00$ \\
				$4$ & $2$ & $10$ & $16$ & $0.00$ & $8.70e-03$ &$10$ & $6$ & $0.01$ & $5.03e-03$ &$10$ & $841$ & $0.11$ & $3.26e-02$ &$10$ & $16$ & $0.04$ & $1.05e-04$ &$-20$ & $0.01$ \\
				$4$ & $3$ & $-23$ & $16$ & $0.00$ & $8.25e-03$ &$-23$ & $7$ & $0.02$ & $1.85e-03$ &$-23$ & $1241$ & $0.19$ & $1.00e-02$ &$-23$ & $17$ & $0.04$ & $8.52e-05$ &$9$ & $0.02$ \\
				$4$ & $4$ & $55$ & $16$ & $0.00$ & $7.86e-03$ &$55$ & $7$ & $0.01$ & $2.03e-03$ &$55$ & $1265$ & $0.16$ & $1.00e-02$ &$55$ & $14$ & $0.02$ & $1.56e-04$ &$-49$ & $0.01$ \\
				$6$ & $2$ & $23$ & $16$ & $0.00$ & $9.87e-03$ &$23$ & $7$ & $0.01$ & $2.17e-03$ &$23$ & $1373$ & $0.14$ & $1.00e-02$ &$23$ & $14$ & $0.02$ & $1.31e-04$ &$-35$ & $0.01$ \\
				$6$ & $3$ & $-32$ & $17$ & $0.00$ & $6.00e-03$ &$-32$ & $7$ & $0.01$ & $2.91e-03$ &$-32$ & $1417$ & $0.17$ & $1.00e-02$ &$-32$ & $18$ & $0.02$ & $1.15e-04$ &$32$ & $0.02$ \\
				$6$ & $4$ & $-84$ & $16$ & $0.00$ & $8.42e-03$ &$-84$ & $6$ & $0.01$ & $4.40e-03$ &$-84$ & $1353$ & $0.16$ & $1.00e-02$ &$-84$ & $16$ & $0.02$ & $9.73e-05$ &$-76$ & $0.02$ \\
				$8$ & $2$ & $-55$ & $17$ & $0.00$ & $5.41e-03$ &$-55$ & $7$ & $0.02$ & $3.00e-03$ &$-55$ & $1425$ & $0.17$ & $1.00e-02$ &$-55$ & $15$ & $0.02$ & $1.24e-04$ &$-89$ & $0.02$ \\
				$8$ & $3$ & $-2$ & $17$ & $0.00$ & $5.89e-03$ &$-2$ & $7$ & $0.02$ & $2.97e-03$ &$-2$ & $1493$ & $0.21$ & $1.00e-02$ &$-2$ & $18$ & $0.03$ & $1.50e-04$ &$-80$ & $0.03$ \\
				$8$ & $4$ & $-101$ & $17$ & $0.00$ & $7.18e-03$ &$-101$ & $7$ & $0.02$ & $3.16e-03$ &$-101$ & $2325$ & $0.37$ & $6.13e-04$ &$-101$ & $16$ & $0.03$ & $7.40e-05$ &$-45$ & $0.08$ \\
				$10$ & $2$ & $5$ & $17$ & $0.00$ & $7.17e-03$ &$5$ & $7$ & $0.02$ & $3.45e-03$ &$5$ & $1545$ & $0.18$ & $1.00e-02$ &$5$ & $15$ & $0.02$ & $1.58e-04$ &$-95$ & $0.02$ \\
				$10$ & $3$ & $-83$ & $18$ & $0.00$ & $5.03e-03$ &$-83$ & $7$ & $0.02$ & $3.45e-03$ &$-83$ & $1509$ & $0.21$ & $1.00e-02$ &$-83$ & $21$ & $0.04$ & $1.27e-04$ &$-89$ & $0.05$ \\
				$10$ & $4$ & $33$ & $18$ & $0.00$ & $5.35e-03$ &$33$ & $7$ & $0.03$ & $3.26e-03$ &$33$ & $1517$ & $0.32$ & $1.00e-02$ &$33$ & $18$ & $0.03$ & $2.13e-04$ &$-225$ & $0.09$ \\
				\hline
				\multicolumn{2}{c|}{average} & & $17$ & $0.00$ & $6.79e-03$ & & $7$ & $0.02$ & $3.02e-03$ & & $1347$ & $0.18$ & $1.09e-02$ & & $17$ & $0.03$ & $1.19e-04$ & & $0.03$ \\
				\hline
		\end{tabular}}
	\end{center}	
\end{table}

\begin{figure*}[htb!]
	\centering
	\subfigure[$\varepsilon=10^{-4}, \texttt{iter}=533, \delta=1.19\times 10^{-2}$]{ 
		\includegraphics[width=0.45\textwidth]{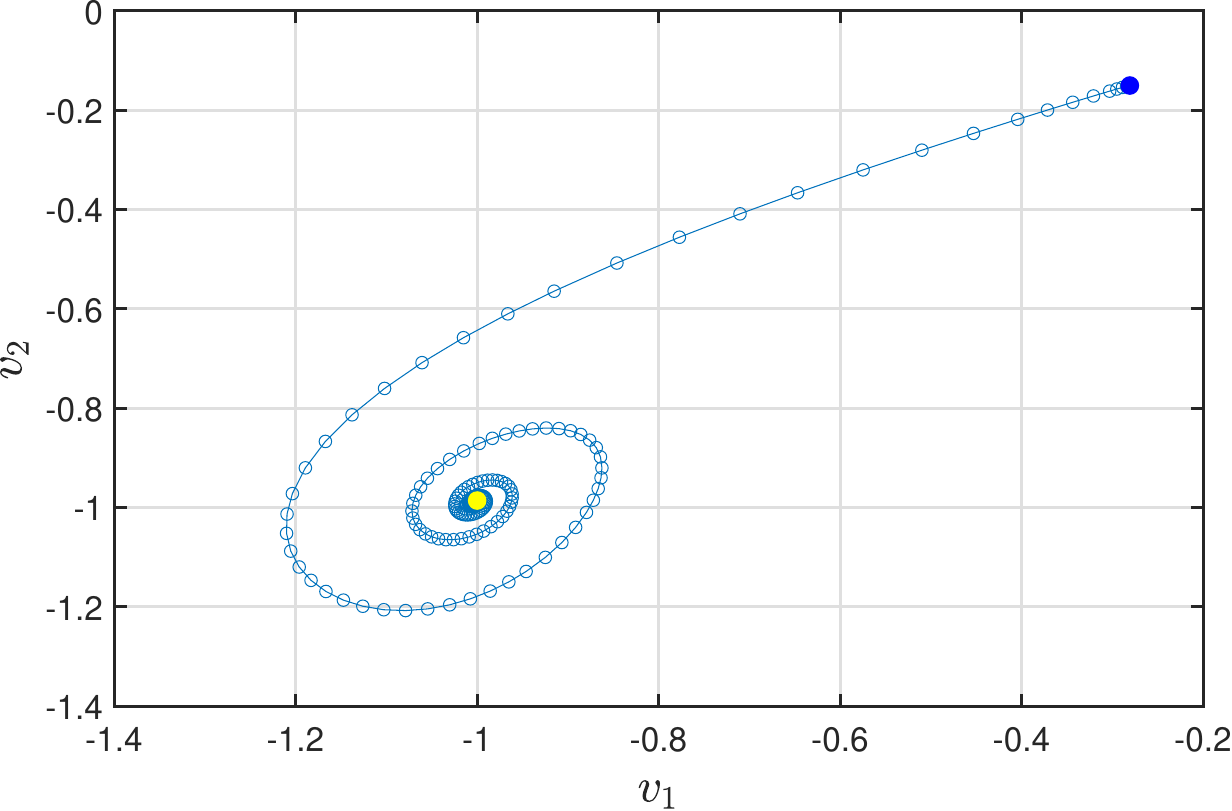}
	}
	\subfigure[$\varepsilon=10^{-5}, \texttt{iter}=1721, \delta=1.19\times 10^{-3}$]{ 
		\includegraphics[width=0.45\textwidth]{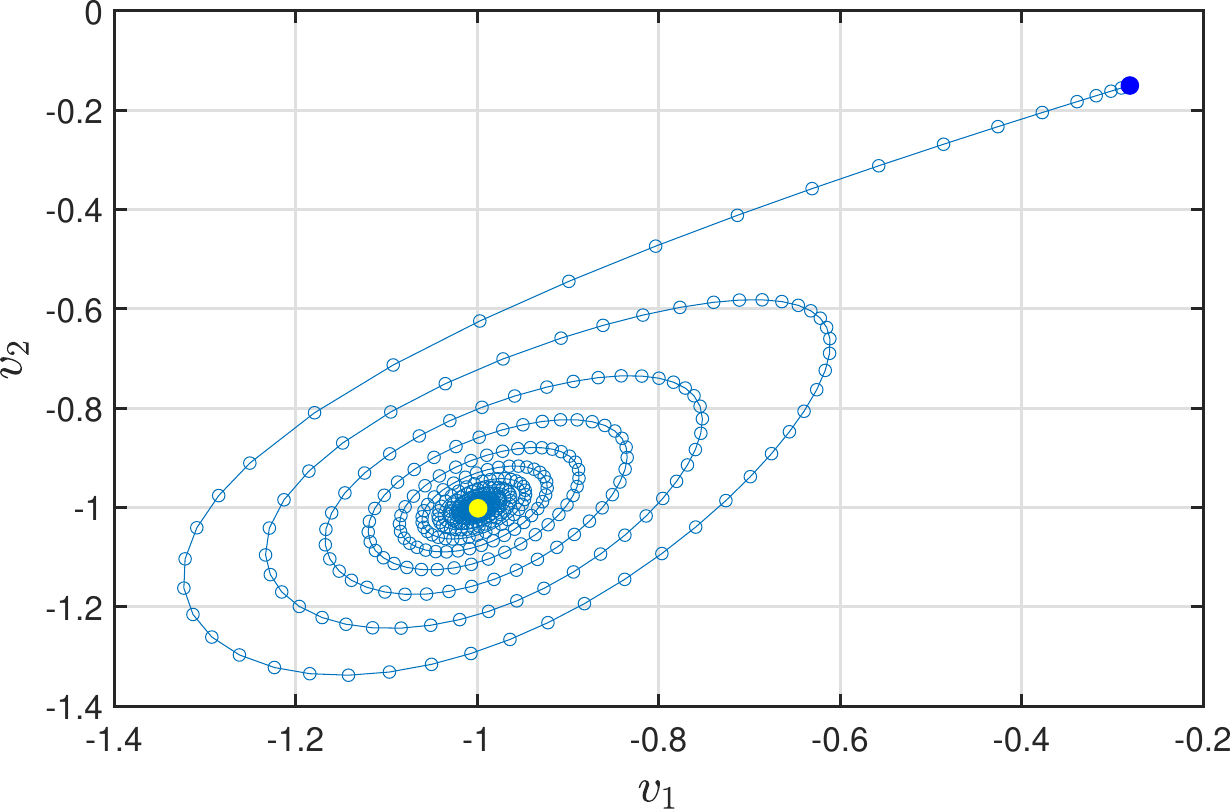}
	}
	\subfigure[$\varepsilon=10^{-7}, \texttt{iter}=7649, \delta=1.71\times 10^{-4}$]{ 
		\includegraphics[width=0.45\textwidth]{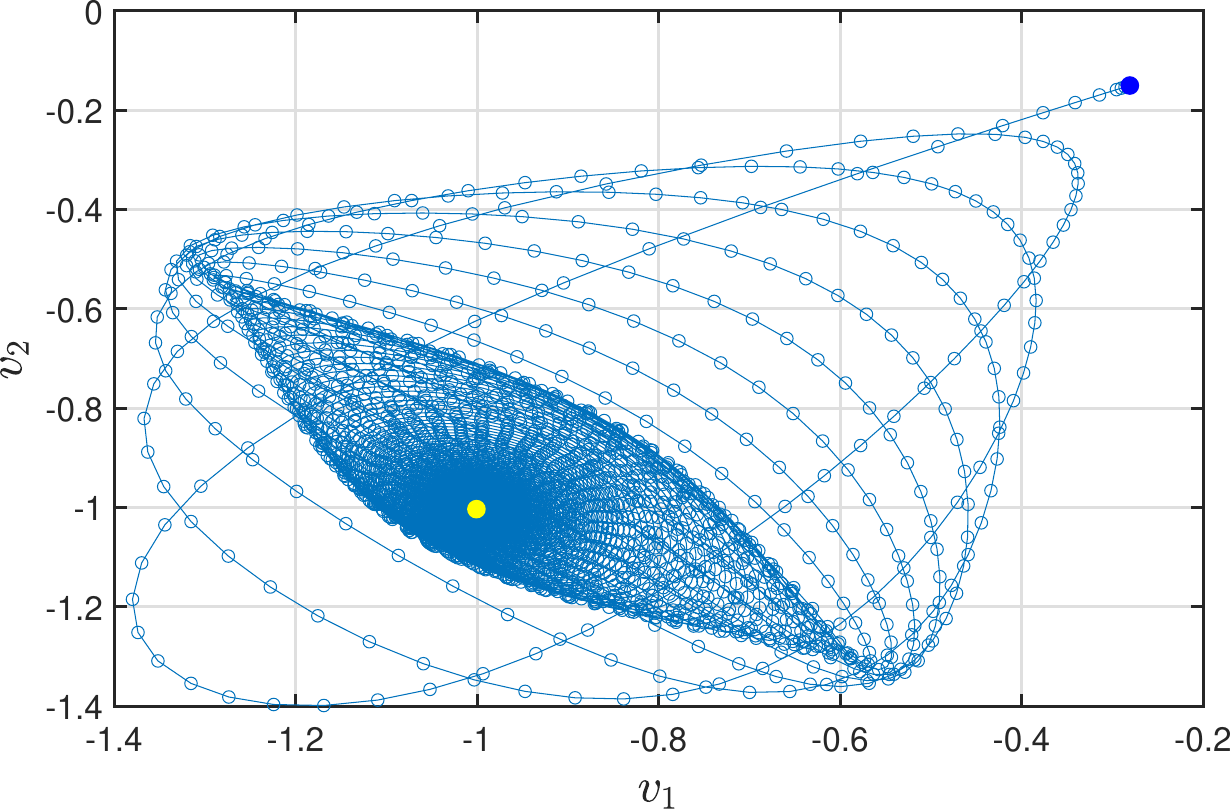} 
	}
	\subfigure[$\varepsilon=10^{-8}, \texttt{iter}=28669, \delta=2.64\times 10^{-5}$]{ 
		\includegraphics[width=0.45\textwidth]{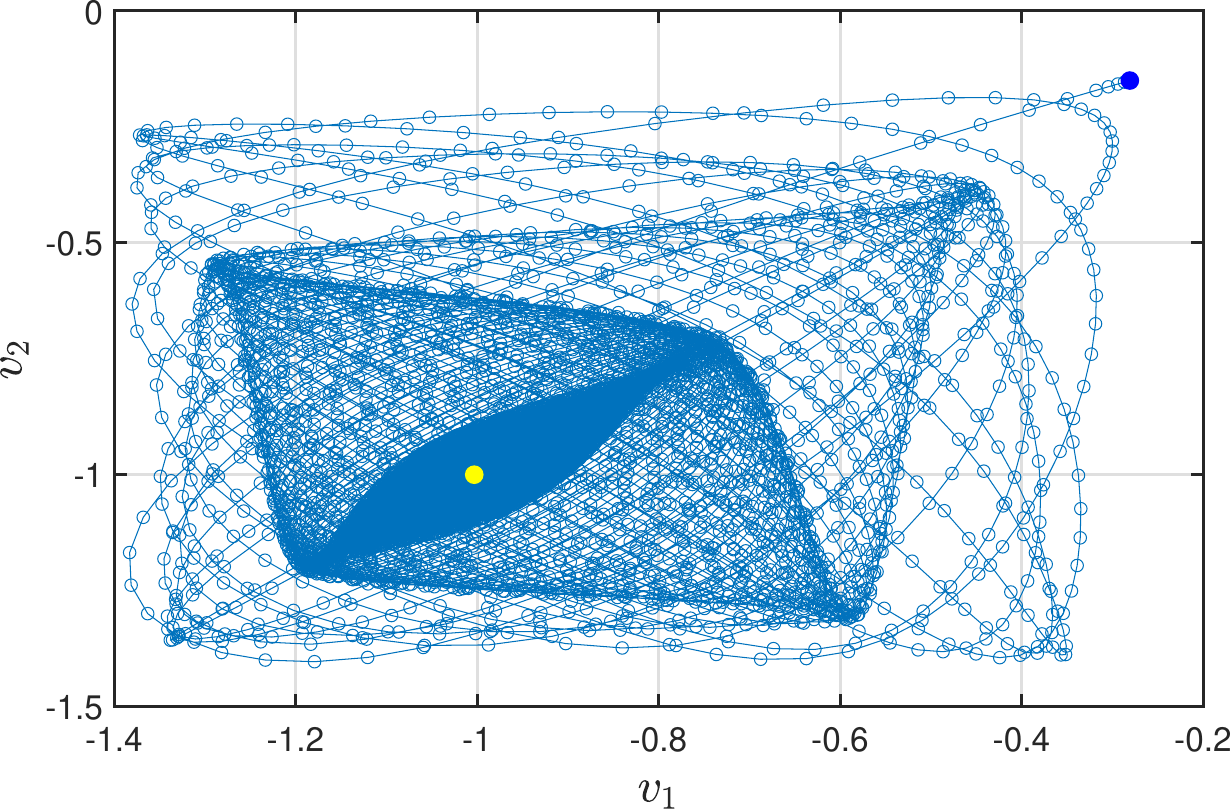} 
	}
	\caption{Impact of $\varepsilon$ to the RK(4,5) scheme.}
	\label{fig:impactofeps_RK45}
\end{figure*}

For a better understanding on the issue of RK(4,5), Figure \ref{fig:impactofeps_RK45} illustrates an example with $(n,d)=(4,3)$ for RK(4,5) where the horizontal and vertical axes are the first and second components of $\V$. We start RK(4,5) from the same initial point (the blue point) with different $\varepsilon\in \{10^{-4}, 10^{-5}, 10^{-7}, 10^{-8}\}$, and observe that the behavior of RK(4,5) scheme varies a lot. A smaller $\varepsilon$ produces more oscillations for the RK(4,5) scheme, thus requires more time steps, but still converge to the same integer solution (the yellow point) with different precision $\delta$.

Note that, we tested on the same dataset in Tables \ref{tab:results1} and \ref{tab:results2}, but obtained different \texttt{obj} values for Houbolt, Lie, RK(4,5) and IPOPT schemes since the initial points are different (generated randomly). Particularly, the \texttt{obj} value of QB-G method will not be affected by the choice of initial point since GUROBI are attempt to find ``global" optimal solutions for the quadratic binary formulation.

\begin{figure*}[htb]
	\centering
	\subfigure[$\gamma/m=10$, number of iterations $221$, solution time $0.03$ seconds]{
		\includegraphics[width=0.46\textwidth]{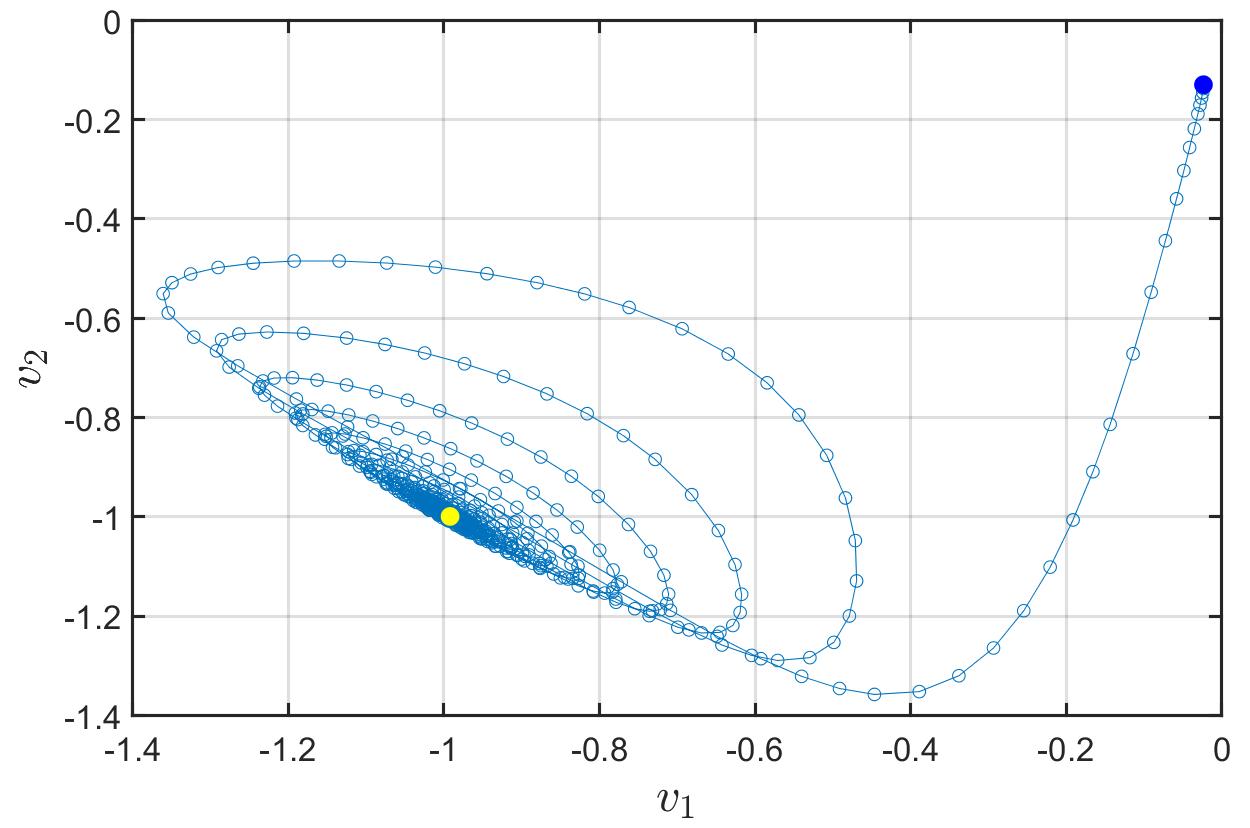}
		\includegraphics[width=0.46\textwidth]{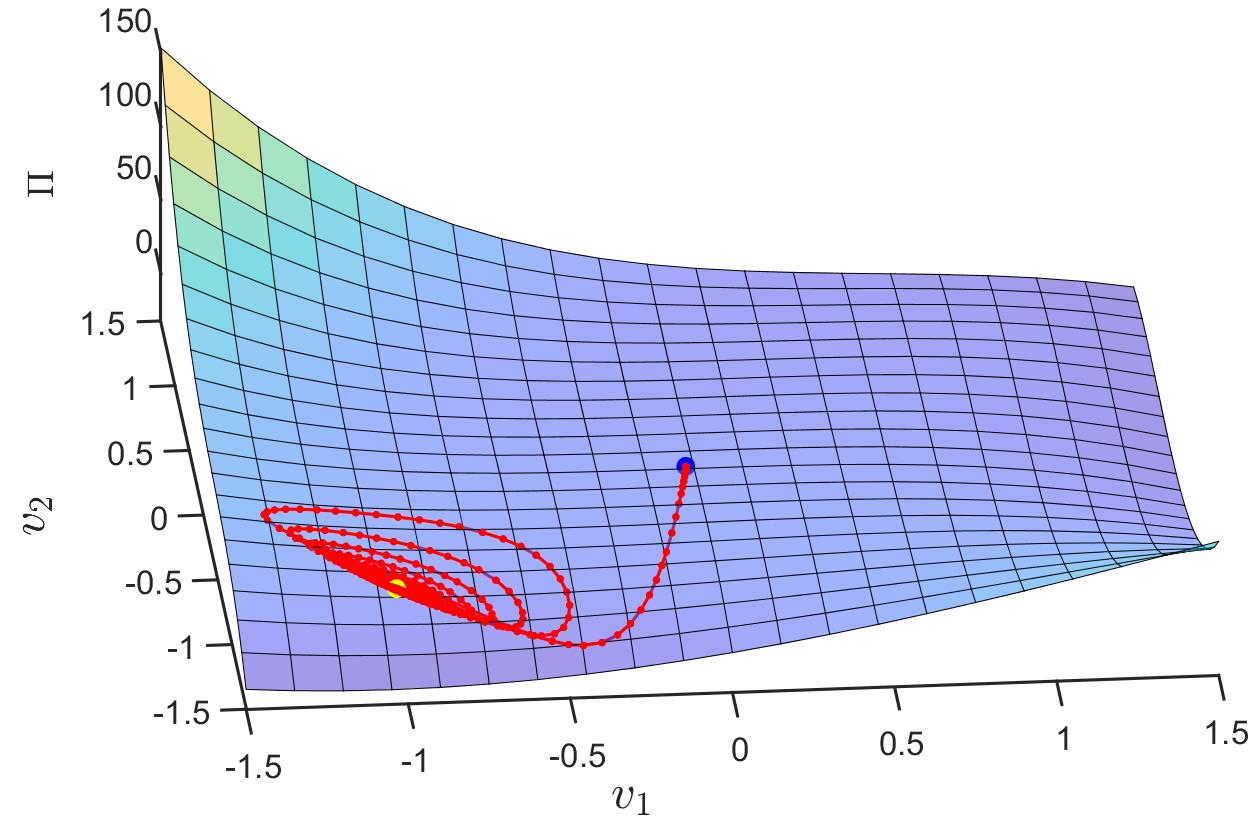}
	}
	\subfigure[$\gamma/m=30$, number of iterations $181$, solution time $0.02$ seconds]{
		\includegraphics[width=0.46\textwidth]{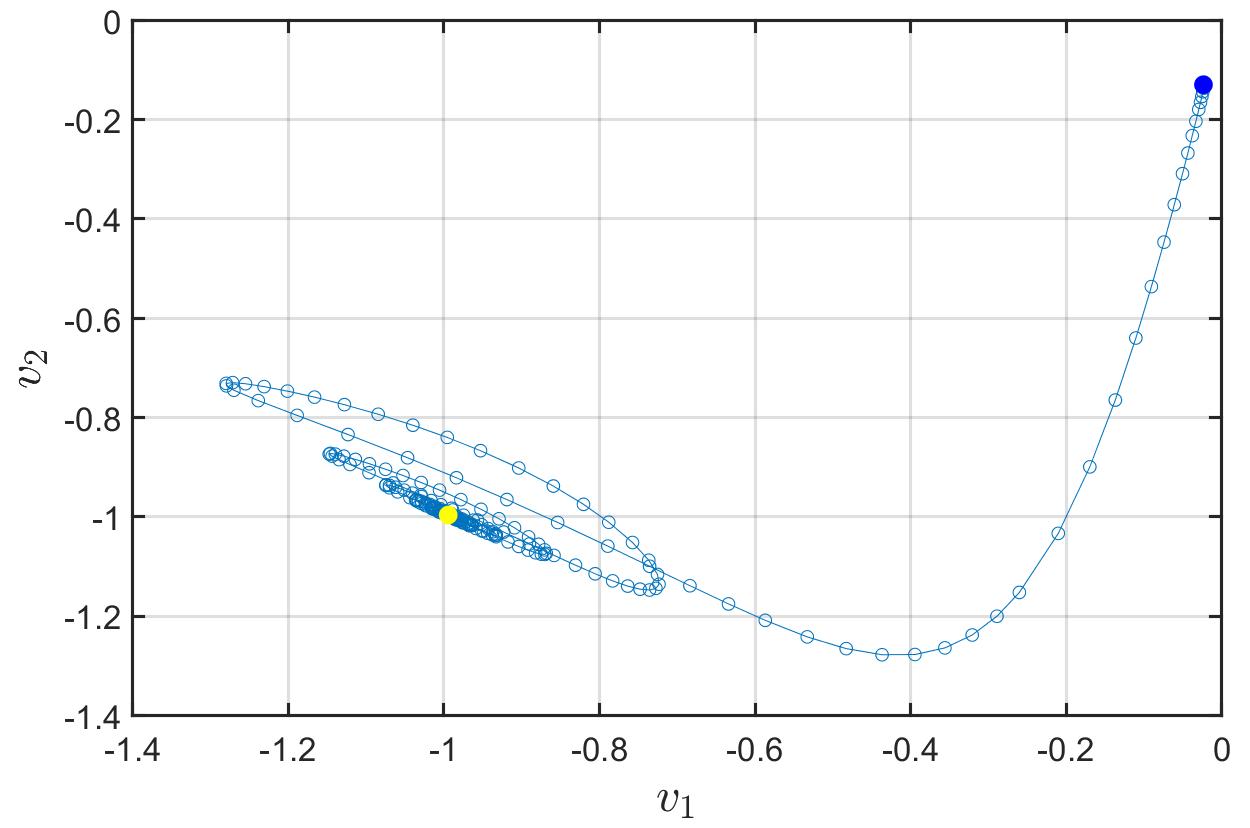}
		\includegraphics[width=0.46\textwidth]{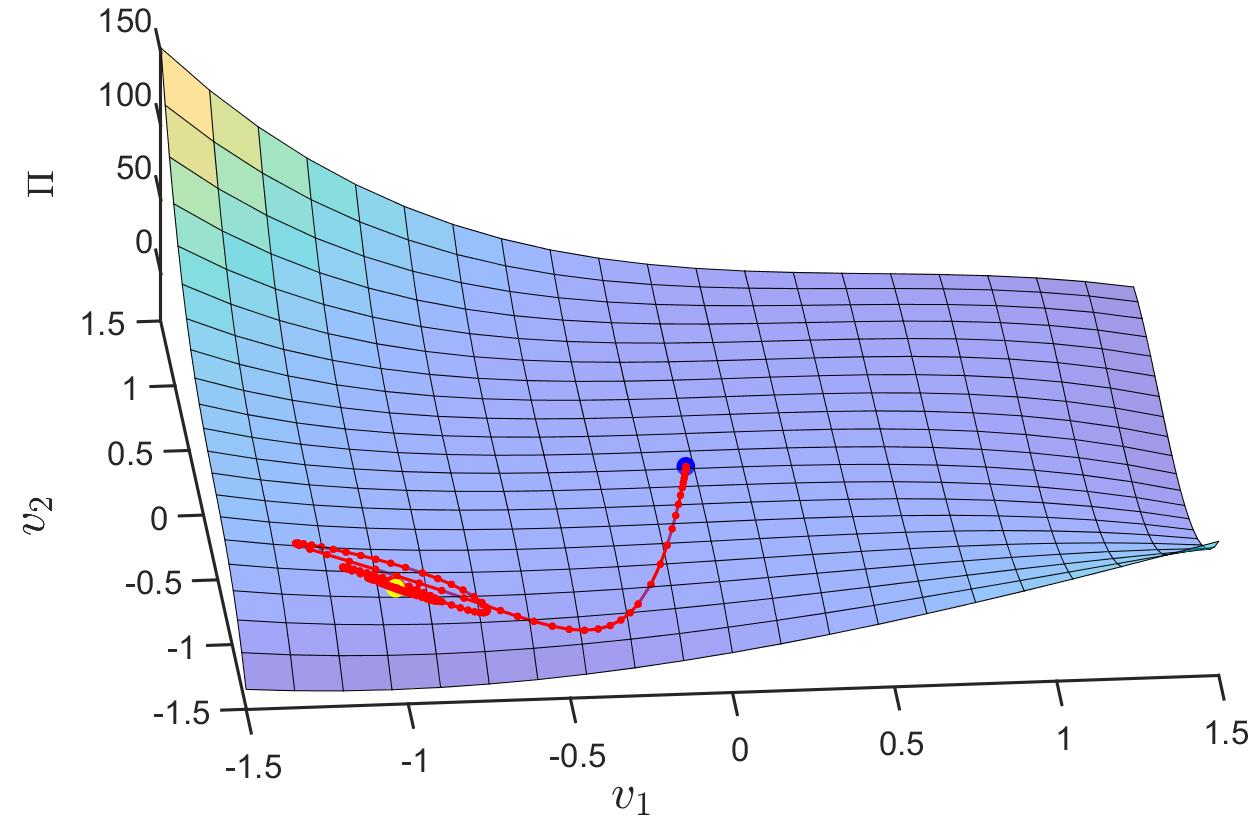}
	}
	\subfigure[$\gamma/m=60$, number of iterations $121$, solution time $0.01$ seconds]{
		\includegraphics[width=0.46\textwidth]{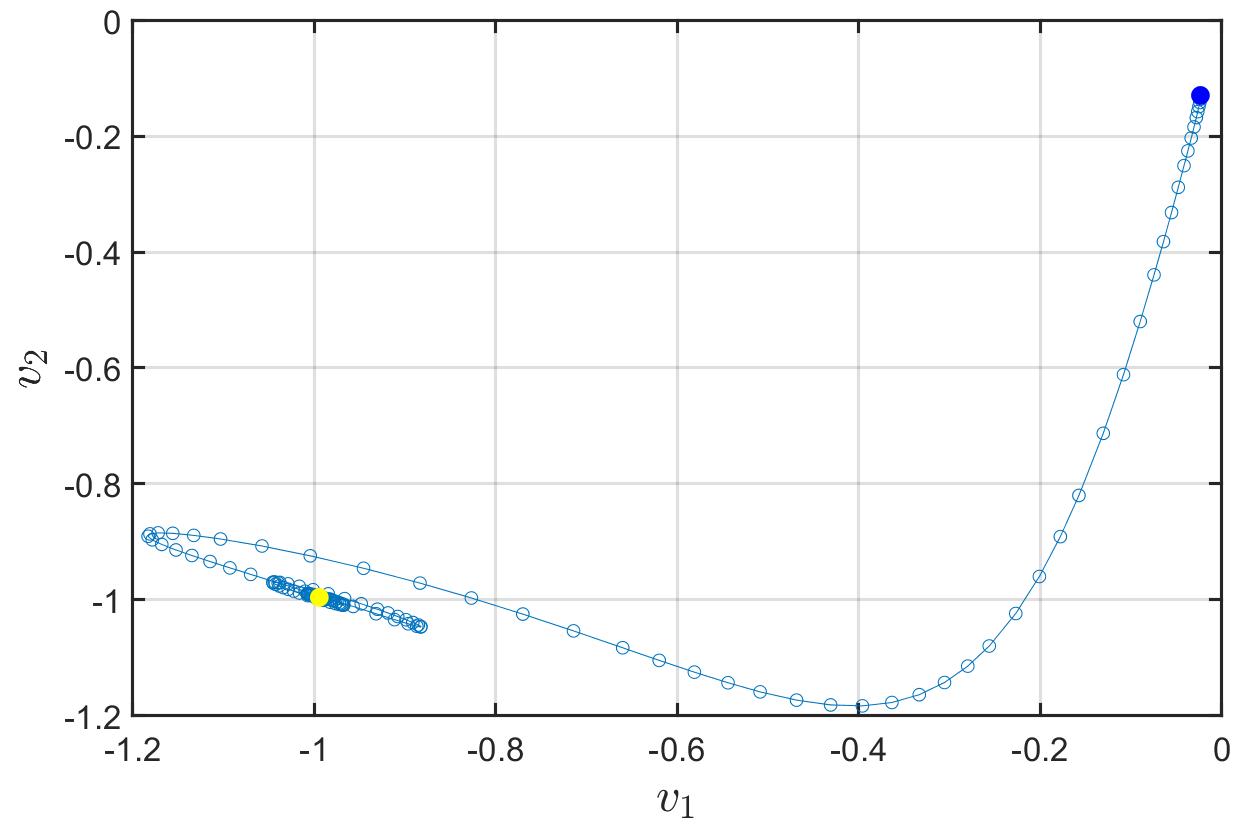}
		\includegraphics[width=0.46\textwidth]{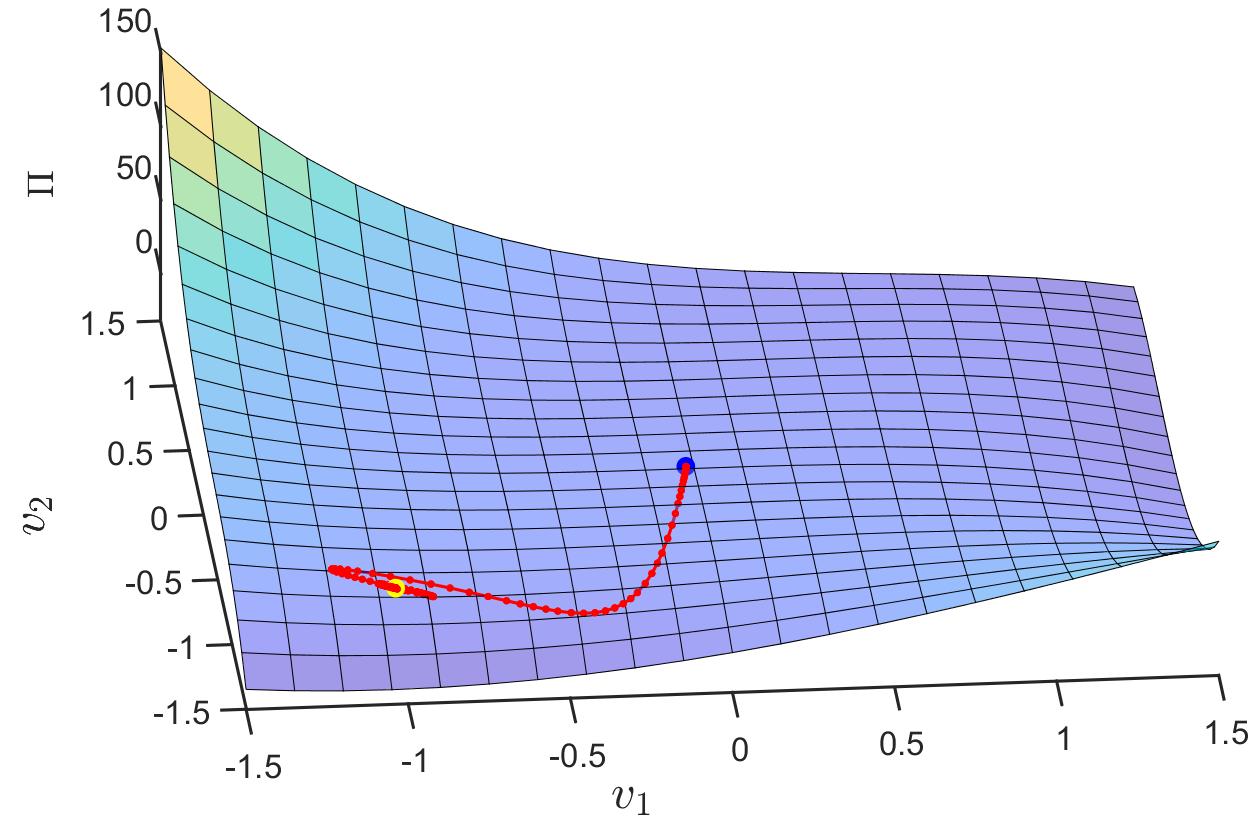}
	}
	\caption{Impact of the damping ratio $\gamma/m$ to the convergence of the RK(4,5) scheme.}
	\label{fig:impactofdampingratio_RK45}
\end{figure*}

\begin{figure*}[htb!]
	\centering
	\subfigure[$\gamma/m=10$, number of iterations $14$, solution time $0.02$ seconds]{
		\includegraphics[width=0.46\textwidth]{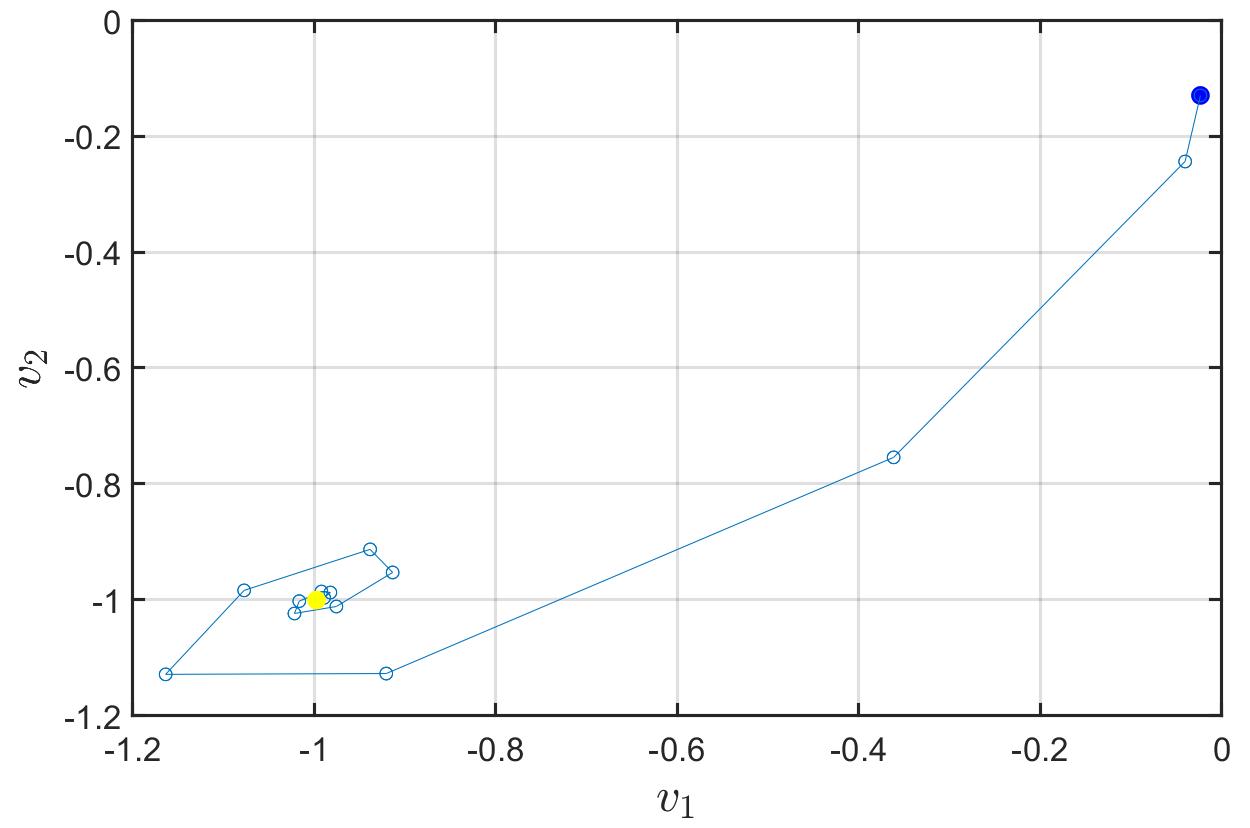}
		\includegraphics[width=0.46\textwidth]{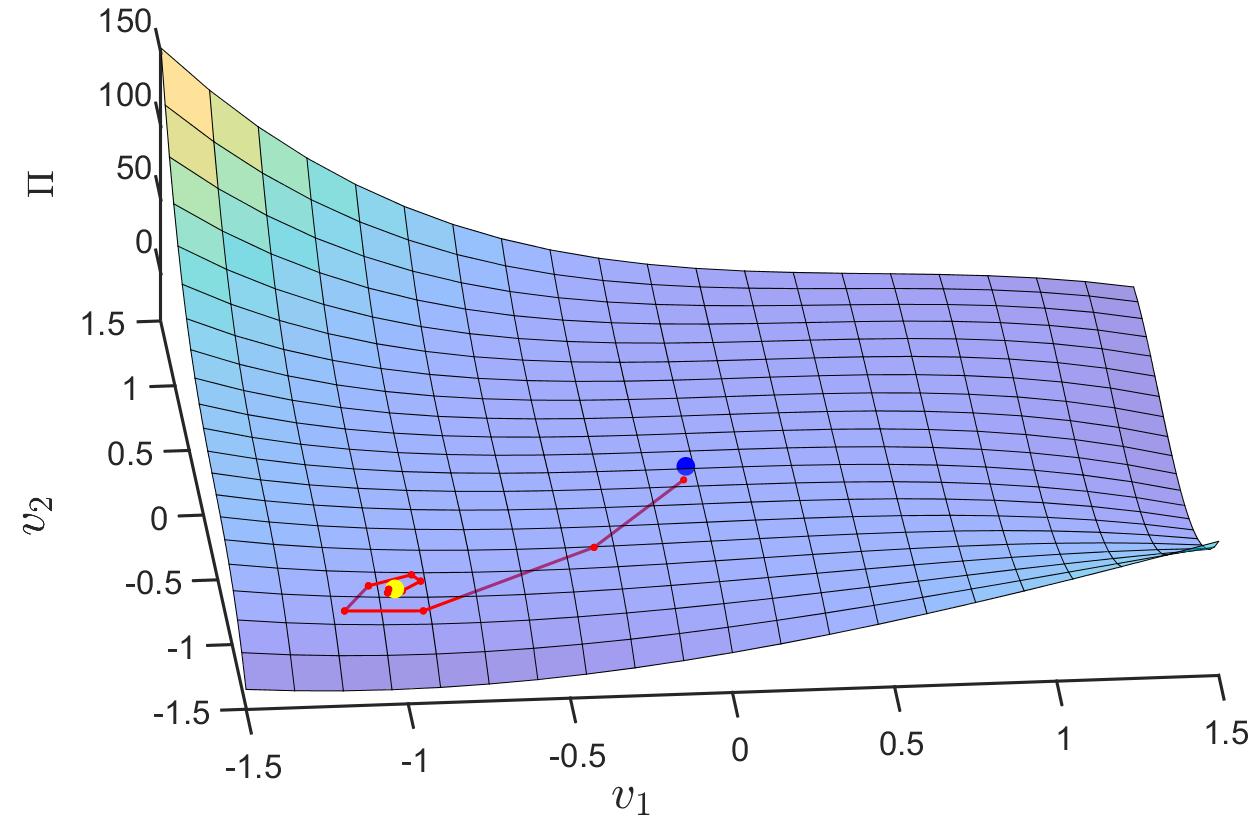}
	}
	\subfigure[$\gamma/m=30$, number of iterations $13$, solution time $0.01$ seconds]{
		\includegraphics[width=0.46\textwidth]{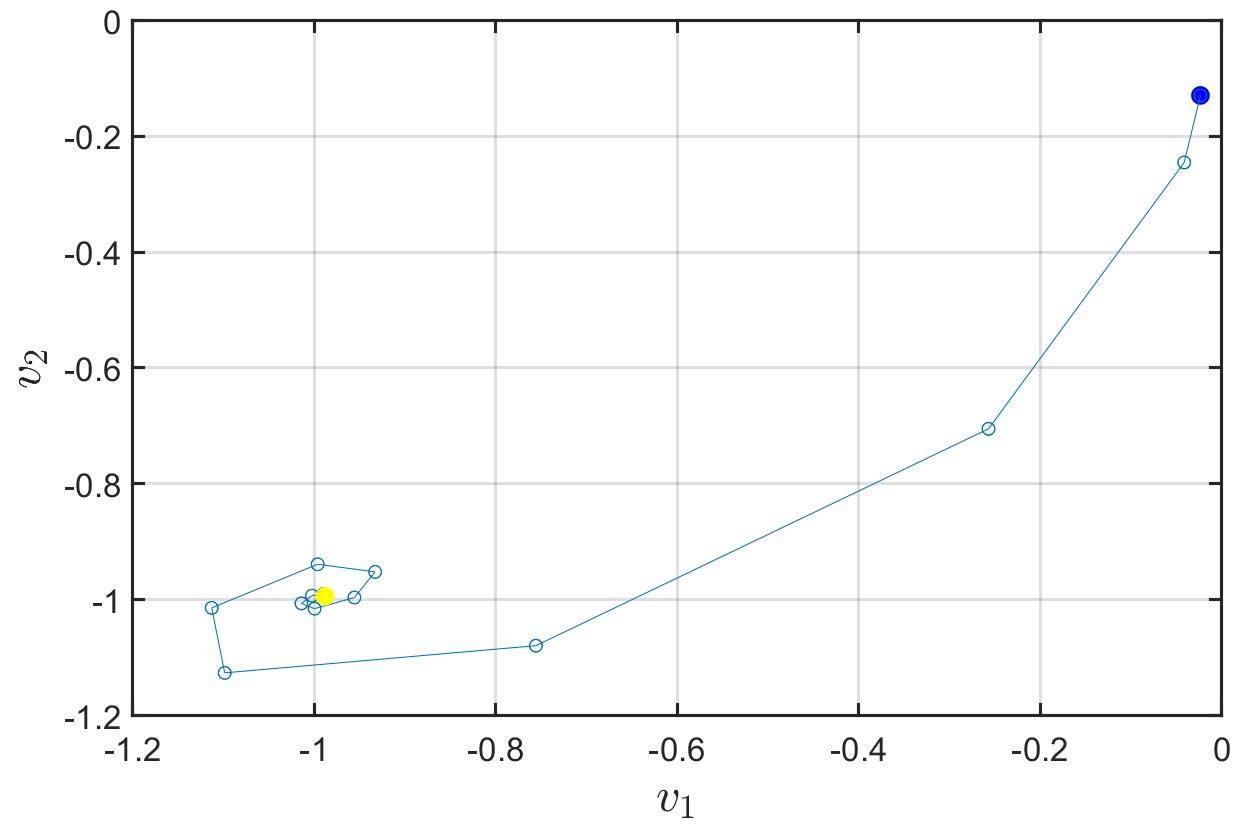}
		\includegraphics[width=0.46\textwidth]{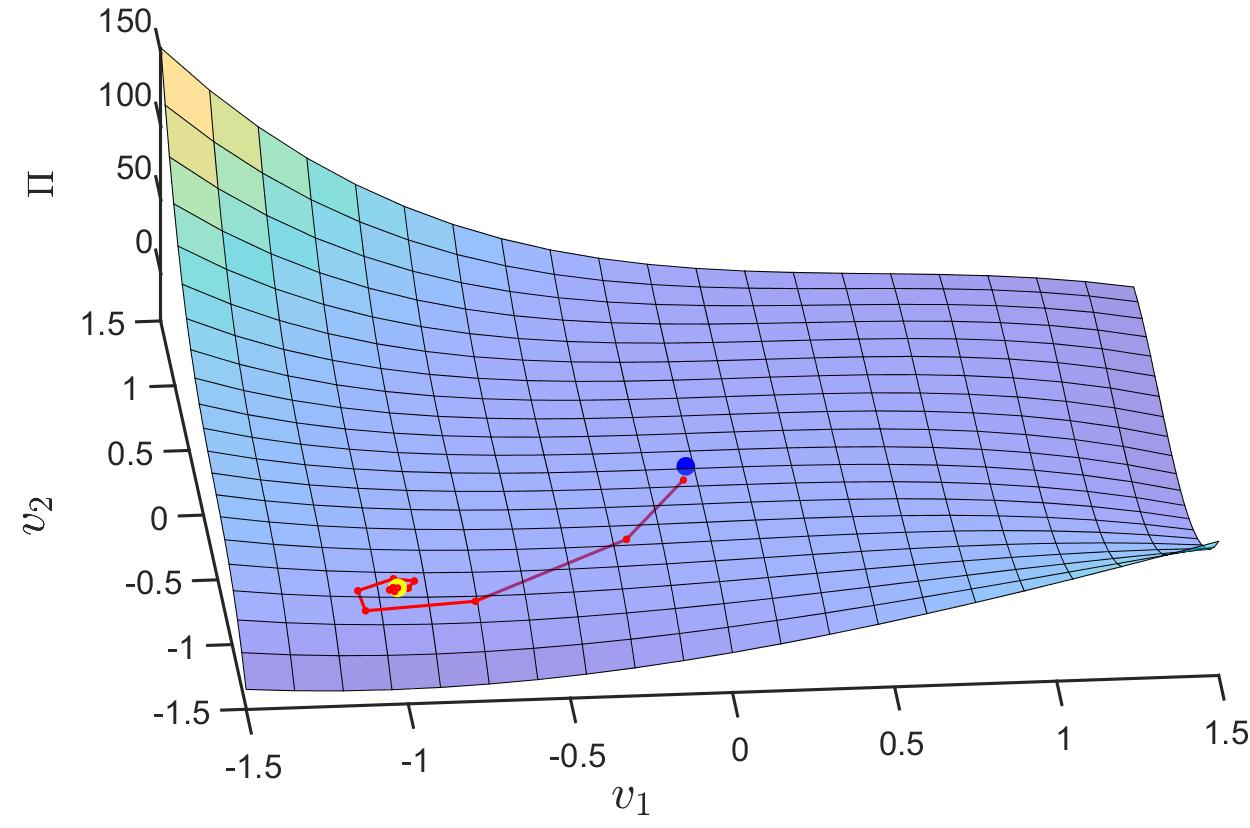}
	}
	\subfigure[$\gamma/m=60$, number of iterations $13$, solution time $0.01$ seconds]{
		\includegraphics[width=0.46\textwidth]{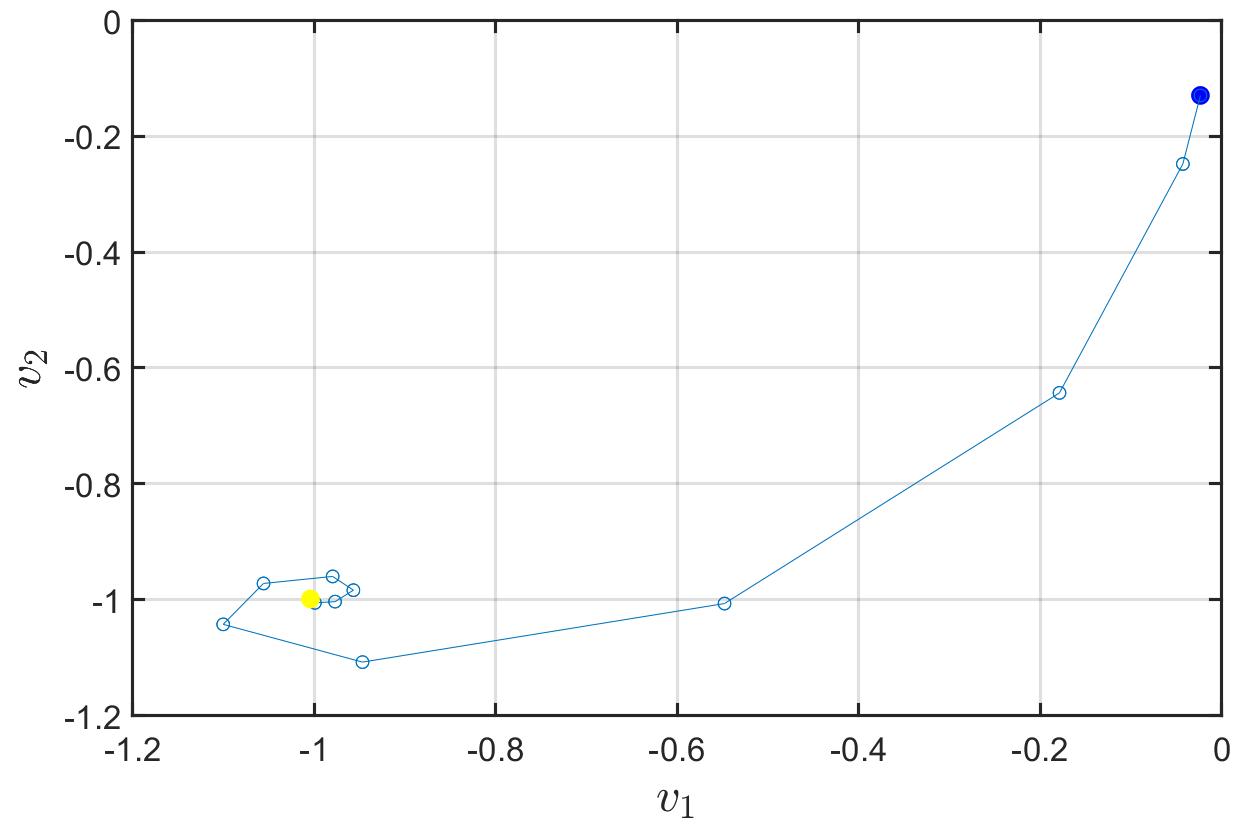}
		\includegraphics[width=0.46\textwidth]{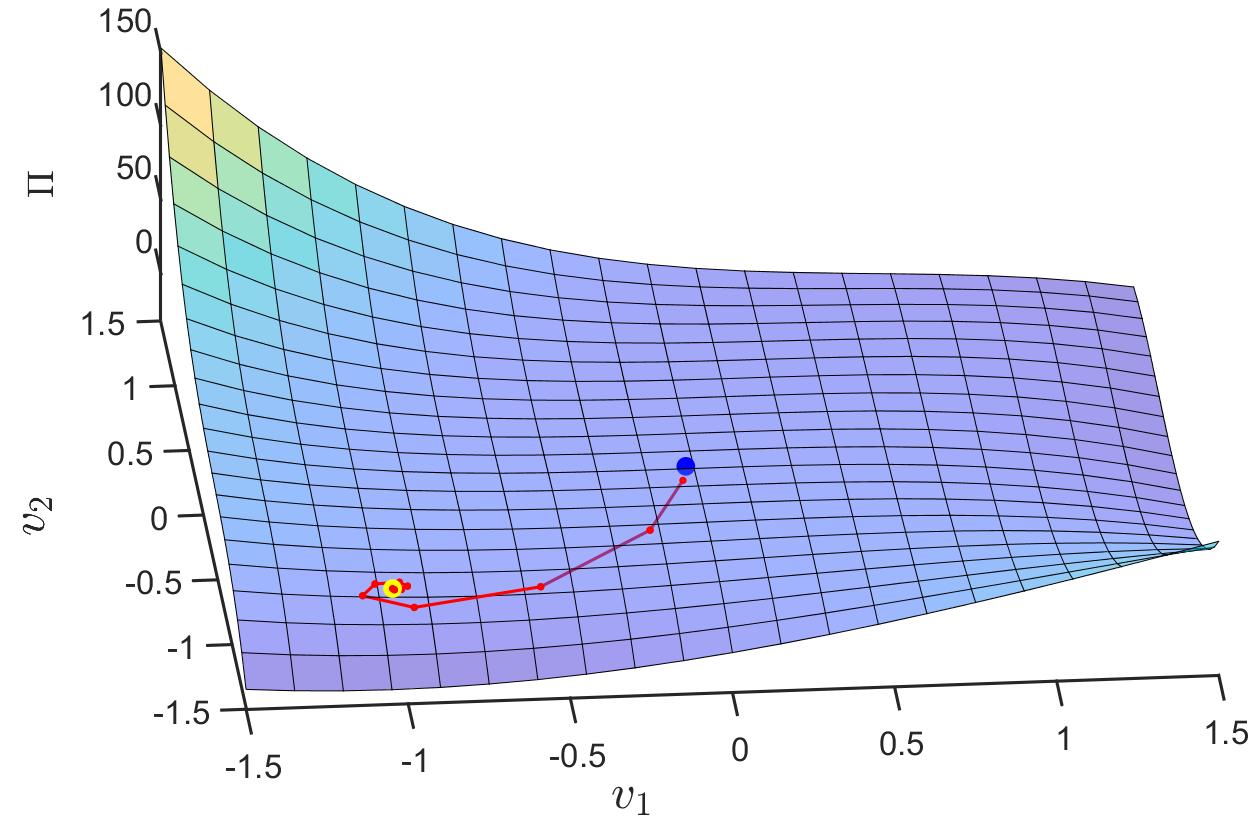}
	}
	\caption{Impact of the damping ratio $\gamma/m$ to the convergence of the Houbolt scheme.}
	\label{fig:impactofdampingratio_Houbolt}
\end{figure*}

\paragraph{\textbf{Impact of $\gamma/m$}}\label{para:impactofdampingratio}
The two parameters $\gamma$ and $m$ are only needed in the Houbolt and RK(4,5) schemes. Figures \ref{fig:impactofdampingratio_RK45} and \ref{fig:impactofdampingratio_Houbolt} illustrate the impact of the ratio $\gamma/m$ to them, which demonstrate the most different behaviors between the implicit and the explicit schemes. These figures are obtained with $n=2$, $d=4$, $\varepsilon=10^{-4}$ and with different ratios $\gamma/m \in\{10,30,60\}$. The initial point (generated randomly) is the same for these tests, and $t\in [0,1]$ for the RK(4,5) scheme. We can observe that both the RK(4,5) and Houbolt schemes converge in very different ways towards the optimal solution $(-1,-1)$, and the convergence rates depend on the ratio $\gamma/m$. Generally speaking, a bigger ratio $\gamma/m$ leads to a faster convergence. It seems that the Houbolt scheme converges more directly and faster to the solution $(-1,-1)$, while the RK(4,5) scheme converges in a more oscillatory way, thus less efficiently than the Houbolt scheme. 

\begin{figure*}[htb!]
	\centering
	\subfigure[$c=0$]{
		\includegraphics[width=0.23\textwidth]{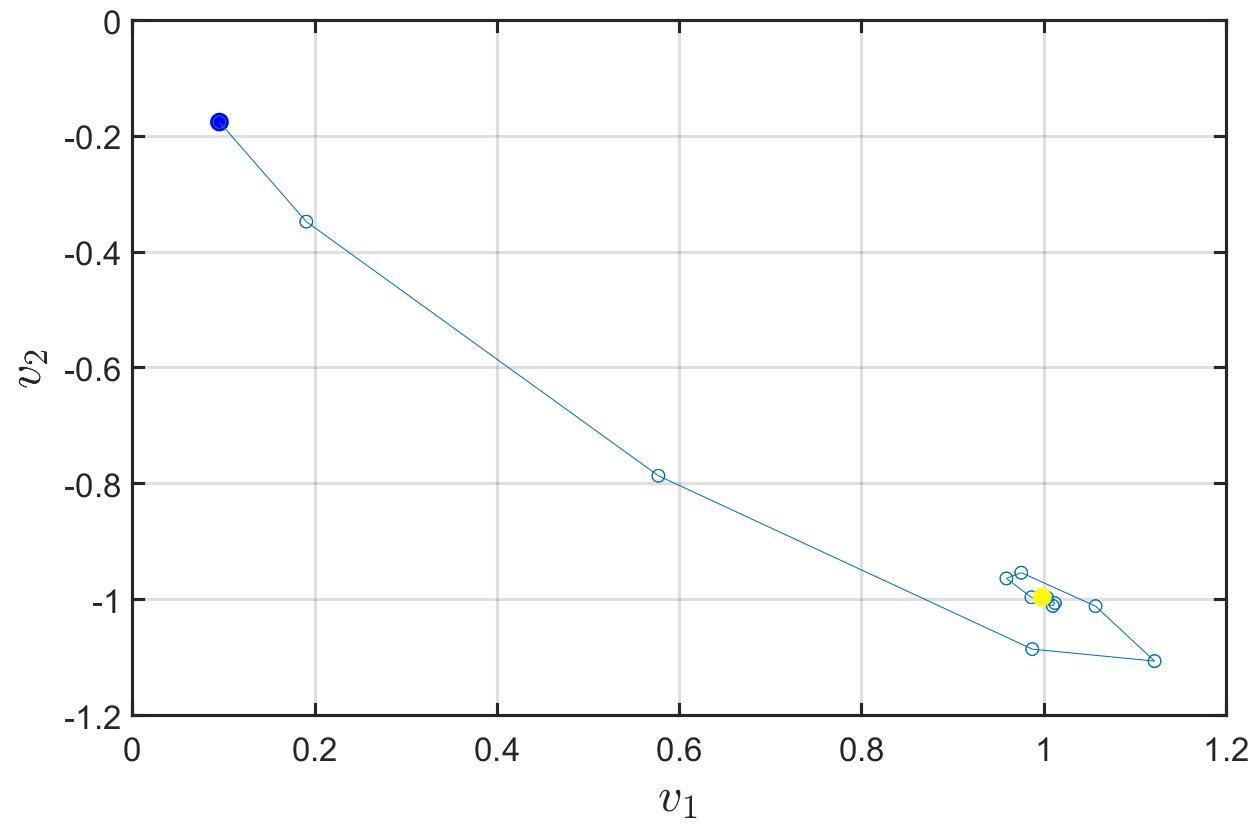}
		\includegraphics[width=0.23\textwidth]{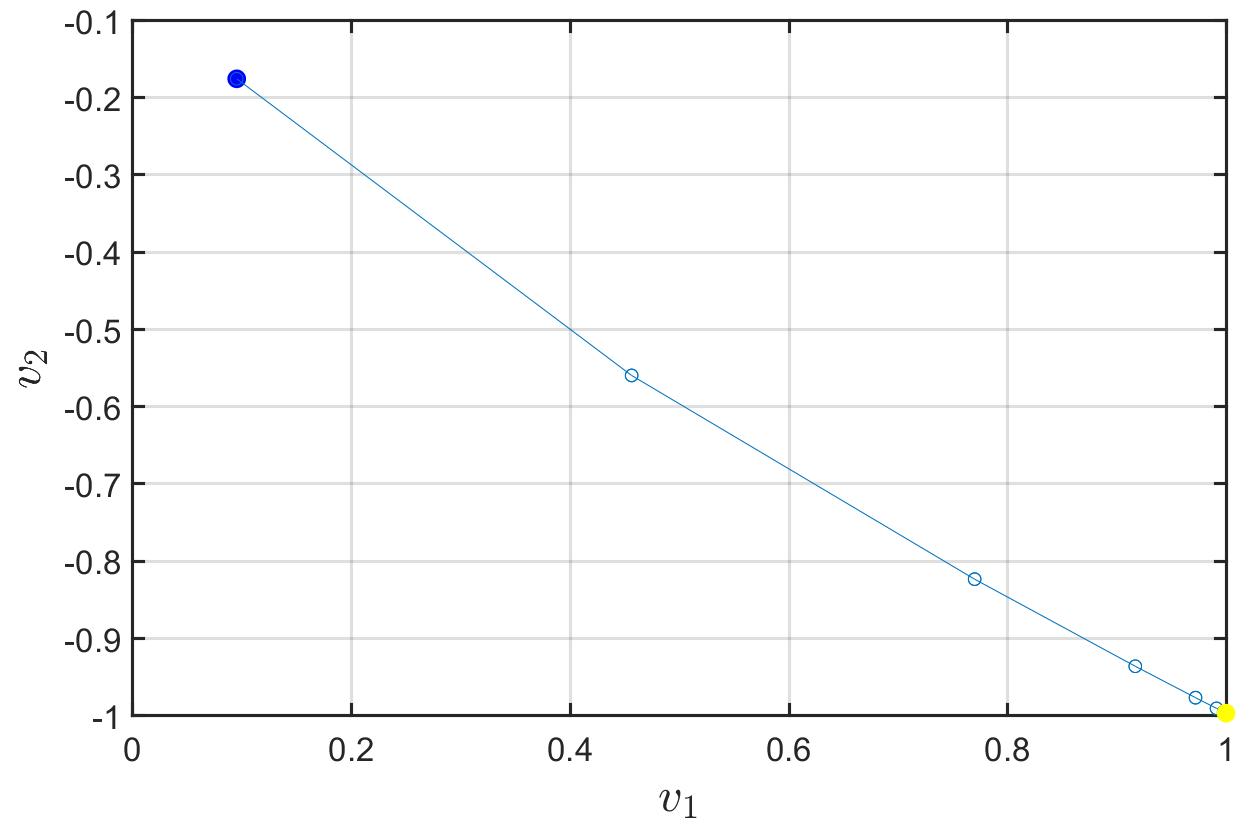}
		\includegraphics[width=0.23\textwidth]{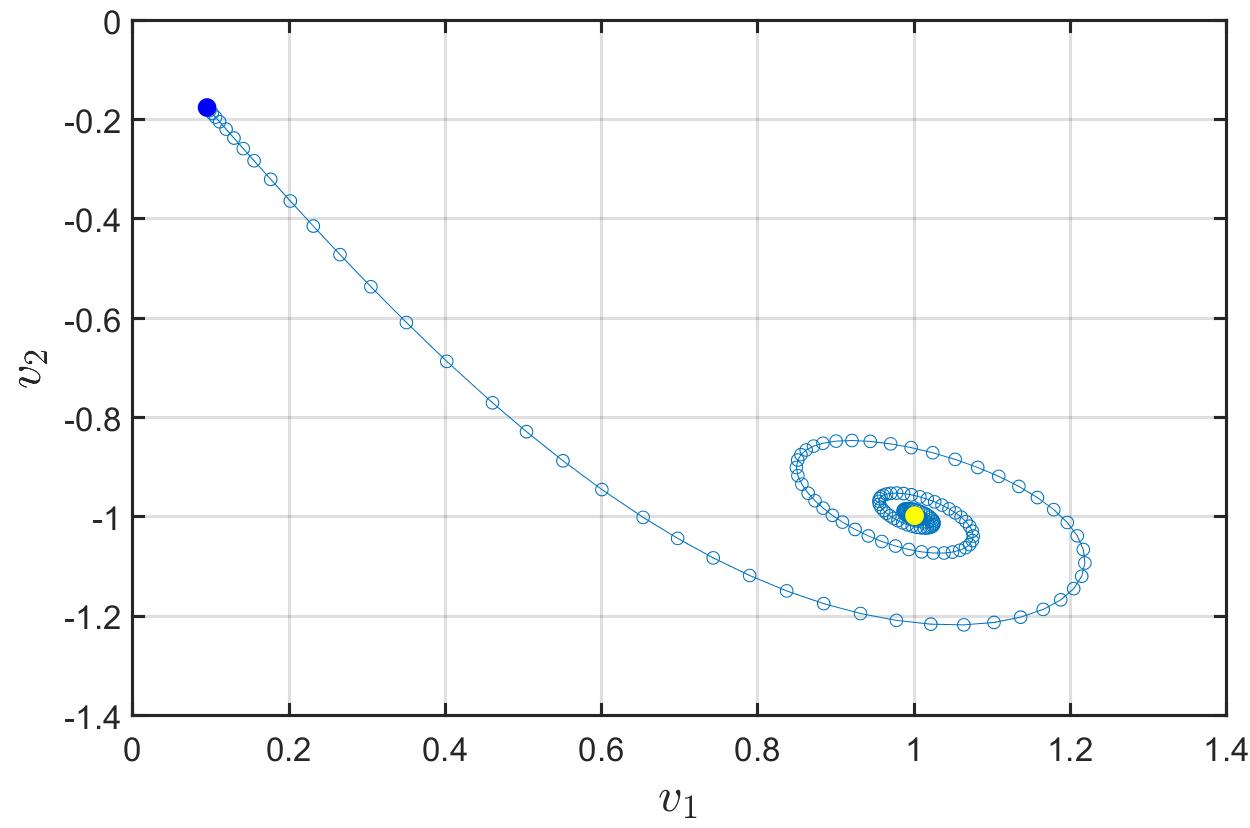}			
		\includegraphics[width=0.23\textwidth]{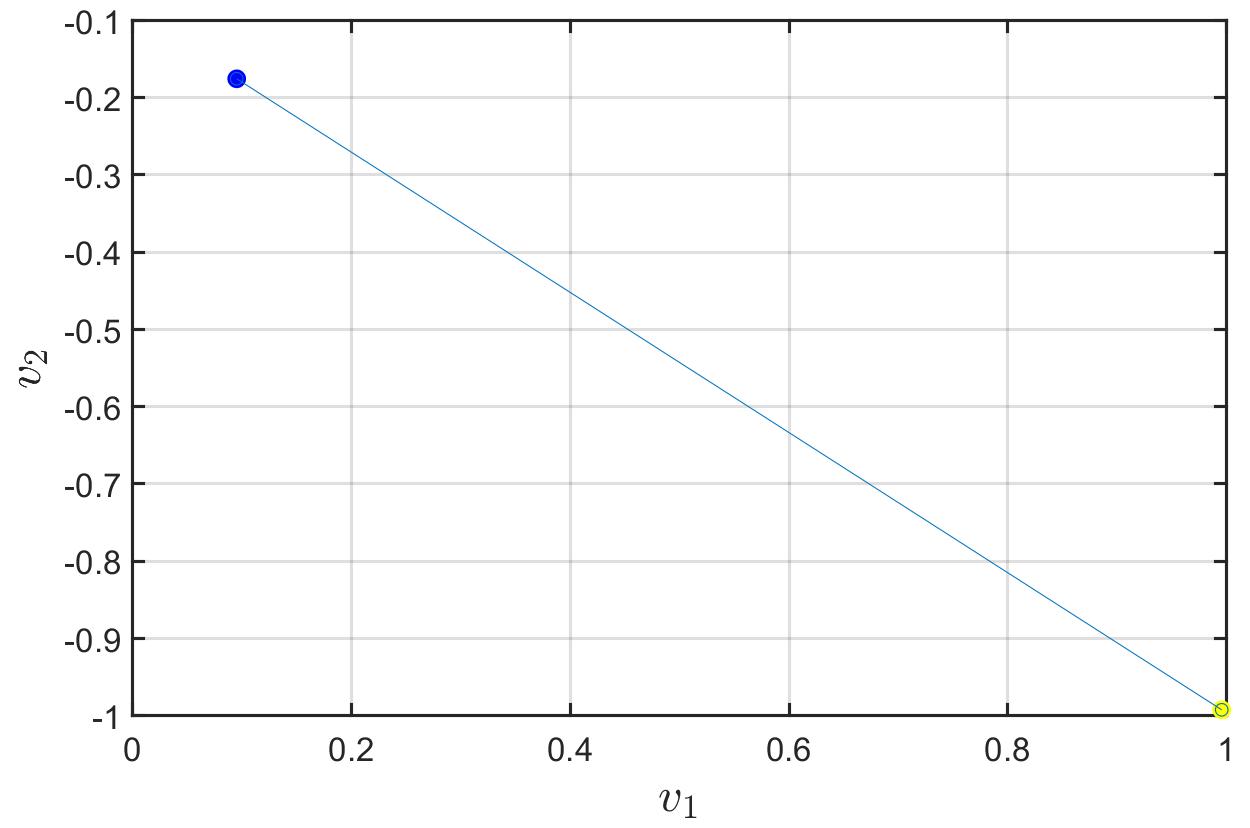}
	}
	\subfigure[$c=100$]{
		\includegraphics[width=0.23\textwidth]{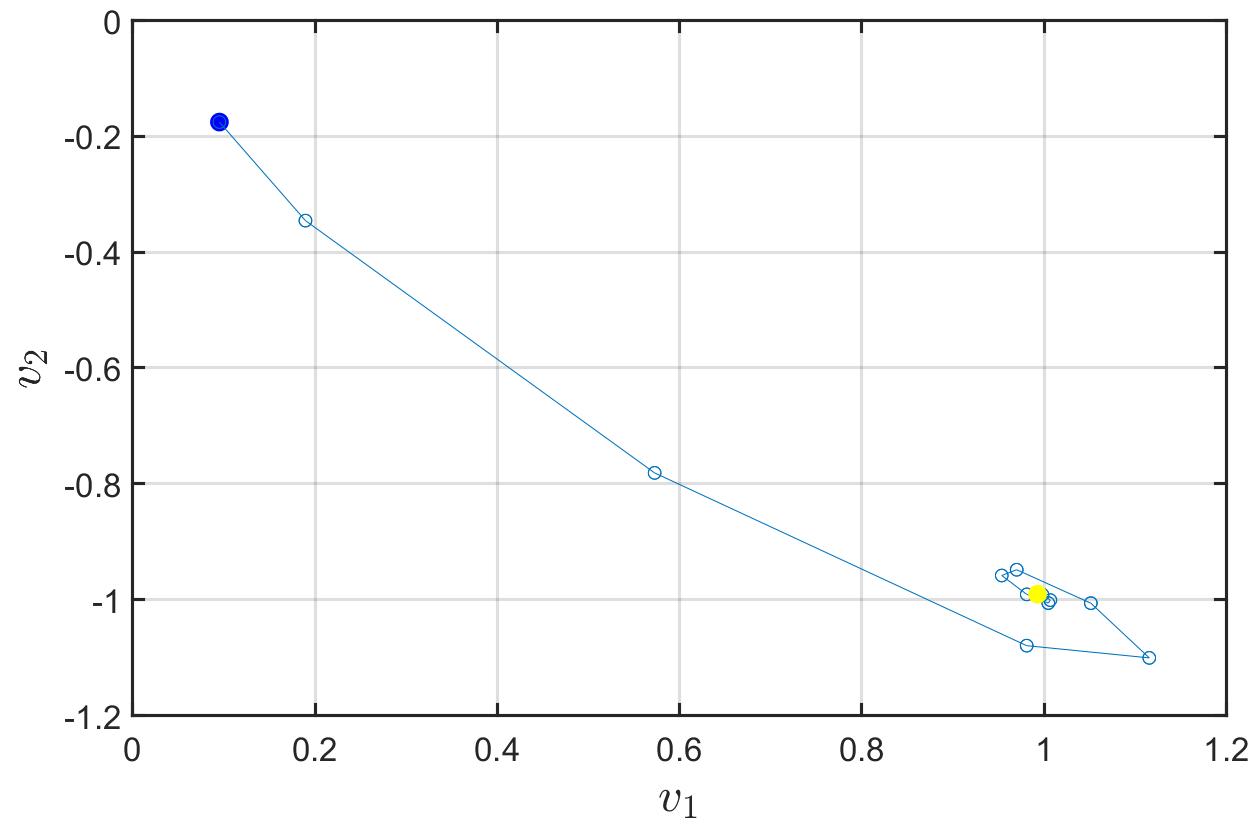}
		\includegraphics[width=0.23\textwidth]{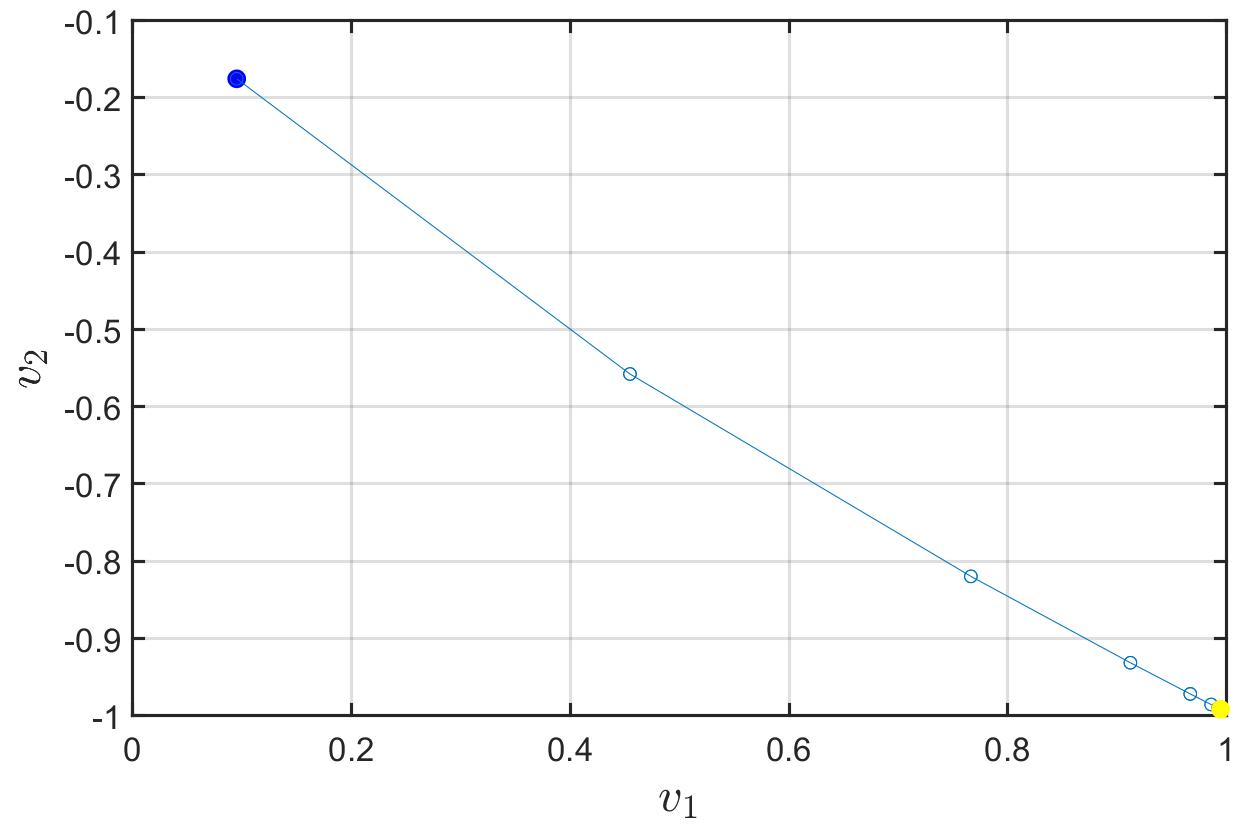}
		\includegraphics[width=0.23\textwidth]{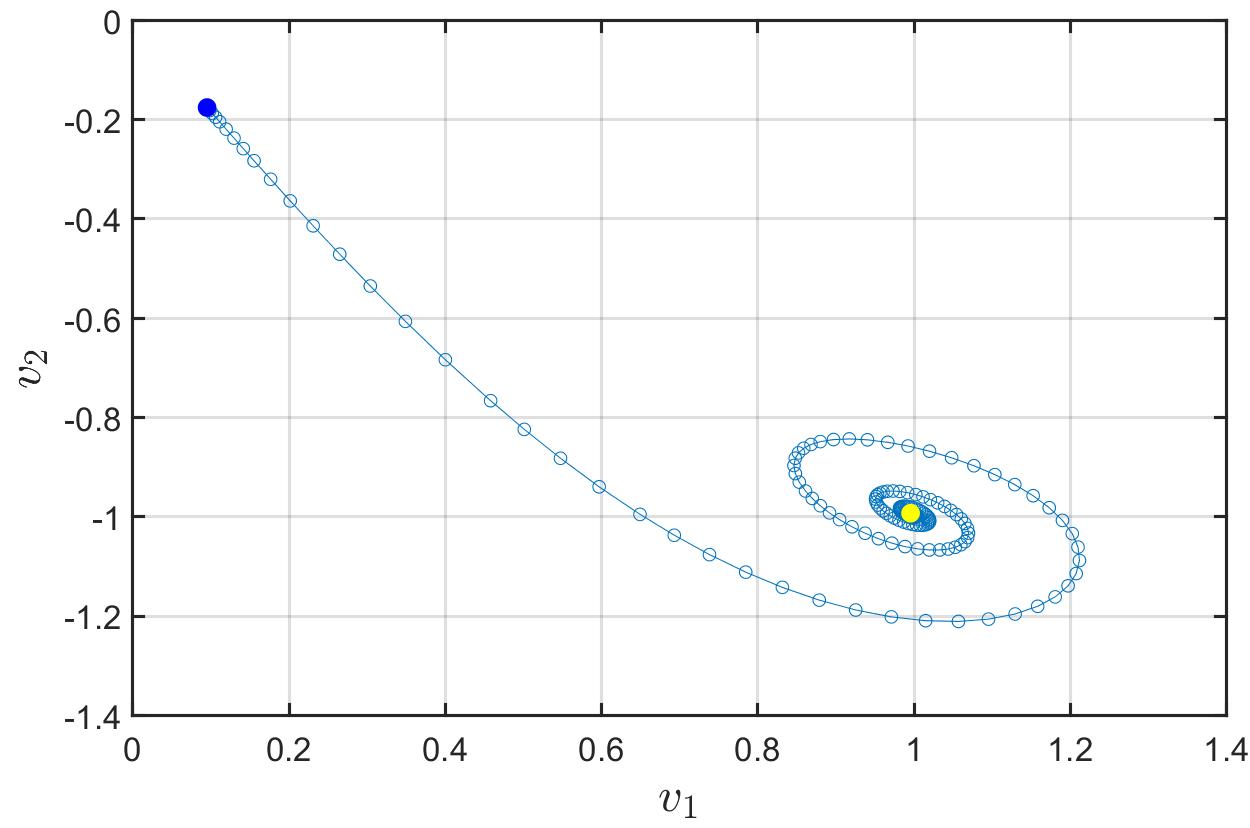}
		\includegraphics[width=0.23\textwidth]{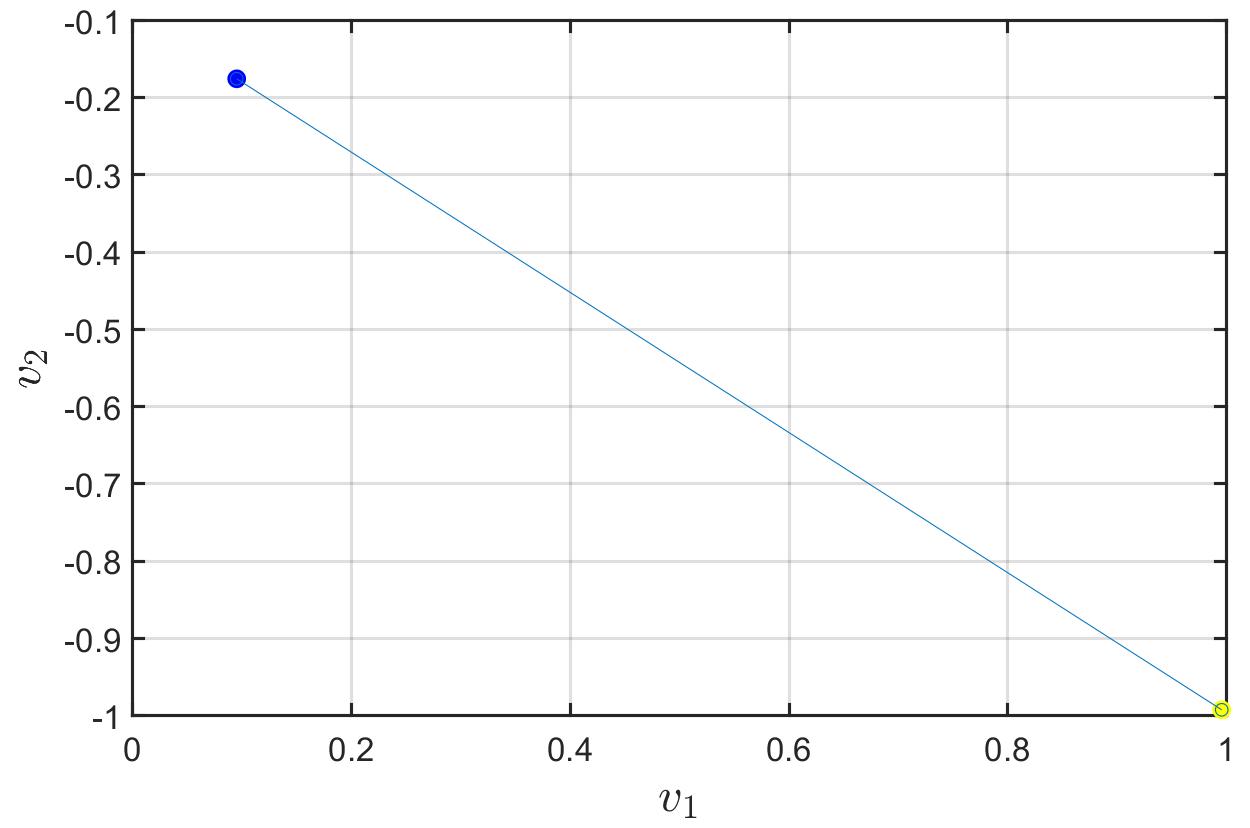}			
	}
	\subfigure[$c=10^{5}$]{
		\label{subfig:instabilityissue}
		\includegraphics[width=0.23\textwidth]{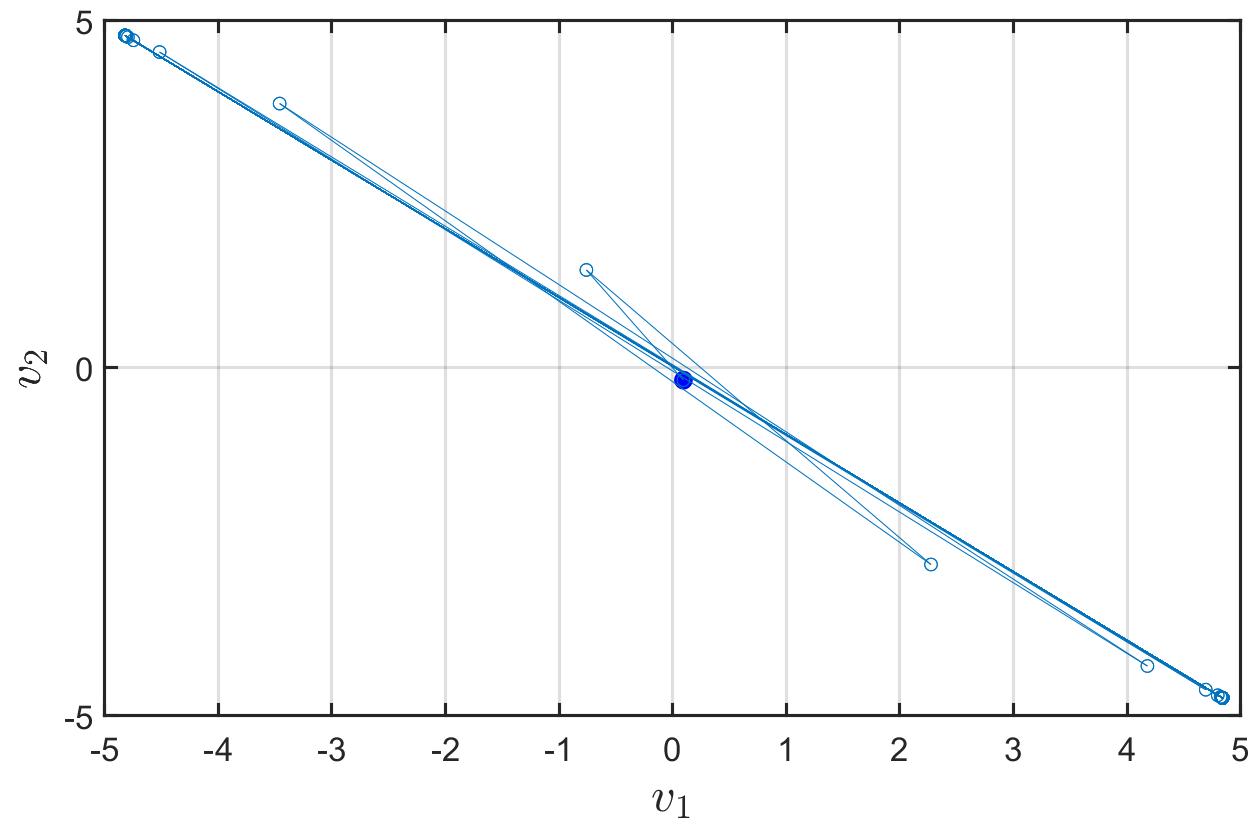}
		\includegraphics[width=0.23\textwidth]{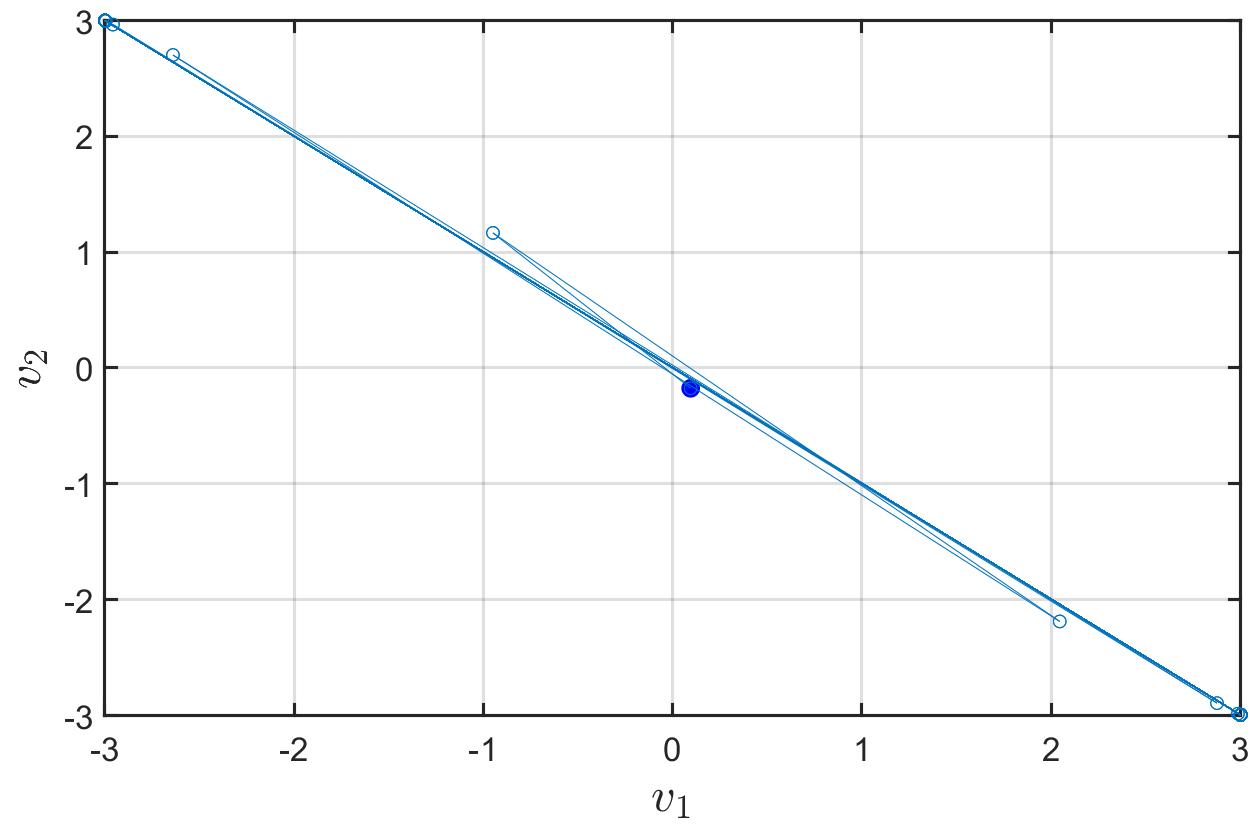}
		\includegraphics[width=0.23\textwidth]{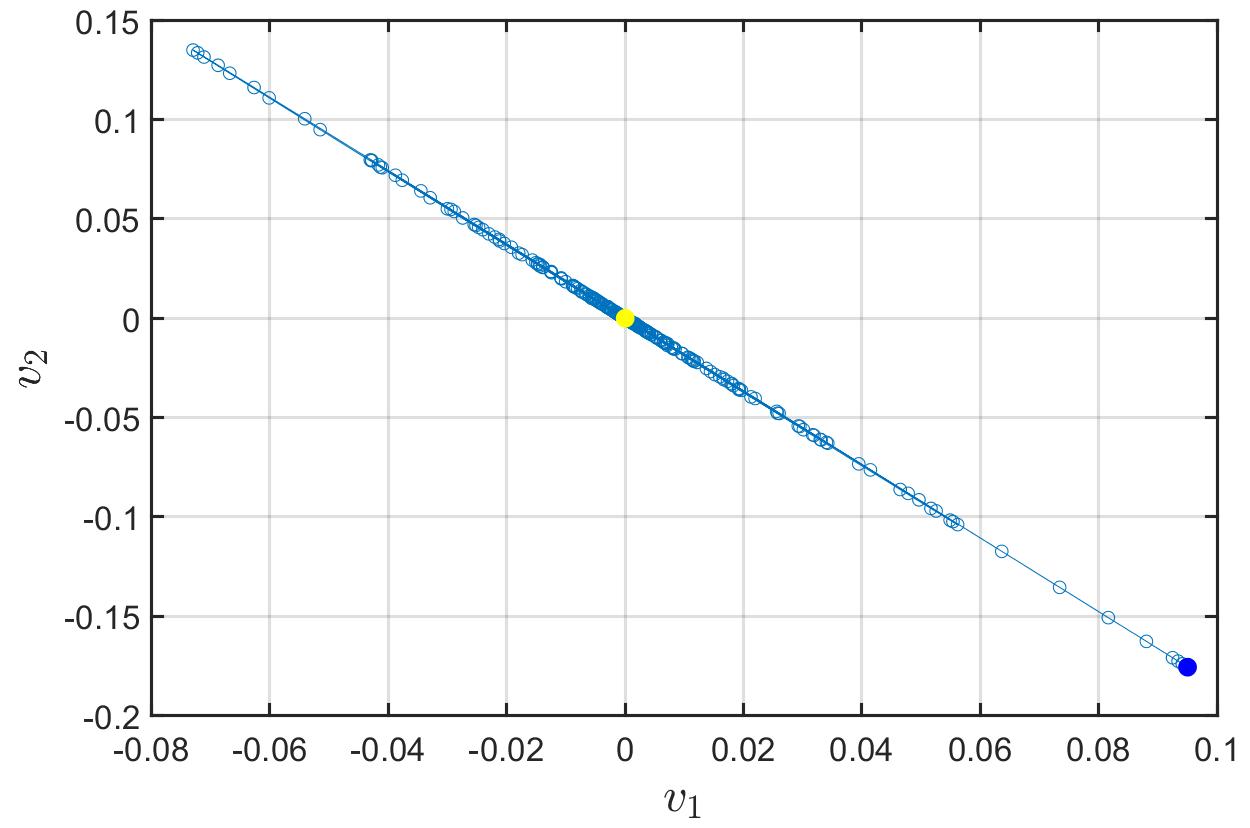}
		\includegraphics[width=0.23\textwidth]{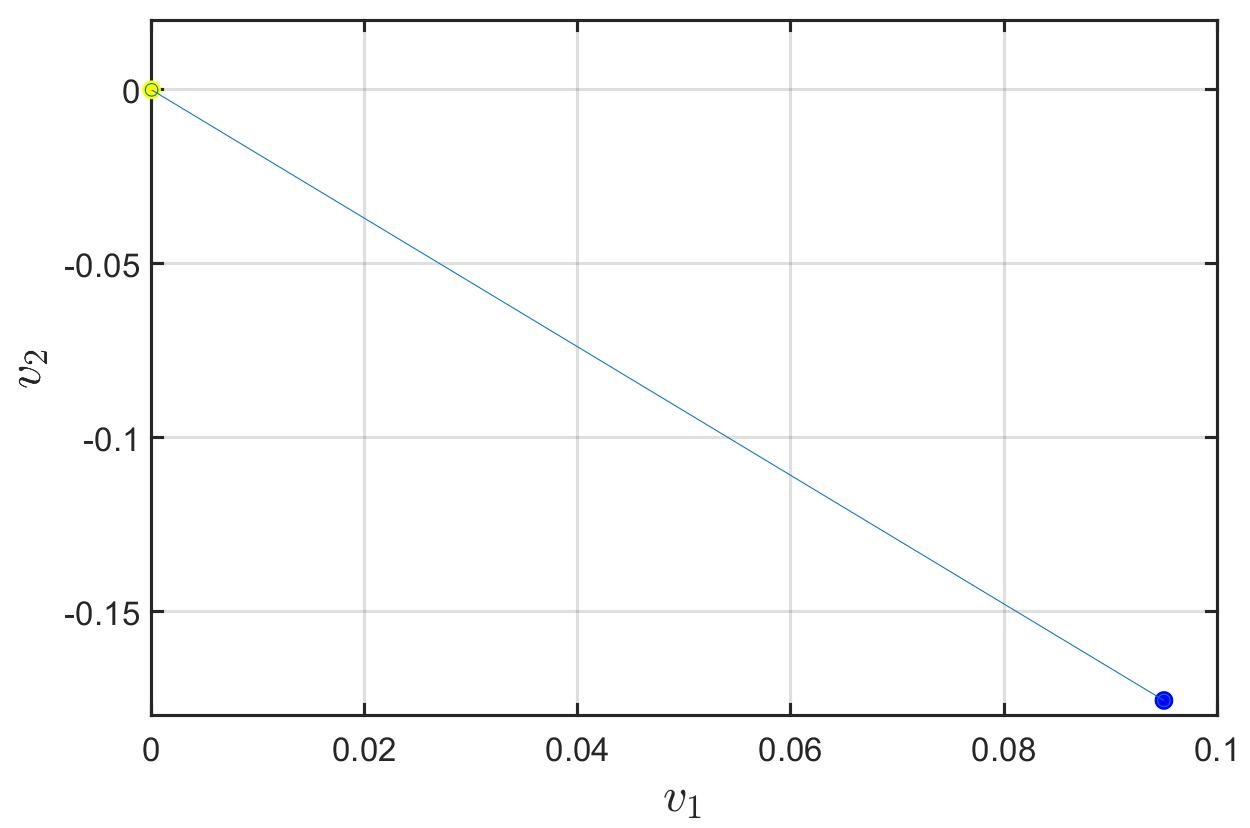}	
	}
	\caption{Impact of $c$ on the convergence of the Houbolt (left), Lie (middle-left) and RK(4,5) (middle-right) and IPOPT (right) schemes.}
	\label{fig:impactofc}
\end{figure*}

\paragraph{\textbf{Impact of $c$}}\label{param:impactofc}
The parameter $c$ is used in all schemes except the BQ-G method. An example is illustrated in Figure \ref{fig:impactofc} where the four schemes (Lie, Houbolt, RK(4,5) and IPOPT) are compared with different values of $c\in\{0,10^2,10^5\}$ and with $n=10$, $d=4$, $\varepsilon=10^{-4}$, $\gamma/m=50$ for the Houbolt and RK(4,5) schemes, $\tau=\sqrt{2m\varepsilon} \approx 0.0141$ for the Houbolt scheme, $\tau = \min\{\varepsilon/(1-c\varepsilon), 0.1\}$ for the Lie scheme, and $t\in [0,0.5]$ for the RK(4,5) scheme.

Numerical tests show that $c$ has very small impact to these methods when $c$ is not too large (even if the functional $\V\mapsto c\|\V\|_2^2/2 + \Pi(\V)$ is nonconvex over $\Bor$ with $r\geq \sqrt{n}$), since there are almost no difference for each method with $c=0$ and $c=100$. However, $c$ cannot be too large, there are several reasons: (i) the term $\Pi(\V)$ will be negligible with respect to $c\|\V\|_2^2/2$, which will make the functional $\V\mapsto c\|\V\|_2^2/2 + \Pi(\V)$ close to the convex quadratic functional $\V\mapsto c\|\V\|_2^2/2$ over $\Bor$. Then the obtained numerical solution for problem \cref{prob:approx_int_cv_opt} will be close to $\oRn$ not $\{\pm 1\}^n$. (ii) when $c$ is too large, the convergence of the Lie scheme will be destroyed since \eqref{eq:boundsoftauandeps} (i.e., $\varepsilon \leq 1/c$ and $0< \tau\leq \varepsilon/(1-c\varepsilon)$) cannot been verified for fixed $\varepsilon$. These issues are also observed in Subfigure \ref{subfig:instabilityissue}, and we can see that when $c=10^5$ (too large), using $\tau = \min\{\varepsilon/(1-c\varepsilon),0.1\}$ for the Lie scheme yields negative time step $\tau=-1.11\times 10^{-5}$ and $\varepsilon = 10^{-4} > 1/c = 10^{-5}$ which destroyed the convergence of the Lie scheme; while the RK(4,5) and IPOPT schemes converge to a point (in yellow) close to $\oRn$. Therefore, the parameter $c$ does not seems helpful to improve the quality of the numerical result and can lead to divergence issue in some cases, for this reason, we suggest fixing $c=0$ in practice and in our numeral tests subsequently.

\paragraph{\textbf{Impact of $\tau$}}
The parameter $\tau$ is only used in the Houbolt and Lie schemes. The convergence is guaranteed when $\tau$ verifying $0<\tau\leq \sqrt{2m\varepsilon}$ for the Houbolt scheme and $c+\frac{1}{\tau}\geq \frac{1}{\varepsilon}$ for the Lie scheme. In previous numerical experiments, we have fixed $\tau$ using formulations \eqref{eq:tau_Houbolt} and \eqref{eq:tau0_Lie}. Numerical results in Tables \ref{tab:results1} and \ref{tab:results2} demonstrate that the Lie and Houbolt schemes perform well enough using fixed $\tau$ for problems in small-scale dataset.

\subsubsection{Numerical Results on Large-scale Randomly Generated Dataset}
The large-scale dataset is randomly generated with $10\leq n\leq 20$ and $5\leq d\leq 6$, where the sparsity of the polynomials is set to $0.8$, which is a situation likely to arise in hard practical applications. Table \ref{tab:results5} illustrates numerical results on large-scale dataset with $\varepsilon = 10^{-6}$, $m=1$, $\gamma=300$, $c=0$, $t\in [0,0.3]$, and with fixed $\tau$ using \eqref{eq:tau_Houbolt} and \eqref{eq:tau0_Lie}. The starting points are the same for different methods and randomly chosen on unit sphere centered at $\oRn$. The stopping tolerances $\tolf=10^{-4}$ and $\tolu=10^{-2}$. Note that for large-scale dataset, we suggest using a big damping ratio $\gamma/m$ for faster convergence, a ratio between $100$ and $1000$ often works fine.

\begin{table}[htb!]
	\begin{center}
		\caption{Numerical results on large-scale dataset with $\varepsilon=10^{-6}$, $m=1$, $\gamma=300$, $c=0$, $\tau = 1.41\times 10^{-3}$ for the Houbolt scheme, $\tau=10^{-6}$ for the Lie scheme, $t\in [0,0.3]$ for the RK(4,5) scheme, \texttt{Tolf}=$10^{-4}$, \texttt{TolU}=$10^{-2}$, and with default parameters for IPOPT and GUROBI solvers.}
		\label{tab:results5}
		\resizebox{\linewidth}{!}{%
			\tabcolsep=2pt
			\begin{tabular}{c|c|cccc|cccc|cccc|cccc|cc} \hline
				\multirow{2}{*}{n} & \multirow{2}{*}{d} & \multicolumn{4}{c}{Houbolt} & \multicolumn{4}{|c}{Lie} & \multicolumn{4}{|c}{RK(4,5)} & \multicolumn{4}{|c}{IPOPT} & \multicolumn{2}{|c}{QB-G}\\
				\cline{3-20}
				& & obj & iter & time & $\delta$ & obj & iter & time & $\delta$ & obj & iter & time & $\delta$ & obj & iter & time & $\delta$ & obj & time \\
				\hline\hline
				$10$ & $5$ & $113$ & $17$ & $0.01$ & $4.92e-03$ &$113$ & $7$ & $0.07$ & $3.08e-03$ &$113$ & $289$ & $0.17$ & $1.00e-02$ &$113$ & $16$ & $0.05$ & $3.11e-04$ &$-115$ & $0.17$ \\
				$10$ & $6$ & $237$ & $17$ & $0.02$ & $4.93e-03$ &$237$ & $7$ & $0.18$ & $3.20e-03$ &$237$ & $289$ & $0.57$ & $1.00e-02$ &$237$ & $17$ & $0.12$ & $4.94e-04$ &$-77$ & $0.31$ \\
				$12$ & $5$ & $-738$ & $17$ & $0.01$ & $6.15e-03$ &$-738$ & $7$ & $0.13$ & $3.46e-03$ &$-738$ & $313$ & $0.37$ & $1.00e-02$ &$-738$ & $16$ & $0.09$ & $5.40e-04$ &$472$ & $0.53$ \\
				$12$ & $6$ & $-264$ & $17$ & $0.05$ & $6.01e-03$ &$-264$ & $7$ & $0.49$ & $3.13e-03$ &$-264$ & $1297$ & $5.78$ & $7.90e-04$ &$-264$ & $16$ & $0.28$ & $7.68e-04$ &$178$ & $1.86$ \\
				$14$ & $5$ & $656$ & $17$ & $0.03$ & $6.01e-03$ &$656$ & $8$ & $0.30$ & $2.17e-03$ &$656$ & $329$ & $0.70$ & $1.00e-02$ &$656$ & $20$ & $0.19$ & $5.67e-04$ &$508$ & $1.79$ \\
				$14$ & $6$ & $-1223$ & $17$ & $0.11$ & $6.91e-03$ &$-1223$ & $7$ & $1.18$ & $3.73e-03$ &$-1223$ & $1313$ & $11.81$ & $1.30e-03$ &$-1223$ & $17$ & $0.66$ & $9.63e-04$ &$-645$ & $10.51$ \\
				$16$ & $5$ & $284$ & $17$ & $0.04$ & $5.86e-03$ &$284$ & $7$ & $0.49$ & $5.39e-03$ &$284$ & $313$ & $1.14$ & $1.00e-02$ &$284$ & $18$ & $0.39$ & $8.55e-04$ &$1492$ & $7.37$ \\
				$16$ & $6$ & $2416$ & $17$ & $0.23$ & $5.78e-03$ &$2416$ & $7$ & $2.96$ & $5.85e-03$ &$2416$ & $1301$ & $25.81$ & $2.21e-03$ &$2416$ & $14$ & $1.80$ & $2.21e-03$ &$-574$ & $71.65$ \\
				$18$ & $5$ & $156$ & $17$ & $0.08$ & $7.92e-03$ &$156$ & $7$ & $1.00$ & $4.83e-03$ &$156$ & $1317$ & $8.00$ & $6.21e-04$ &$156$ & $14$ & $0.92$ & $6.21e-04$ &$-122$ & $26.07$ \\
				$18$ & $6$ & $-713$ & $17$ & $0.48$ & $8.88e-03$ &$-713$ & $8$ & $8.01$ & $2.01e-03$ &$-713$ & $1317$ & $54.46$ & $1.96e-03$ &$-713$ & $15$ & $4.38$ & $1.87e-03$ &$391$ & $340.38$ \\
				$20$ & $5$ & $-1023$ & $17$ & $0.16$ & $8.84e-03$ &$-1023$ & $8$ & $2.06$ & $1.66e-03$ &$-1023$ & $317$ & $3.78$ & $1.00e-02$ &$499$ & $17$ & $2.12$ & $1.62e-03$ &$-131$ & $78.83$ \\
				$20$ & $6$ & $-2379$ & $17$ & $0.83$ & $9.65e-03$ &$-2379$ & $8$ & $14.84$ & $2.15e-03$ &$-2379$ & $1325$ & $93.82$ & $2.32e-03$ &$-2379$ & $16$ & $8.54$ & $1.99e-03$ &$1881$ & $1028.46$ \\
				\hline
				\multicolumn{2}{c|}{average} & & $17$ & $0.17$ & $6.82e-03$ & & $7$ & $2.64$ & $3.39e-03$ & & $810$ & $17.20$ & $5.77e-03$ & & $16$ & $1.63$ & $1.07e-03$ & & $130.66$ \\
				\hline
		\end{tabular}}
	\end{center}	
\end{table}

\paragraph{\textbf{Observations}}
We can observe in Table \ref{tab:results5} that the three ODE methods and IPOPT solver converge much faster than the QB-G method on large-scale dataset. Moreover, the ODE methods and IPOPT solver quite often provide same computed solutions (with average $\delta$ of order $O(10^{-3})$), which seem frequently better than the solutions of the QB-G method ($7$ better cases over $12$ are observed in Table \ref{tab:results5}). The Houbolt scheme is the fastest solver ($0.17$ seconds in average) among the others, then IPOPT ($1.63$ seconds in average) and the Lie scheme ($2.64$ seconds in average) are slightly slower than the Houbolt scheme, the RK(4,5) scheme seems a little bit more slower ($17.2$ seconds in average), and the QB-G method is the slowest one ($130.66$ seconds in average). Note again that our ODE approaches are codes in MATLAB and IPOPT and Gurobi are coded in C++ and Fortran, so that the computation time for our ODE methods should be reduced again in C++ implementation.

\begin{remark}\label{rmk:7.1}
	It is worth noting that: (i) using zero initial point will often slow down the ODE approaches, particularly to the Lie scheme, the average number of iterations could increase to more than $600$, which is reduced to $8$ with non-zero initial point randomly chosen on unit sphere centered at $\oRn$. The reason seems that zero initial point is too far from the optimum, therefore, we strongly suggest not using zero initial point for ODE schemes, particularly for large-scale cases. (ii) the Lie scheme requires solving two nonlinear systems, \cref{eq:prob_uk0.5_Lie} and \cref{eq:prob_uk+1_Lie}, to obtain $U_{k+1}$ from $U_k$, where the system \cref{eq:prob_uk+1_Lie} is solved explicitly using Cardano's formula (inexpensive) and the system \cref{eq:prob_uk0.5_Lie} is solved using a nonlinear equation solver (expensive in general for large-scale cases since the evaluation of the gradient $\nabla \Pi$ is involving in each iteration of the Netwton's method); While the Houbolt scheme needs only one system \cref{eq:prob_uk+1} which can be solved explicitly using Cardano's formula again (inexpensive). So that if there are large number of iterations in the Lie scheme with fixed $\tau$, then varying $\tau$ using formulation \eqref{eq:updatetau} is suggested to reduce the number of iterations and thus improve its numerical performance.
\end{remark}

\subsection{Optimality of The Computed Solutions}
The penalty problems \cref{prob:approx_int_cv_opt} and \cref{prob:penaltyoverball} are in general nonconvex, and our proposed methods are first-order optimization approaches which can only find a stationary point (local minimizer at best). The quality of the computed solutions depends on the initial conditions, i.e., $\U_0$ for the Lie scheme, and $(\U_0,\V_0)$ for the Houbolt and RK(4,5) schemes.

In order to justify the quality of the computed solutions, we will compare the results with the exact global optimal solutions provided by an exhaustive method (checking all points in $\{\pm 1\}^n$) for solving problem \cref{prob:int_opt}. This method only works for not too large size problems. 
Keep in mind that the total number of feasible solutions for \cref{prob:int_opt} is $2^n$. Due to the complicated and dense structure of general polynomial function $\Pi(\V)$ (with $n$ variables and degree $d$) which has at most $\binom{n+d}{d}$ monomials, the evaluation of $J_{\varepsilon}(\V)$ and its gradient $\nabla J_{\varepsilon}(\V)$ could be very time consuming. In POLYLAB, the evaluation procedure is coded in C language and called in MATLAB through \verb|mex| function for better speed. Note that the polynomial evaluations will be significantly slower in pure MATLAB routine than in \verb|mex| function. Moreover, all polynomial coefficients and exponents are designed as sparse matrices in POLYLAB which aims at handling more efficiently large-scale problems with sparse structure.   

In order to deal with large-scale instances more effectively, we use parallel computing techniques in exhaustive method, and tested on our cluster using $80$ CPUs. Moreover, we also designed a parallel framework for our ODE methods (Lie, Houbolt and RK(4,5) schemes) and IPOPT in the way that: each problem is solved simultaneously $80$ times (one task per CPU) from randomly chosen initial conditions on the unit sphere centered at $\oRn$, then we compare the obtained solutions and return the best one possessing the smallest objective value. Taking advantage of high performance devices, we are able to quickly update local solutions with random multi-starts.

Table \ref{tab:results3} illustrates some parallel computing numerical results on both small and large-scale datasets for ODE methods and IPOPT solver comparing to the exact global optimal solutions provided by the exhaustive method. The average number of iterations for $80$ runs \texttt{avgiter}, and the total wall-clock computation time \texttt{tt} for solving $80$ problems from different initial points using an SPMD (Single Program Multiple Data) parallel scheme on MATLAB are summarized. The optimality gaps to the global optimal solutions is measured by the errors (namely \texttt{err}) between the smallest value of $\Pi$ at the rounding solutions (\text{round}$(\Ueps)$) provided by our methods and at the exact global optimal solutions ($\U^*$) obtained by the exhaustive method via the formulation:
$$\texttt{err} = \frac{|\Pi(\text{round}(\U_{\varepsilon})) - \Pi(\U^*)|}{1+|\Pi(\U^*)|}.$$
This error represents a compromise between the relative error and the absolute error, which could be greater than $1$. Clearly, if $|\Pi(\U^*)|=0$, then \texttt{err} equals to the absolute error; otherwise, \texttt{err} approximates the relative error. Note that we introduce $1$ in the denominator to assume that the error is well defined.

\begin{table}[htb!]
	\begin{center}
		\caption{Parallel results on large and small scale datasets with $80$ CPUs, $\varepsilon=10^{-5}$, $m=1$, $\gamma=300$, $c=0$, $\tau = 4.5\times 10^{-3}$ for the Houbolt scheme, $\tau=10^{-5}$ for the Lie scheme, $t\in [0,0.3]$ for the RK(4,5) scheme, \texttt{Tolf}=$10^{-4}$, \texttt{TolU}=$10^{-2}$, with default parameters for IPOPT and GUROBI solvers, and random initial points on unit sphere centered at $\oRn$.}
		\label{tab:results3}
		\resizebox{\linewidth}{!}{%
			\tabcolsep=2pt
			\begin{tabular}{c|c|ccccc|ccccc|ccccc|ccccc|cc} \hline
				\multirow{2}{*}{n} & \multirow{2}{*}{d} & \multicolumn{5}{c}{Houbolt} & \multicolumn{5}{|c}{Lie} & \multicolumn{5}{|c}{RK(4,5)} & \multicolumn{5}{|c}{IPOPT} & \multicolumn{2}{|c}{Exhaustive}\\
				\cline{3-24}
				& & obj & avgiter & tt & $\delta$ & err & obj & avgiter & tt & $\delta$ & err & obj & avgiter & tt & $\delta$ & err & obj & avgiter & tt & $\delta$ & err & obj & time \\
				\hline\hline
				$2$ & $2$ & $-19$ & $12$ & $0.23$ & $6.40e-03$ & $0.00$ &$-19$ & $6$ & $0.18$ & $4.30e-03$ & $0.00$ &$-19$ & $122$ & $0.21$ & $1.00e-02$ & $0.00$ &$-19$ & $17$ & $0.28$ & $3.54e-05$ & $0.00$ &$-19$ & $0.06$ \\
				$2$ & $3$ & $21$ & $11$ & $0.14$ & $5.46e-03$ & $0.00$ &$21$ & $6$ & $0.19$ & $2.22e-03$ & $0.00$ &$21$ & $105$ & $0.18$ & $1.00e-02$ & $0.00$ &$21$ & $16$ & $0.21$ & $1.55e-04$ & $0.00$ &$21$ & $0.06$ \\
				$2$ & $4$ & $-6$ & $11$ & $0.16$ & $4.78e-03$ & $0.00$ &$-6$ & $6$ & $0.18$ & $1.78e-03$ & $0.00$ &$-6$ & $133$ & $0.18$ & $6.52e-05$ & $0.00$ &$-6$ & $16$ & $0.24$ & $6.52e-05$ & $0.00$ &$-6$ & $0.11$ \\
				$4$ & $2$ & $-20$ & $13$ & $0.18$ & $4.13e-03$ & $0.00$ &$-20$ & $6$ & $0.17$ & $1.79e-03$ & $0.00$ &$-20$ & $247$ & $0.20$ & $1.00e-02$ & $0.00$ &$-20$ & $16$ & $0.22$ & $5.72e-05$ & $0.00$ &$-20$ & $0.06$ \\
				$4$ & $3$ & $-35$ & $13$ & $0.16$ & $6.57e-03$ & $0.00$ &$-35$ & $6$ & $0.16$ & $1.41e-03$ & $0.00$ &$-35$ & $260$ & $0.19$ & $1.00e-02$ & $0.00$ &$-35$ & $16$ & $0.25$ & $1.48e-04$ & $0.00$ &$-35$ & $0.07$ \\
				$4$ & $4$ & $-49$ & $13$ & $0.17$ & $5.57e-03$ & $0.00$ &$-49$ & $6$ & $0.17$ & $2.92e-03$ & $0.00$ &$-49$ & $255$ & $0.19$ & $1.00e-02$ & $0.00$ &$-49$ & $16$ & $0.22$ & $3.72e-04$ & $0.00$ &$-49$ & $0.06$ \\
				$6$ & $2$ & $-35$ & $13$ & $0.17$ & $5.32e-03$ & $0.00$ &$-31$ & $7$ & $0.16$ & $2.95e-03$ & $0.11$ &$-35$ & $265$ & $0.21$ & $1.00e-02$ & $0.00$ &$-33$ & $18$ & $0.24$ & $5.83e-05$ & $0.06$ &$-35$ & $0.06$ \\
				$6$ & $3$ & $-78$ & $14$ & $0.16$ & $5.04e-03$ & $0.17$ &$-94$ & $7$ & $0.17$ & $1.53e-03$ & $0.00$ &$-78$ & $267$ & $0.23$ & $1.00e-02$ & $0.17$ &$-94$ & $16$ & $0.24$ & $3.05e-04$ & $0.00$ &$-94$ & $0.06$ \\
				$6$ & $4$ & $-154$ & $14$ & $0.16$ & $4.31e-03$ & $0.06$ &$-164$ & $7$ & $0.17$ & $1.35e-03$ & $0.00$ &$-154$ & $267$ & $0.22$ & $1.00e-02$ & $0.06$ &$-154$ & $16$ & $0.22$ & $8.27e-04$ & $0.06$ &$-164$ & $0.06$ \\
				$8$ & $2$ & $-89$ & $14$ & $0.17$ & $7.51e-03$ & $0.00$ &$-85$ & $7$ & $0.18$ & $3.88e-03$ & $0.04$ &$-89$ & $290$ & $0.30$ & $3.16e-04$ & $0.00$ &$-75$ & $17$ & $0.22$ & $8.12e-05$ & $0.16$ &$-89$ & $0.09$ \\
				$8$ & $3$ & $-122$ & $14$ & $0.17$ & $8.79e-03$ & $0.02$ &$-124$ & $7$ & $0.19$ & $2.30e-03$ & $0.00$ &$-122$ & $283$ & $0.23$ & $1.00e-02$ & $0.02$ &$-106$ & $17$ & $0.22$ & $6.28e-04$ & $0.14$ &$-124$ & $0.06$ \\
				$8$ & $4$ & $-217$ & $14$ & $0.19$ & $5.28e-03$ & $0.23$ &$-281$ & $7$ & $0.19$ & $1.17e-03$ & $0.00$ &$-189$ & $311$ & $0.24$ & $1.00e-02$ & $0.33$ &$-191$ & $16$ & $0.24$ & $4.63e-04$ & $0.32$ &$-281$ & $0.06$ \\
				$10$ & $2$ & $-67$ & $14$ & $3.15$ & $5.43e-03$ & $0.29$ &$-73$ & $7$ & $2.80$ & $2.77e-03$ & $0.23$ &$-75$ & $266$ & $0.55$ & $1.00e-02$ & $0.21$ &$-67$ & $18$ & $0.61$ & $1.70e-04$ & $0.29$ &$-95$ & $0.40$ \\
				$10$ & $3$ & $-169$ & $14$ & $0.24$ & $2.06e-02$ & $0.08$ &$-103$ & $7$ & $0.30$ & $3.36e-03$ & $0.43$ &$-183$ & $297$ & $0.33$ & $9.36e-04$ & $0.00$ &$-155$ & $17$ & $0.27$ & $6.56e-04$ & $0.15$ &$-183$ & $0.11$ \\
				$10$ & $4$ & $-331$ & $14$ & $0.20$ & $4.35e-03$ & $0.26$ &$-363$ & $7$ & $0.24$ & $2.61e-03$ & $0.19$ &$-341$ & $316$ & $0.36$ & $2.59e-03$ & $0.24$ &$-351$ & $17$ & $0.28$ & $1.25e-03$ & $0.22$ &$-449$ & $0.11$ \\
				$10$ & $5$ & $-533$ & $14$ & $0.30$ & $8.66e-03$ & $0.10$ &$-521$ & $7$ & $0.42$ & $3.53e-03$ & $0.12$ &$-423$ & $314$ & $0.54$ & $3.46e-03$ & $0.29$ &$-493$ & $17$ & $0.33$ & $1.66e-03$ & $0.17$ &$-595$ & $0.09$ \\
				$10$ & $6$ & $-821$ & $14$ & $0.32$ & $1.82e-02$ & $0.23$ &$-1033$ & $7$ & $0.63$ & $8.19e-03$ & $0.03$ &$-1061$ & $323$ & $1.34$ & $1.09e-02$ & $0.00$ &$-1061$ & $17$ & $0.43$ & $5.27e-03$ & $0.00$ &$-1061$ & $0.17$ \\
				$12$ & $5$ & $-586$ & $15$ & $0.27$ & $8.65e-03$ & $0.47$ &$-800$ & $7$ & $0.50$ & $4.48e-03$ & $0.27$ &$-894$ & $319$ & $0.84$ & $1.00e-02$ & $0.19$ &$-756$ & $17$ & $0.38$ & $3.88e-03$ & $0.31$ &$-1102$ & $0.19$ \\
				$12$ & $6$ & $-1506$ & $15$ & $0.59$ & $1.51e-02$ & $0.37$ &$-1632$ & $7$ & $1.45$ & $1.29e-02$ & $0.32$ &$-1340$ & $302$ & $2.49$ & $1.42e-02$ & $0.44$ &$-1420$ & $17$ & $0.65$ & $5.36e-03$ & $0.41$ &$-2402$ & $0.52$ \\
				$14$ & $5$ & $-1204$ & $15$ & $0.39$ & $1.48e-02$ & $0.42$ &$-1350$ & $7$ & $0.76$ & $8.84e-03$ & $0.35$ &$-1288$ & $334$ & $1.58$ & $1.00e-02$ & $0.38$ &$-1456$ & $18$ & $0.49$ & $4.11e-03$ & $0.30$ &$-2088$ & $0.99$ \\
				$14$ & $6$ & $-1669$ & $15$ & $1.16$ & $1.78e-02$ & $0.48$ &$-2093$ & $7$ & $2.89$ & $1.69e-02$ & $0.35$ &$-2371$ & $348$ & $5.36$ & $1.95e-02$ & $0.27$ &$-1859$ & $18$ & $1.49$ & $8.34e-03$ & $0.43$ &$-3235$ & $3.54$ \\
				$16$ & $5$ & $-1354$ & $15$ & $0.58$ & $1.22e-02$ & $0.49$ &$-1348$ & $7$ & $1.29$ & $8.11e-03$ & $0.49$ &$-1034$ & $336$ & $2.44$ & $1.00e-02$ & $0.61$ &$-1216$ & $17$ & $0.86$ & $4.64e-03$ & $0.54$ &$-2632$ & $6.02$ \\
				$16$ & $6$ & $-2738$ & $15$ & $2.61$ & $2.26e-02$ & $0.50$ &$-3652$ & $7$ & $7.21$ & $3.06e-02$ & $0.33$ &$-2846$ & $410$ & $12.05$ & $2.83e-02$ & $0.48$ &$-3784$ & $17$ & $2.73$ & $1.26e-02$ & $0.30$ &$-5432$ & $25.22$ \\
				$18$ & $5$ & $-1906$ & $15$ & $0.91$ & $1.31e-02$ & $0.46$ &$-1538$ & $8$ & $2.51$ & $8.66e-03$ & $0.57$ &$-2206$ & $316$ & $3.60$ & $1.60e-02$ & $0.38$ &$-2048$ & $18$ & $1.32$ & $1.03e-02$ & $0.42$ &$-3558$ & $37.12$ \\
				$18$ & $6$ & $-5963$ & $15$ & $6.04$ & $5.25e-02$ & $0.20$ &$-3929$ & $8$ & $15.48$ & $3.52e-02$ & $0.47$ &$-3981$ & $415$ & $23.68$ & $3.41e-02$ & $0.47$ &$-4133$ & $19$ & $6.08$ & $1.92e-02$ & $0.44$ &$-7443$ & $169.15$ \\
				$20$ & $5$ & $-2495$ & $15$ & $1.54$ & $1.93e-02$ & $0.59$ &$-3005$ & $8$ & $4.07$ & $1.78e-02$ & $0.50$ &$-2597$ & $360$ & $6.04$ & $1.86e-02$ & $0.57$ &$-2185$ & $19$ & $1.83$ & $1.13e-02$ & $0.64$ &$-6043$ & $214.77$ \\
				$20$ & $6$ & $-3995$ & $15$ & $12.85$ & $4.10e-02$ & $0.58$ &$-5383$ & $8$ & $31.89$ & $5.24e-02$ & $0.44$ &$-4773$ & $417$ & $40.30$ & $4.73e-02$ & $0.50$ &$-4241$ & $19$ & $13.01$ & $2.95e-02$ & $0.56$ &$-9545$ & $1082.85$ \\
				\hline
				\multicolumn{2}{c|}{average} & & $14$ & $1.23$ & $1.27e-02$ & $0.22$ & & $7$ & $2.76$ & $9.04e-03$ & $0.19$ & & $292$ & $3.85$ & $1.25e-02$ & $0.21$ & & $17$ & $1.24$ & $4.50e-03$ & $0.22$ & & $57.11$ \\
				\hline
		\end{tabular}}
	\end{center}	
\end{table}

We can observe in Table \ref{tab:results3} that for $n\leq12$, exhaustive method performs very fast (among the fastest one), since we only need to check simultaneously the values of $\Pi$ on the set $\{\pm 1\}^n$ which consists of $2^{12}$, i.e., $4096$ integer points at most, that is not a difficult task. For $n>12$, the computation time for exhaustive method increases dramatically. However, it is interesting to observe that our ODE methods and IPOPT seem to be quite stable in number of iterations regardless the increase of problem size. Note that when $n>d$, the increase of $d$ is more sensitive to the computation time than the increase of $n$. Because if $n$ is increased by $k$, then the number of monomials is increased up to $\binom{(n+k)+d}{d}$, but if $d$ is increased by $k$, then the number of monomials is $\binom{n+(d+k)}{d+k} > \binom{(n+k)+d}{d}, \forall k\in \N^*$ and $n>d$. This can explains why the case $(n,d)=(18,6)$ is more computational expensive than the case $(n,d)=(20,5)$. Note that sometimes the proposed ODE methods and IPOPT cannot get an optimal solution. Because these methods are local approaches based on the first order local optimality condition, which is necessary but not sufficient local optimality condition if the problem \eqref{prob:approx_int_cv_opt} is nonconvex.

Concerning the optimality gap, the errors for all methods (Houbolt, Lie, RK(4,5) and IPOPT) are upper bounded of order $O(10^{-1})$ (less than $0.22$) in average and particularly equal to zero for many small-scale cases and upto $0.59$ for some large-scale instances. The errors seem to increase when the size of problems increases. In general, there are very slightly difference in average \texttt{err} between ODE approaches and IPOPT solver, which confirms the fact that our ODE methods could provide local optimal solutions as good as IPOPT.

\subsection{Tests on MQLib Large-scale Benchmark Dataset}\label{subsec:testonMQLib}
In this subsection, we are going to present the numerical results on a large-scale benchmark dataset \texttt{MQLib} using our ODE methods and compare with IPOPT. The Max-Cut problem (in form of formulation \eqref{prob:int_opt}) is defined by:
	$$\max_{\V\in \{\pm 1\}^n} \frac{1}{2} \sum_{1\leq i<j\leq n} w_{ij} (1-v_i v_j),$$
	where $W=(w_{ij}) \in \R^{n\times n}$ is given weighted adjacency matrix. This problem is well-known as one of the Karp's $21$ NP-complete problems \cite{Karp1972}. The QUBO problem (in form of formulation \eqref{prob:Boolean_opt}) is defined similarly as:
	$$\max_{\Y\in \{0,1\}^n} \Y^{\top} S \Y$$
	with given matrix $S\in \R^{n\times n}$. As we shown in Section \ref{sec:formulations} that these two problems are indeed equivalent and particular cases of \eqref{prob:Boolean_opt} and \eqref{prob:int_opt} with quadratic functional $\Pi$.

For the particularity of the quadratic functional $\Pi$ where its gradient is known explicitly, we can use matrix form instead of POLYLAB to model these instances for better numerical performance in the computation of the values of $\Pi$ and its gradient. In this case, the Lie scheme is simplified as an explicit scheme.
\paragraph{\textbf{Simplification of the Lie scheme for Max-Cut and QUBO problems}} Let us denote the general form of the Max-Cut and QUBO problems as:
	\begin{equation}
	\label{prob:QP}
	\min_{\V\in \{\pm 1\}^n} \Pi(\V) = \frac{1}{2}\V^{\top}Q\V + \vl^{\top}\V + s,
	\end{equation}
	where $Q\in \R^{n\times n}, \vl\in \R^n, s\in \R$. Particularly, for Max-Cut problem, we have $$Q = \frac{W}{2}, ~\vl = \oRn, ~s=-\frac{1}{2}\sum_{1\leq i<j\leq n} w_{ij};$$ and for QUBO problem, we have $$Q = -\frac{S}{2}, ~ \vl = -\frac{1}{2}\one^{\top}S, ~ s=-\frac{1}{4}\one^{\top}S\one.$$ In the quadratic formulation \eqref{prob:QP}, the gradient $\nabla \Pi(\V)$ is explicitly computed by:
	$$\nabla \Pi(\V) = Q\V + \vl,$$
	then the system \eqref{eq:system_uk0.5_lie} in the Lie scheme is simplified as:
	$$\U^{k+\frac{1}{2}} + \tau \nabla \Pi \left(\U^{k+\frac{1}{2}}\right) = \U^{k+\frac{1}{2}} + \tau Q\U^{k+\frac{1}{2}} + \tau\vl = \U^k,$$
	and $\U^{k+\frac{1}{2}}$ is obtained by solving the linear system:
	\begin{equation}
	\label{eq:linsys}
	(I_n + \tau Q)\U^{k+\frac{1}{2}} = \U^k - \tau\vl,
	\end{equation}
	where $I_n$ denotes the identity matrix of size $n\times n$. We may use a fast linear system solver, e.g., the MATLAB functions \texttt{decomposition} (for creating reusable matrix decompositions: LU, LDL, Cholesky, QR and more), \texttt{mldivide} (for sparse matrix), or \texttt{linsolve} (for dense matrix) for efficiently solving this linear system with real symmetric matrix $I_n+\tau Q$. Note that the matrix $I_n+\tau Q$ is never changed during the iterations for methods with fixed $\tau$, thus creating a reusable matrix through \texttt{decomposition} is much faster than the other methods since it avoids decomposition of the same matrix $I_n + \tau Q$ multiple times in different iterations. Finally, we can see that the Lie scheme is simplified as an explicit scheme for Max-Cut and QUBO problems.

\begin{figure*}[h!]
	\centering
	\subfigure[GKA with $44$ cases $21 \leq n \leq 501$]{ 
		\includegraphics[width=0.45\textwidth]{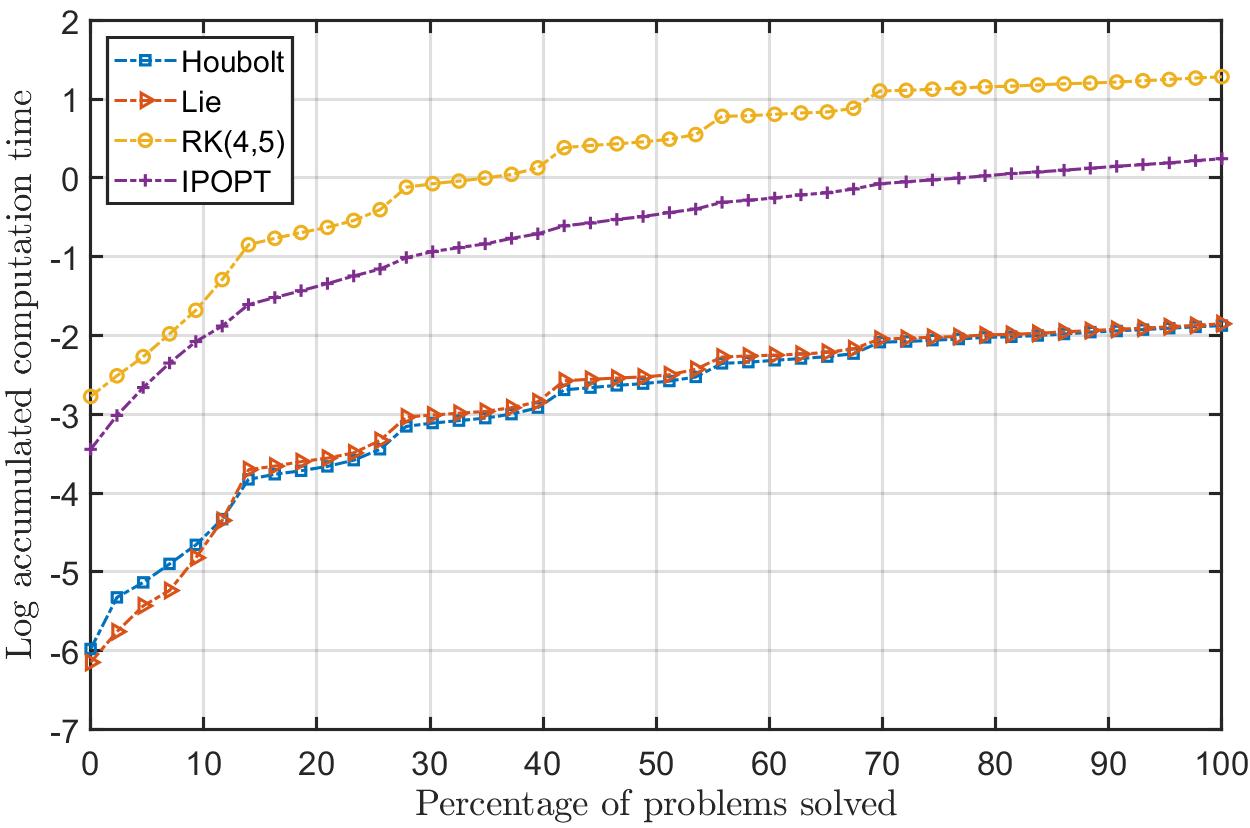}
	}
	\subfigure[G-set with $71$ cases $800 \leq n \leq 20000$]{ 
		\includegraphics[width=0.45\textwidth]{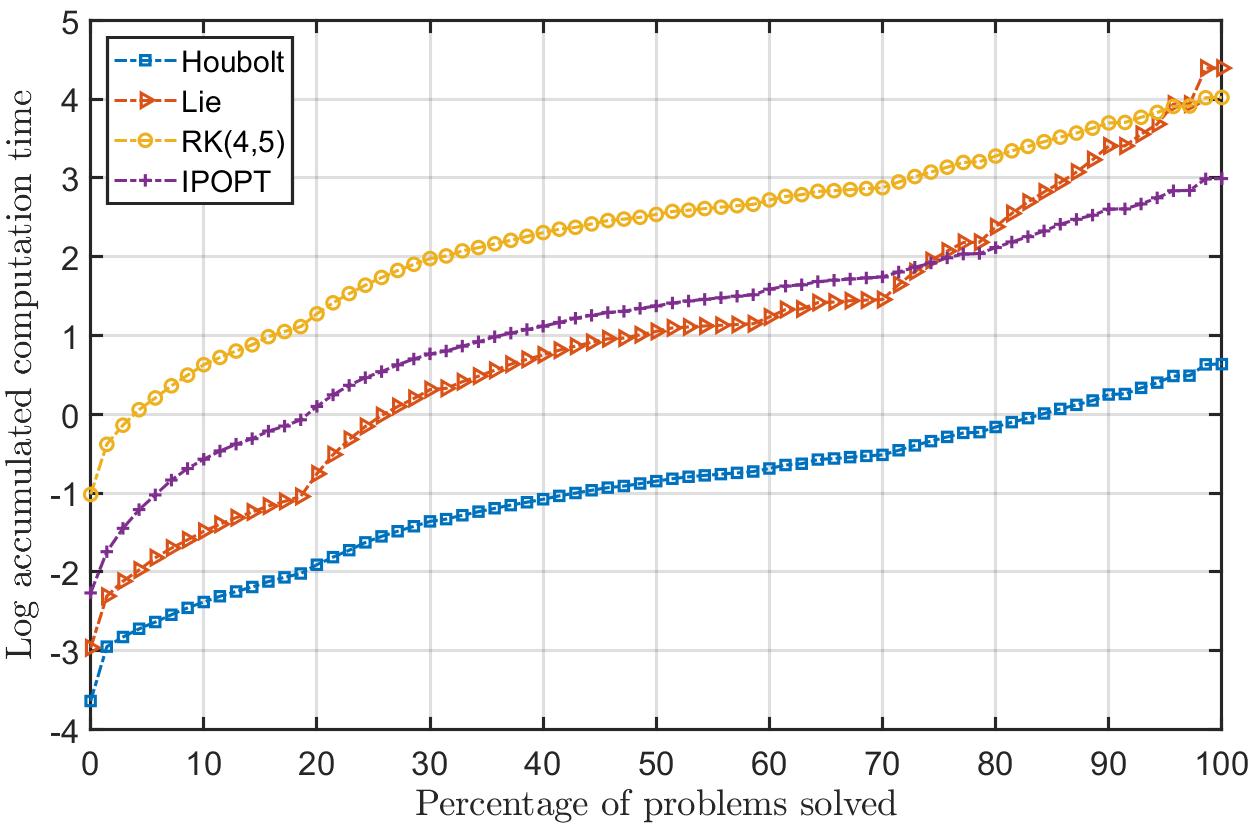}
	}
	\subfigure[Beasley with $59$ cases $51 \leq n \leq 2501$]{ 
		\includegraphics[width=0.45\textwidth]{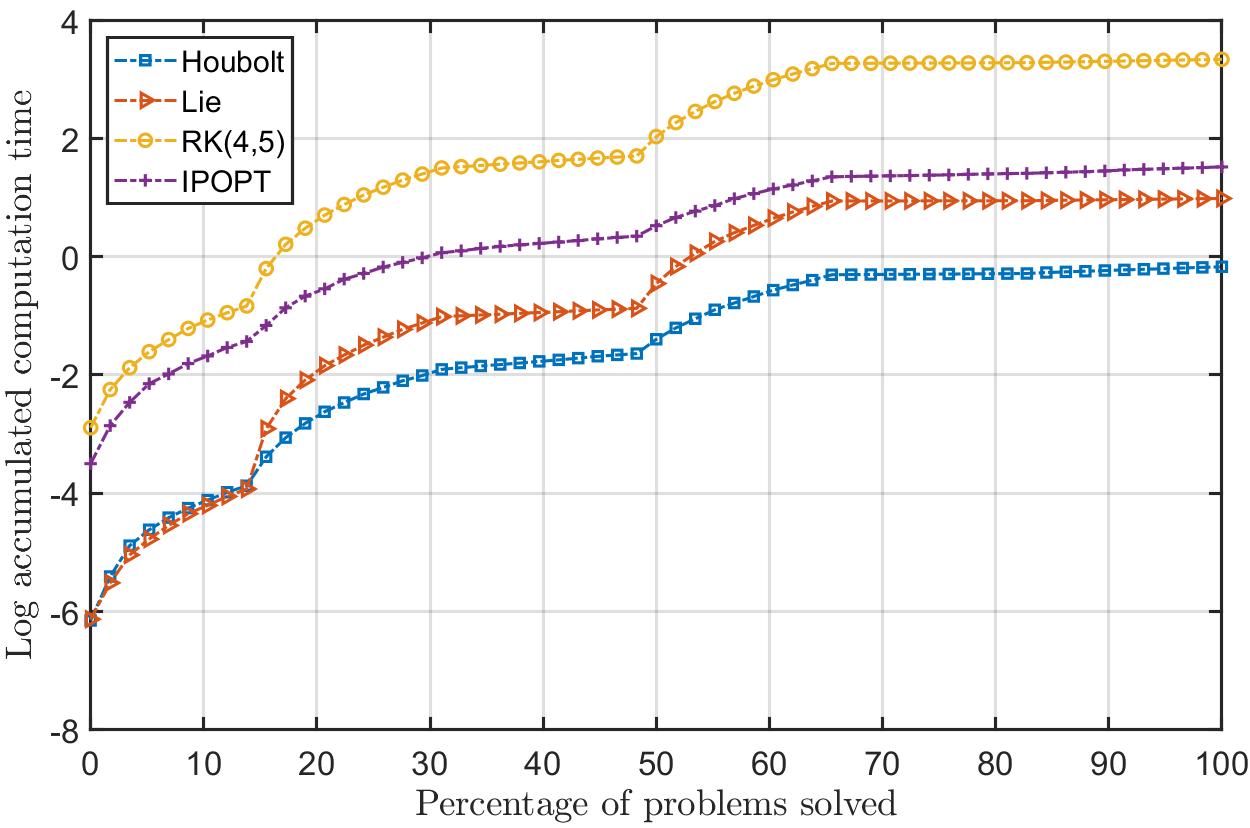} 
	}
	\subfigure[Culberson with $107$ cases $1015 \leq n \leq 4997$]{ 
		\includegraphics[width=0.45\textwidth]{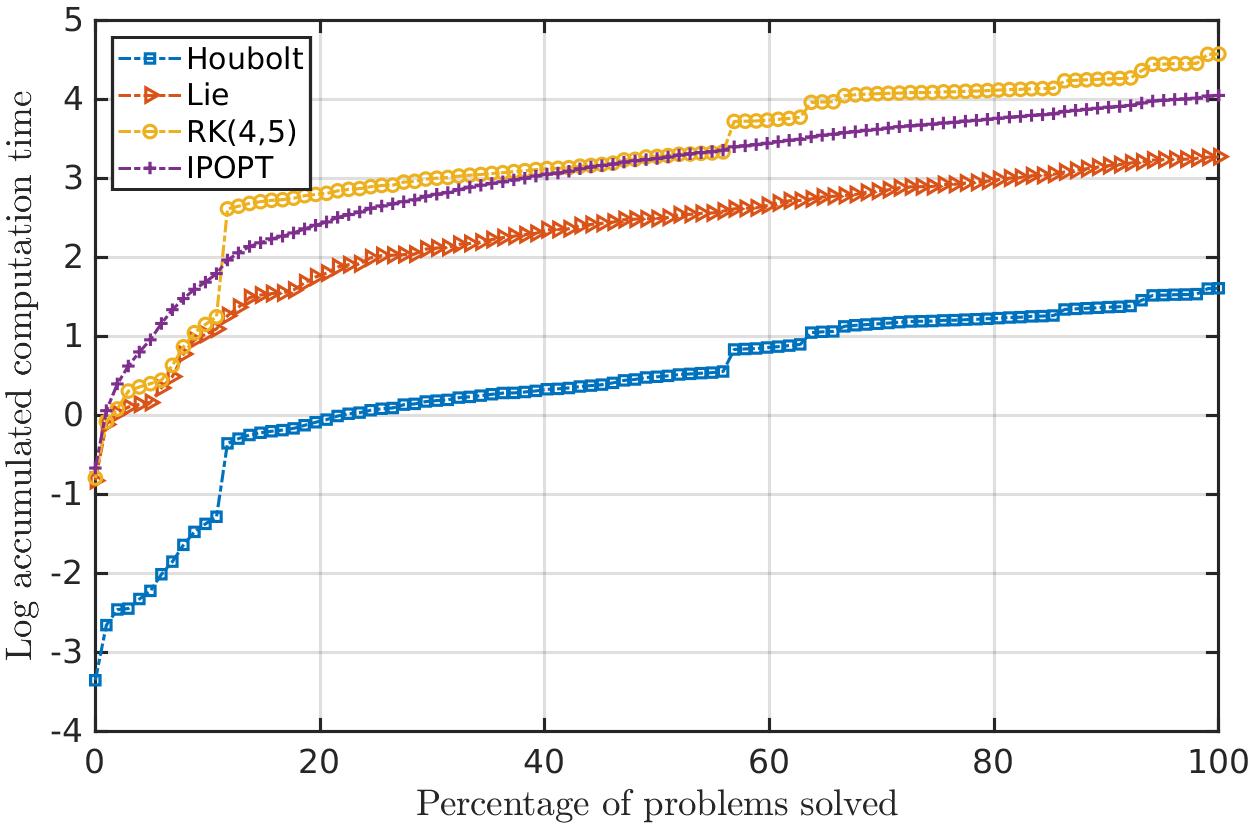} 
	}
	\caption{Percentage of problems solved v.s. Log accumulated computation time for the Houbolt, Lie, RK(4,5) schemes and IPOPT solver on $4$ subsets of MQLib (tests on single CPU with parameters $\varepsilon=10^{-5}$, $m=1$, $\gamma=300$, $c=0$, $\tau = 4.5\times 10^{-3}$ for the Houbolt scheme, $\tau=10^{-5}$ for the Lie scheme, $t\in [0,0.3]$ for the RK(4,5) scheme, \texttt{Tolf}=$10^{-4}$, \texttt{TolU}=$10^{-2}$, with default parameters for IPOPT solver, and random initial points on unit sphere centered at $\oRn$).}
	\label{fig:MQLib_computingtime}
\end{figure*}

\paragraph{\textbf{Numerical results on MQLib dataset}} Now, we tested on $4$ subsets of MaxCut problems (namely, GKA, G-set, Beasley and Culberson) in the MQLib dataset where GKA and Beasley are QUBO instances but converted to Max-Cut. These subsets consists of $281$ instances with variables from $21$ to $20000$. Figure \ref{fig:MQLib_computingtime} illustrates the comparative results on the computation time for ODE methods and IPOPT. 
We can observe in Figure \ref{fig:MQLib_computingtime} that the fastest method is always the Houbolt scheme, then the Lie scheme is in the second rank quite often faster than IPOPT solver. The slowest one is almost always the RK(4,5) scheme. One exception is observed on the G-set where the Lie scheme slowed down at the end, which is due to the fact that there are very large-scale instances with $20000$ variables leading to the solution of the linear system \eqref{eq:linsys} in the Lie scheme quite slow.

Note that for problems with ill-conditioned matrices $Q$, the ODE methods may difficult to converge. Decreasing $\varepsilon$ and increasing the ratio $\gamma/m$ may improve the convergence of ODE schemes. An example of p-set in MQLib, namely \texttt{p4000\_2} with $4001$ variables where the condition number of the matrix $Q$ is $2.92\times 10^{7}$, shows that with $\varepsilon=10^{-5}$ and $\gamma/m=100$, the Houbolt scheme converges with $26$ iterations in $0.83$ seconds, but the Lie scheme do not converge and the RK(4,5) scheme converges very slowly (with $1201$ iterations in $43.35$ seconds). If we decrease $\varepsilon$ to $10^{-8}$, then the Lie scheme converges with $11$ iterations in $1.01$ seconds; If we increase also $\gamma/m$ to $10000$, then the RK(4,5) scheme converges much faster with $261$ iterations in $10.31$ seconds. Therefore, for large-scale ill-conditioned problems, we suggest setting $\varepsilon$ smaller than $10^{-7}$ and $\gamma/m$ bigger than $10^{4}$.

Based on above numerical tests, we can conclude that our ODE methods appear to be promising approaches which perform stable and fast convergence, especially in relative large-scale cases, to obtain local minimizers often close to global minimizers of the Boolean program, the relative approximation error being of order $O(10^{-1})$ in average for $\varepsilon=10^{-5}$.

\section{Conclusion and Perspective}\label{sec:conclusion}
In this paper, we investigate how to deal with Boolean polynomial program using numerical solution methods for differential equations, namely the Houbolt scheme, the Lie scheme, and the Runge-Kutta scheme. Our methods are not limited to problem with polynomials but also applicable to any Boolean optimization problem involving differentiable non-polynomial functions. Based on the equivalence between Boolean optimization and integer optimization with decision variables in $\{\pm 1\}^n$, we introduce a quartic penalty approximation for Boolean polynomial program, and prove that any converging sub-sequence of minimizers for penalty problems converge to a global minimizer of the Boolean problem. It is interesting to observe that the distance between minimizers of penalty formulations and the set $\{\pm 1\}^n$ is bounded of order $O(\sqrt{n}\varepsilon)$ with penalty parameter $\varepsilon$. Then, we introduce three numerical methods (Houbolt, Lie and RK(4,5) schemes) to find local minimizers for penalty problem. The choice of parameters involved in these schemes is discussed. Numerical simulations on both small and large-scale synthetic datasets and on the MQLib large-scale benchmark dataset of Max-Cut and QUBO problems, testing with our ODE approaches and comparing to the nonlinear optimization solver IPOPT and a quadratic binary formulation approach (QB-G), demonstrate good performance of our ODE methods, which yield stable and fast convergence to obtain numerical solutions closing to integer ones with average error of order $O(10^{-1})$ with $\varepsilon=10^{-5}$ to global minimizers. This paper show us that it is possible to apply numerical methods for differential equation to solve hard nonconvex optimization problems even discrete ones. 

Many potential future works deserve more attentions. First, we can consider modifying our approaches to adapt optimization problems involving some constraints. It should be easy for some non-empty well-qualified constrained, such as linear constraints and convex constraints under Slater condition, based on the corresponding KKT system. Since the KKT system is still a polynomial equation, thus our methods could be applied directly without any difficulty. Moreover, it is interesting to consider more general cases to allow feasible regions defined by inequalities in polynomials (i.e., semi-algebraic sets), in which constraint qualification may not be verified (i.e., KKT system is no-longer necessary optimality condition anymore). Second, it may be interesting to compare our methods with the Lasserre's method \cite{Lasserre2002} for particular structured (e.g., sparse or low-rank) Boolean polynomial problems arising in some important real-world applications. Third, it deserves more attention on how to link the stopping criteria to the error bounds in order to obtain a computed solution with desired error bound precision. Moreover, we can extend our approaches to solve other hard combinatorial or continuous optimization problems, e.g., when the integer variables are restricted to some discrete set different from Boolean ones, or when the functional $\Pi$ is not polynomial (e.g., convex/concave functions, Lipschitz functions, DC (difference-of-convex) functions etc.). We truly believe that many other numerical schemes for differential equations should be also useful to deal with hard optimization problems in a similar way once an ingenious differentiable optimization reformulation is established, and we hope that our article will share some lights to further more extensive researches in these directions.

\begin{acknowledgements}
The authors are partially supported by the National Natural Science Foundation of China (Grant 11601327) and the Hong Kong Kennedy Wong Foundation. Special thanks to the anonymous reviewers for their careful reading of our manuscript and their insightful comments and suggestions which improved the quality of our paper.
\end{acknowledgements}

%
%

\bibliographystyle{spmpsci}      
\bibliography{references}   

\end{document}